\documentclass[oneside,12pt,a4paper]{article}
\usepackage[colorlinks,citecolor=DarkRed,urlcolor=DarkRed]{hyperref} % (find-es "tex" "hyperref")
\usepackage{tocloft}                   % (find-es "tex" "tocloft")
\usepackage{indentfirst}
\usepackage{amsmath}
\usepackage{amsfonts}
\usepackage{amssymb}
\usepackage{pict2e}
\usepackage[x11names,svgnames]{xcolor} % (find-es "tex" "xcolor")
\usepackage{ifluatex}
\ifluatex
  \input edrxaccents.tex            % (find-LATEX "edrxaccents.tex")
  \input edrx21chars.tex            % (find-LATEX "edrxchars.tex")
  \input edrxheadfoot.tex           % (find-LATEX "edrxheadfoot.tex")
\else
  \usepackage[utf8]{inputenc}
\DeclareUnicodeCharacter{00A1}{\text{\textexclamdown}} % ¡
\DeclareUnicodeCharacter{00A7}{\S}                 % §
\DeclareUnicodeCharacter{00AC}{\neg}               % ¬
\DeclareUnicodeCharacter{00B0}{^\circ}             % °
\DeclareUnicodeCharacter{00B1}{\pm}                % ±
\DeclareUnicodeCharacter{00B2}{^2}                 % ²
\DeclareUnicodeCharacter{00B3}{^3}                 % ³
\DeclareUnicodeCharacter{00B7}{\cdot}              % ·
\DeclareUnicodeCharacter{00B9}{^{-1}}              % ¹
\DeclareUnicodeCharacter{00D7}{\times}             % ×
\DeclareUnicodeCharacter{00F7}{\div}               % ÷
\DeclareUnicodeCharacter{0393}{\Gamma}             % Γ
\DeclareUnicodeCharacter{0394}{\Delta}             % Δ
\DeclareUnicodeCharacter{0398}{\Theta}             % Θ
\DeclareUnicodeCharacter{039B}{\Lambda}            % Λ
\DeclareUnicodeCharacter{03A0}{\Pi}                % Π
\DeclareUnicodeCharacter{03A3}{\Sigma}             % Σ
\DeclareUnicodeCharacter{03A6}{\Phi}               % Φ
\DeclareUnicodeCharacter{03A8}{\Psi}               % Ψ
\DeclareUnicodeCharacter{03A9}{\Omega}             % Ω
\DeclareUnicodeCharacter{03B1}{\alpha}             % α
\DeclareUnicodeCharacter{03B2}{\beta}              % β
\DeclareUnicodeCharacter{03B3}{\gamma}             % γ
\DeclareUnicodeCharacter{03B4}{\delta}             % δ
\DeclareUnicodeCharacter{03B5}{\epsilon}           % ε
\DeclareUnicodeCharacter{03B6}{\zeta}              % ζ
\DeclareUnicodeCharacter{03B7}{\eta}               % η
\DeclareUnicodeCharacter{03B8}{\theta}             % θ
\DeclareUnicodeCharacter{03B9}{\iota}              % ι
\DeclareUnicodeCharacter{03BA}{\kappa}             % κ
\DeclareUnicodeCharacter{03BB}{\lambda}            % λ
\DeclareUnicodeCharacter{03BC}{\mu}                % μ
\DeclareUnicodeCharacter{03BD}{\nu}                % ν
\DeclareUnicodeCharacter{03BE}{\xi}                % ξ
\DeclareUnicodeCharacter{03C0}{\pi}                % π
\DeclareUnicodeCharacter{03C1}{\rho}               % ρ
\DeclareUnicodeCharacter{03C3}{\sigma}             % σ
\DeclareUnicodeCharacter{03C4}{\tau}               % τ
\DeclareUnicodeCharacter{03C6}{\phi}               % φ
\DeclareUnicodeCharacter{03C7}{\chi}               % χ
\DeclareUnicodeCharacter{03C8}{\psi}               % ψ
\DeclareUnicodeCharacter{03C9}{\omega}             % ω
\DeclareUnicodeCharacter{03D5}{\origphi}           % ϕ
\DeclareUnicodeCharacter{1D400}{\catA}             % 𝐀
\DeclareUnicodeCharacter{1D401}{\catB}             % 𝐁
\DeclareUnicodeCharacter{1D402}{\catC}             % 𝐂
\DeclareUnicodeCharacter{1D403}{\catD}             % 𝐃
\DeclareUnicodeCharacter{1D404}{\catE}             % 𝐄
\DeclareUnicodeCharacter{1D405}{\catF}             % 𝐅
\DeclareUnicodeCharacter{1D406}{\catG}             % 𝐆
\DeclareUnicodeCharacter{1D407}{\catH}             % 𝐇
\DeclareUnicodeCharacter{1D408}{\catI}             % 𝐈
\DeclareUnicodeCharacter{1D409}{\catJ}             % 𝐉
\DeclareUnicodeCharacter{1D40A}{\catK}             % 𝐊
\DeclareUnicodeCharacter{1D40B}{\catL}             % 𝐋
\DeclareUnicodeCharacter{1D40C}{\catM}             % 𝐌
\DeclareUnicodeCharacter{1D40D}{\catN}             % 𝐍
\DeclareUnicodeCharacter{1D40E}{\catO}             % 𝐎
\DeclareUnicodeCharacter{1D40F}{\catP}             % 𝐏
\DeclareUnicodeCharacter{1D410}{\catQ}             % 𝐐
\DeclareUnicodeCharacter{1D411}{\catR}             % 𝐑
\DeclareUnicodeCharacter{1D412}{\catS}             % 𝐒
\DeclareUnicodeCharacter{1D413}{\catT}             % 𝐓
\DeclareUnicodeCharacter{1D414}{\catU}             % 𝐔
\DeclareUnicodeCharacter{1D415}{\catV}             % 𝐕
\DeclareUnicodeCharacter{1D416}{\catW}             % 𝐖
\DeclareUnicodeCharacter{1D417}{\catX}             % 𝐗
\DeclareUnicodeCharacter{1D418}{\catY}             % 𝐘
\DeclareUnicodeCharacter{1D419}{\catZ}             % 𝐙
\DeclareUnicodeCharacter{1D41B}{\mathbf}           % 𝐛
\DeclareUnicodeCharacter{1D422}{\textsl}           % 𝐢
\DeclareUnicodeCharacter{1D42B}{\mathrm}           % 𝐫
\DeclareUnicodeCharacter{1D42C}{\mathsf}           % 𝐬
\DeclareUnicodeCharacter{1D42D}{\text}             % 𝐭
\DeclareUnicodeCharacter{1D4D0}{\mathcal{A}}       % 𝓐
\DeclareUnicodeCharacter{1D4D1}{\mathcal{B}}       % 𝓑
\DeclareUnicodeCharacter{1D4D2}{\mathcal{C}}       % 𝓒
\DeclareUnicodeCharacter{1D4D3}{\mathcal{D}}       % 𝓓
\DeclareUnicodeCharacter{1D4D4}{\mathcal{E}}       % 𝓔
\DeclareUnicodeCharacter{1D4D5}{\mathcal{F}}       % 𝓕
\DeclareUnicodeCharacter{1D4D6}{\mathcal{G}}       % 𝓖
\DeclareUnicodeCharacter{1D4D7}{\mathcal{H}}       % 𝓗
\DeclareUnicodeCharacter{1D4D8}{\mathcal{I}}       % 𝓘
\DeclareUnicodeCharacter{1D4D9}{\mathcal{J}}       % 𝓙
\DeclareUnicodeCharacter{1D4DA}{\mathcal{K}}       % 𝓚
\DeclareUnicodeCharacter{1D4DB}{\mathcal{L}}       % 𝓛
\DeclareUnicodeCharacter{1D4DC}{\mathcal{M}}       % 𝓜
\DeclareUnicodeCharacter{1D4DD}{\mathcal{N}}       % 𝓝
\DeclareUnicodeCharacter{1D4DE}{\mathcal{O}}       % 𝓞
\DeclareUnicodeCharacter{1D4DF}{\mathcal{P}}       % 𝓟
\DeclareUnicodeCharacter{1D4E0}{\mathcal{Q}}       % 𝓠
\DeclareUnicodeCharacter{1D4E1}{\mathcal{R}}       % 𝓡
\DeclareUnicodeCharacter{1D4E2}{\mathcal{S}}       % 𝓢
\DeclareUnicodeCharacter{1D4E3}{\mathcal{T}}       % 𝓣
\DeclareUnicodeCharacter{1D4E4}{\mathcal{U}}       % 𝓤
\DeclareUnicodeCharacter{1D4E5}{\mathcal{V}}       % 𝓥
\DeclareUnicodeCharacter{1D4E6}{\mathcal{W}}       % 𝓦
\DeclareUnicodeCharacter{1D4E7}{\mathcal{X}}       % 𝓧
\DeclareUnicodeCharacter{1D4E8}{\mathcal{Y}}       % 𝓨
\DeclareUnicodeCharacter{1D4E9}{\mathcal{Z}}       % 𝓩
\DeclareUnicodeCharacter{2022}{\bullet}            % •
\DeclareUnicodeCharacter{2026}{\ldots}             % …
\DeclareUnicodeCharacter{2113}{\ell}               % ℓ
\DeclareUnicodeCharacter{214B}{\bindnasrepma}      % ⅋
\DeclareUnicodeCharacter{2190}{\ot}                % ←
\DeclareUnicodeCharacter{2191}{\upto}              % ↑
\DeclareUnicodeCharacter{2192}{\to}                % →
\DeclareUnicodeCharacter{2193}{\dnto}              % ↓
\DeclareUnicodeCharacter{2194}{\bij}               % ↔
\DeclareUnicodeCharacter{2195}{\updownarrow}       % ↕
\DeclareUnicodeCharacter{2196}{\nwarrow}           % ↖
\DeclareUnicodeCharacter{2197}{\nearrow}           % ↗
\DeclareUnicodeCharacter{2198}{\searrow}           % ↘
\DeclareUnicodeCharacter{2199}{\swarrow}           % ↙
\DeclareUnicodeCharacter{21A3}{\twoheadrightarrow} % ↣
\DeclareUnicodeCharacter{21A4}{\mapsfrom}          % ↤
\DeclareUnicodeCharacter{21A6}{\mapsto}            % ↦
\DeclareUnicodeCharacter{21C0}{\rightharpoonup}    % ⇀
\DeclareUnicodeCharacter{21D0}{\Leftarrow}         % ⇐
\DeclareUnicodeCharacter{21D2}{\funto}             % ⇒
\DeclareUnicodeCharacter{21D4}{\Leftrightarrow}    % ⇔
\DeclareUnicodeCharacter{2200}{\forall}            % ∀
\DeclareUnicodeCharacter{2202}{\partial}           % ∂
\DeclareUnicodeCharacter{2203}{\exists}            % ∃
\DeclareUnicodeCharacter{2205}{\emptyset}          % ∅
\DeclareUnicodeCharacter{2207}{\nabla}             % ∇
\DeclareUnicodeCharacter{2208}{\in}                % ∈
\DeclareUnicodeCharacter{2216}{\backslash}         % ∖
\DeclareUnicodeCharacter{2218}{\circ}              % ∘
\DeclareUnicodeCharacter{221A}{\sqrt}              % √
\DeclareUnicodeCharacter{221E}{\infty}             % ∞
\DeclareUnicodeCharacter{2227}{\land}              % ∧
\DeclareUnicodeCharacter{2228}{\lor}               % ∨
\DeclareUnicodeCharacter{2229}{\cap}               % ∩
\DeclareUnicodeCharacter{222A}{\cup}               % ∪
\DeclareUnicodeCharacter{222B}{\int}               % ∫
\DeclareUnicodeCharacter{223C}{\sim}               % ∼
\DeclareUnicodeCharacter{2243}{\simeq}             % ≃
\DeclareUnicodeCharacter{2245}{\cong}              % ≅
\DeclareUnicodeCharacter{2248}{\approx}            % ≈
\DeclareUnicodeCharacter{2260}{\neq}               % ≠
\DeclareUnicodeCharacter{2261}{\equiv}             % ≡
\DeclareUnicodeCharacter{2264}{\le}                % ≤
\DeclareUnicodeCharacter{2265}{\ge}                % ≥
\DeclareUnicodeCharacter{2282}{\subset}            % ⊂
\DeclareUnicodeCharacter{2283}{\supset}            % ⊃
\DeclareUnicodeCharacter{2286}{\subseteq}          % ⊆
\DeclareUnicodeCharacter{2287}{\supseteq}          % ⊇
\DeclareUnicodeCharacter{2293}{\sqcap}             % ⊓
\DeclareUnicodeCharacter{2294}{\sqcup}             % ⊔
\DeclareUnicodeCharacter{2295}{\oplus}             % ⊕
\DeclareUnicodeCharacter{2296}{\ominus}            % ⊖
\DeclareUnicodeCharacter{2297}{\otimes}            % ⊗
\DeclareUnicodeCharacter{2298}{\oslash}            % ⊘
\DeclareUnicodeCharacter{2299}{\odot}              % ⊙
\DeclareUnicodeCharacter{22A2}{\vdash}             % ⊢
\DeclareUnicodeCharacter{22A3}{\dashv}             % ⊣
\DeclareUnicodeCharacter{22A4}{\top}               % ⊤
\DeclareUnicodeCharacter{22A5}{\bot}               % ⊥
\DeclareUnicodeCharacter{22A8}{\vDash}             % ⊨
\DeclareUnicodeCharacter{22C0}{\bigwedge}          % ⋀
\DeclareUnicodeCharacter{22C1}{\bigvee}            % ⋁
\DeclareUnicodeCharacter{22C2}{\bigcap}            % ⋂
\DeclareUnicodeCharacter{22C3}{\bigcup}            % ⋃
\DeclareUnicodeCharacter{22C4}{\lozenge}           % ⋄
\DeclareUnicodeCharacter{22C5}{\Box}               % ⋅
\DeclareUnicodeCharacter{2329}{\langle}            % 〈
\DeclareUnicodeCharacter{232A}{\rangle}            % 〉
\DeclareUnicodeCharacter{2581}{\_}                 % ▁
\DeclareUnicodeCharacter{25A1}{\Box}               % □
\DeclareUnicodeCharacter{25FB}{\Box}               % ◻
\DeclareUnicodeCharacter{266D}{\flat}              % ♭
\DeclareUnicodeCharacter{266E}{\natural}           % ♮
\DeclareUnicodeCharacter{266F}{\sharp}             % ♯
\DeclareUnicodeCharacter{27E6}{\llbracket}         % ⟦
\DeclareUnicodeCharacter{27E7}{\rrbracket}         % ⟧
\DeclareUnicodeCharacter{2806}{{:}}                % ⠆

\fi
\usepackage{edrx21}               % (find-LATEX "edrx21.sty")
\def\eqP{\sim_P}
\def\eqJ{\sim_J}
\def\eqL{\sim_L}
\def\eqR{\sim_R}
\def\eqS{\sim_S}

\def\eqQ{\sim_Q}

\def\pile            {\mathsf{pile}}
\def\pileelements#1#2{{\{#1▁, \ldots, 1▁, \; ▁1, \ldots, ▁#2\}}}

\def\Int{{\operatorname{int}}}
\def\Int{{\operatorname{\mathsf{int}}}}

\def\pile{\mathsf{pile}}

\def\myresizebox#1{%
  \noindent\hbox to \textwidth{\hss
    \resizebox{1.0\textwidth}{!}{#1}%
    \hss}
  }

\def\squigbij{\newsquigbij}
\def\newsquigbij{\;\; \squigbijbody \;\;}
\def\squigbijy{-1.2}
\def\squigbijbody{\squigbijbodywithparams{1.5pt}{0.3pt}{1.0}}
\def\squigbijtriangle(#1,#2)#3{\polygon*(#1,0)(#2,#3)(#2,-#3)}
\def\squigbijbodywithparams#1#2#3{{%
  \unitlength=#1
  \linethickness{#2}
  \begin{picture}(22.4,2.4)(-5.2,\squigbijy)%
    \polyline(-3,0)(0,0)(1,1)(3,-1)(5,1)(7,-1)(9,1)(11,-1)(12,0)(14,0)
    \squigbijtriangle(-5,-2){#3}
    \squigbijtriangle(17,14){#3}
  \end{picture}%
  }}

\usepackage{proof}   % For derivation trees ("%:" lines)
\input diagxy        % For 2D diagrams ("%D" lines)
\xyoption{curve}     % For the ".curve=" feature in 2D diagrams
\usepackage[backend=biber,
   style=alphabetic]{biblatex}            % (find-es "tex" "biber")
\begin{document}

\ifluatex
  \catcode`\^^J=10
  \directlua{dofile "dednat6load.lua"}
\else
  \input\jobname.dnt   % (find-LATEXfile "2021lindenhovius-j-to-X.dnt")
  \def\pu{}
\fi

%\directlua{dofile "2020dn6-error-handling.lua"} % (find-LATEX "2020dn6-error-handling.lua")

%L forths["<-'"] = function () pusharrow("<-^{) }") end
%L forths["<-->"] = function () pusharrow("<-->") end

%L kite  = ".1.|2.3|.4.|.5."
%L mp = MixedPicture.new({def="dagKite", meta="s", scale="6pt"}, z):zfunction(kite):output()
%L -- mp = MixedPicture.new({def="dagKite", meta="s", scale="10pt"}, z):zfunction(kite):output()
%L -- mp = MixedPicture.new({def="dagKote", meta="s", scale="10pt"}, z):zfunction(kite):output()
\pu

% %L dofile "edrxtikz.lua"  -- (find-LATEX "edrxtikz.lua")
% %L dofile "edrxpict.lua"  -- (find-LATEX "edrxpict.lua")
% \pu

% «defs»  (to ".defs")

\def\Incs     {\mathsf{Incs}}
\def\inc      {\mathsf{inc}}
\def\CanSub   {\mathsf{CanSub}}
\def\SubPoints{\mathsf{SubPoints}}
\def\Subsets  {\mathsf{Subsets}}
\def\Can      {\mathsf{Can}}
\def\can      {\mathsf{can}}
\def\Eq       {\mathsf{Eq}}
\def\eq       {\mathsf{eq}}
\def\Logic    {\mathsf{Logic}}
\def\Po       {\mathsf{Po}}
\def\Ob       {\mathsf{Ob}}
\def\catF     {\mathbf{F}}
\def\catH     {\mathbf{H}}
\def\catK     {\mathbf{K}}
\def\catT     {\mathbf{T}}
\def\catW     {\mathbf{W}}
\def\Clo      {\textsc{Cl}}
\def\univ     {\text{univ}}
\def\pcdot    {(·)}
\def\opcdot   {\ovl{(·)}}
\def\clop     {\ovl{(·)}}
\def\SetC     {\Set^\catC}
\def\SetD     {\Set^\catD}
\def\SetH     {\Set^\catH}
\def\SetK     {\Set^\catK}
\def\SetT     {\Set^\catT}
\def\calW     {\mathcal{W}}
\def\toM      {\ton{\text{M}}}

\def\Nucs  {\mathsf{Nucs}}
\def\GrTops{\mathsf{GrTops}}
\def\LTTops{\mathsf{LTTops}}
\def\Clops {\mathsf{Clops}}
\def\dnu  {{↓}u}
\def\dnou {{↓}^\circ u}
\def\dnous{(\dnou)^*}
\def\dnus {(\dnu)^*}

\def\Ddp   {\Downs({↓}p)}
\def\dnp   {{↓}p}
\def\dnpo  {{↓}p∖\{p\}}

\def\nuc      {(·)^*}
\def\nus      {(·)^𝓢}

\def\Jcan  {{J_\mathrm{can}}}
\def\hasmax{\mathsf{hasmax}}
\def\trans {\mathsf{trans}}
\def\stab  {\mathsf{stab}}

\def\Above   {\mathop{\textsf{above}}}
\def\Covers  {\mathop{\textsf{covers}}}
\def\JCovers {\mathop{\textsf{$J$-covers}}}
\def\RCovers {\mathop{\textsf{r-covers}}}
\def\RJCovers{\mathop{\textsf{r-$J$-covers}}}

\def\Downs    {\mathsf{D}}
\def\Cst      {\mathsf{CST}}

\def\rotl#1{\rotatebox{90}{$#1$}}

\def\eqP{\sim_P}
\def\eqJ{\sim_J}
\def\eqL{\sim_L}
\def\eqR{\sim_R}
\def\eqS{\sim_S}

\def\TODO{(TODO)}

\def\standout#1{\fcolorbox{red}{yellow}{#1}}

% «theorem»  (to ".theorem")
% (find-es "tex" "newtheorem")
% (find-es "tex" "newcounter")
%
\newtheorem{thmsection}      {My Theorem}[section]
\newtheorem{thmsubsection}   {My Theorem}[subsection]
\newtheorem{thmsubsubsection}{My Theorem}[subsubsection]
\def\stepThm#1{\refstepcounter{#1}\csname the#1\endcsname}
\def\Theoremsection         {{\bf Theorem \stepThm{thmsection}. }}
\def\Theoremsubsection      {{\bf Theorem \stepThm{thmsubsection}. }}
\def\Theoremsubsubsection   {{\bf Theorem \stepThm{thmsubsubsection}. }}
\def\Lemmasection           {{\bf Lemma   \stepThm{thmsection}. }}
\def\Lemmasubsection        {{\bf Lemma   \stepThm{thmsubsection}. }}
\def\Lemmasubsubsection     {{\bf Lemma   \stepThm{thmsubsubsection}. }}
\def\Conjecturesection      {{\bf Conjecture   \stepThm{thmsection}. }}
\def\Conjecturesubsection   {{\bf Conjecture   \stepThm{thmsubsection}. }}
\def\Conjecturesubsubsection{{\bf Conjecture   \stepThm{thmsubsubsection}. }}
\def\Examplesection         {{\bf Example \stepThm{thmsection}. }}
\def\Examplesubsection      {{\bf Example \stepThm{thmsubsection}. }}
\def\Examplesubsubsection   {{\bf Example \stepThm{thmsubsubsection}. }}

% \newcounter{Counts}[section]
% \newcounter{Countss}[subsection]
% \newcounter{Countsss}[subsubsection]
% \def\stepCount#1{\refstepcounter{#1}\csname the#1\endcsname}
% \def\Theoremsection      {{\bf Theorem \stepThm{Counts}.}}
% \def\Theoremsubsection   {{\bf Theorem \stepThm{Countss}.}}
% \def\Theoremsubsubsection{{\bf Theorem \stepThm{Countsss}.}}

%  _____ _ _   _      
% |_   _(_) |_| | ___ 
%   | | | | __| |/ _ \
%   | | | | |_| |  __/
%   |_| |_|\__|_|\___|
%                     
% «title»  (to ".title")

\title{Each closure operator induces a topology \\ and vice-versa
  (``version for children'')}

\author{Eduardo Ochs}

% http://angg.twu.net/math-b.html\#favco

\maketitle

%     _    _         _                  _   
%    / \  | |__  ___| |_ _ __ __ _  ___| |_ 
%   / _ \ | '_ \/ __| __| '__/ _` |/ __| __|
%  / ___ \| |_) \__ \ |_| | | (_| | (__| |_ 
% /_/   \_\_.__/|___/\__|_|  \__,_|\___|\__|
%                                           
% «abstract»  (to ".abstract")
% (favp 1 "abstract")
% (fav    "abstract")

\begin{abstract}

  % (find-books "__cats/__cats.el" "johnstone-topostheory")
  % (find-toposthepubpage 108  "3.1. Topologies")
  % (find-toposthepubtext 108  "3.1. Topologies")
  % (find-books "__cats/__cats.el" "mclarty")
  % (find-mclartypage (+ 4 196) "21. Topologies")
  % (find-books "__cats/__cats.el" "maclane-moerdijk")
  % (find-maclanemoerdijkpage (+   11 219) "V.1 Lawvere-Tierney Topologies")
  % (find-books "__cats/__cats.el" "lambek-scott")
  % (find-lambekscottpage (+ 8 200) "Toposes with canonical subobjects")
  % (find-books "__cats/__cats.el" "freyd-aspects")
  % (find-freydatpage 16 "1.4. Modal operators in Heyting Algebras")
  % (find-books "__cats/__cats.el" "fourman-scott")
  % (find-slnm0753page (+ 16 324) "J-operators")

  One of the main prerequisites for understanding sheaves on
  elementary toposes is the proof that a (Lawvere-Tierney) topology on
  a topos induces a closure operator on it, and vice-versa. That
  standard theorem is usually presented in a relatively brief way,
  with most details being left to the reader --- see for example
  \cite[section 3.1]{Johnstone}, \cite[chapter 21]{McLarty},
  \cite[section V.1]{MacLaneMoerdijk}, \cite[chapter 5]{BellLST} ---
  and with no hints on how to visualize some of the hardest axioms and
  proofs.

  These notes are, on a first level, an attempt to present that
  standard theorem in all details and in a visual way, following the
  conventions in \cite{FavC}; in particular, some properties, like
  stability by pullbacks, are always drawn in the same ``shape''.

  On a second level these notes are also an experiment on doing these
  proofs on ``archetypal cases'' in ways that makes all the proofs
  easy to lift to the ``general case''. Our first archetypal case is a
  ``topos with inclusions''. This is a variant of the ``toposes with
  canonical subobjects'' from \cite[section 2.15]{LambekScott}; all
  toposes of the form $\SetC$, where $\catC$ is a small category, are
  toposes with inclusions, and when we work with toposes with
  inclusions concepts like subsets and intersections are very easy to
  formalize. We do all our proofs on the correspondence between
  closure operators and topologies in toposes with inclusions, and
  then we show how to lift all our proofs to proofs that work on any
  topos. Our second archetypal case is toposes of the form $\SetD$,
  where $𝐃$ is a finite two-column graph. We show a way to visualize
  all the Lawvere-Tierney topologies on toposes of the form $\SetD$,
  and we define formally a way to ``add visual intuition to a proof
  about toposes''; this is related to the several techniques for doing
  ``Category Theory for children'' that are explained in the first
  sections of \cite{FavC}.

  % We also use the idea that ``$\Set$ is the archetypal topos'' (from
  % \cite[section 16]{IDARCT}) and a variant of the ``canonical
  % subobjects'' from \cite[section 2.15]{LambekScott} to do a version
  % ``for children'' of the proof of the correspondence between
  % topologies and closure operators; this proof ``for children'' can
  % be lifted without much pain to a proof that works on toposes
  % without canonical subobjects.

  % These notes are an attempt to present that standard theorem in all
  % details and in a visual way, following the conventions in
  % \cite{FavC}; in particular, some properties, like stability by
  % pullbacks, are always drawn in the same ``shape''. We also use the
  % idea that ``$\Set$ is the archetypal topos'' (from \cite[section
  % 16]{IDARCT}) and a variant of the ``canonical subobjects'' from
  % \cite[section 2.15]{LambekScott} to do a version ``for children''
  % of the proof of the correspondence between topologies and closure
  % operators; this proof ``for children'' can be lifted without much
  % pain to a proof that works on toposes without canonical
  % subobjects.

  % The last sections of these notes show how, for certain toposes,
  % the operation that restricts a closure operation on a topos to its
  % action on $\Sub(1)$ --- i.e., to a ``modal operator on its Heyting
  % Algebra of truth values'' (\cite[section 1.4]{FreydAspects}), also
  % called a ``J-operator'' in \cite[definition 2.11]{Fourman} and
  % \cite{PH2} --- is a bijection, and shows how to visualize this.
  % {\sl I haven't been able to find mentions of this bijection in the
  % literature... if you know any, please let me know!}

\end{abstract}

\newpage

% «toc»  (to ".toc")
% (cltp 2 "toc")
% (clta   "toc")

% (find-es "tex" "tocloft")
% (favp 2 "toc")
% (fav    "toc")
\renewcommand{\cfttoctitlefont}{\bfseries}
\setlength{\cftbeforesecskip}{2.5pt}
\tableofcontents

\vspace*{2.0cm}

{\sl Status of these notes:} this is not yet in final form. The last
section still needs to be written, and there are a few proofs in the
last sections are marked with ``TODO'' --- but all things that still
need to be done are clearly marked, so that's not so bad. For the most
recent version look here:

\url{http://angg.twu.net/math-b.html\#clops-and-tops}

% \fbox{\parbox{12.5cm}{
% 
% {\sl Status of these notes:} this is a work in progress. The sections
% that still need changes before I upload a version of this to Arxiv are
% marked with ``(TODO)''. The compilation date of this PDF is at the
% footer. For the most recent version see:
% 
% \url{http://angg.twu.net/math-b.html\#clops-and-tops}
% 
% \ssk
% 
% My e-mail is eduardoochs@gmail.com.
% }}

\newpage

%  ____            _            _       __     
% | __ )  __ _ ___(_) ___    __| | ___ / _|___ 
% |  _ \ / _` / __| |/ __|  / _` |/ _ \ |_/ __|
% | |_) | (_| \__ \ | (__  | (_| |  __/  _\__ \
% |____/ \__,_|___/_|\___|  \__,_|\___|_| |___/
%                                              
% «inclusions»  (to ".inclusions")
% (cltp 4 "inclusions")
% (clta   "inclusions")
\section{Subobjects and inclusions}
\label{inclusions}

The {\sl subobjects} of an object $D$ of a topos $\catE$ are the
monics with codomain $D$ {\sl modulo isomorphism}. Here is an example
in $\Set$:
%
%D diagram subobjects-example
%D 2Dx     100  +35  +35  +35  +35  +35
%D 2D  100 A0 - A1 - A2   B0 - B1 - B2 
%D 2D         \ |  /         \ |  /    
%D 2D  +35      A3             B3      
%D 2D
%D ren A0 A1 A2 A3 ==> A B C D
%D ren B0 B1 B2 B3 ==> \{20,40\} \{2,4\} \{2,4\} \{1,2,3,4\}
%D
%D (( A0 A1 <-> A1 A2 <->
%D    A0 A3 >-> .plabel= l f
%D    A1 A3 `-> .plabel= r g
%D    A2 A3 >-> .plabel= r h
%D ))
%D (( B0 B1 -> sl^ .plabel= a \sm{20↦2\\40↦4}
%D    B0 B1 <- sl_
%D    B1 B2 -> sl^ .plabel= a \sm{2↦4\\4↦2}
%D    B1 B2 <- sl_
%D    B0 B3 >-> .plabel= l \sm{20↦2\\20↦4}
%D    B1 B3 `->
%D    B2 B3 >-> .plabel= r \sm{2↦4\\4↦2}
%D ))
%D enddiagram
%D
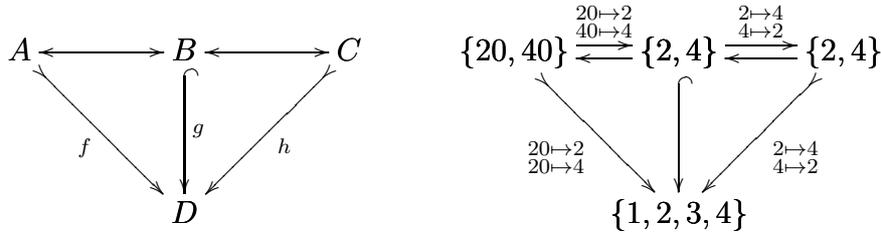
$$\pu
  \diag{subobjects-example}
$$

Here the monics $f:A \monicto D$, $g:B \monicto D$, $h:C \monicto D$
are all equivalent; in some texts they are ``the same subobject''.
Let's make that precise. For us the elements of $\Sub(D)$ are the
monics with codomain $D$. If $(f:A \monicto D)$, $(g:B \monicto D)$
are elements of $\Sub(D)$ then they are {\sl equivalent} (notation:
$f≡g$) iff there is an iso $A↔B$ making the obvious triangle commute.
We write $[f]$ for the equivalence class made of an $f∈\Sub(D)$ and
all other monics in $\Sub(D)$ equivalent to $f$, and we write
$\ovl{\Sub}(D)$ for $\Sub(D)$ modulo equivalence: so
$[f]∈\ovl{\Sub}(D)$.

A monic $g:B \monicto D$ in $\Set$ is an {\sl inclusion} if it obeys:
$$∀b∈\dom(g).\, g(b) = b.$$

% (find-books "__cats/__cats.el" "lambek-scott")
% (find-lambekscottpage (+ 8 200) "Toposes with canonical subobjects")

The usual way to formalize inclusions in toposes is via canonical
subobjects. A topos $\catE$ {\sl has canonical subobjects} is it comes
equiped with a class $\CanSub(\catE)$ of monics that obey a certain
list of properties --- see \cite[p.200 onwards]{LambekScott} --- that
are also obeyed by the inclusions in $\Set$. Here we will do something
similar but with a different list of properties, and in section
\ref{without-inclusions} we will see how to translate our proofs, done
in toposes with inclusions, to proofs in arbitrary toposes.

% Here we will do that backwards: we will work with inclusions using
% some of the properties that they obviously have in $\Set$, and assure
% the reader that everything here works with ``inclusion'' replaced by
% ``canonical subobject'' everywhere; checking that is left as an
% exercise.

% If $f: A \ito B$ is an inclusion than we say that $f$ is an {\sl
%   inclusion map} and that $A$ is a {\sl subset} of $B$. In $\Set$ for
% any two objects $A$ and $B$ we have at most one inclusion map from $A$
% to $B$; this property --- that is not in \cite{LambekScott} --- lets
% us talk of {\sl the} intersection and {\sl the} union of subsets
% $A,B⊆C$ and to denote them by $A∩B$ and $A∪C$, omitting all the names
% of the arrows in the diagrams:
% %
% %D diagram cap-and-cup
% %D 2Dx     100 +25 +25 +25 +25
% %D 2D  100 A0  A1  B0  B1  B2
% %D 2D
% %D 2D  +25 A2  A3      B3
% %D 2D
% %D ren A0 A1 A2 A3 ==> A∩B B A C
% %D ren B0 B1 B2 B3 ==> A A∪B B C
% %D
% %D (( A0 A1 `->
% %D    A0 A2 `->
% %D    A1 A3 `->
% %D    A2 A3 `->
% %D    A0 relplace 7 7 \pbsymbol{7}
% %D
% %D    B0 B1 `->
% %D    B1 B2 <-'
% %D    B0 B3 `->
% %D    B1 B3 `->
% %D    B2 B3 `->
% %D ))
% %D enddiagram
% %D
% $$\pu
%   \diag{cap-and-cup}
% $$

\msk

When $f: A \monicto C$ and $g: B \monicto C$ are subobjects of $C$ we
say that {\sl $f$ is contained in $g$} (notation: $f⊆g$) when there is
a monic $m: A \monicto B$ making the obvious triangle commute. We call
$m$ the ``mediating map''.

\msk

In $\Set$ we have two different operations that take two maps $f,g$
with a common codomain and produce pullbacks:
%
%D diagram two-different-pullback-ops
%D 2Dx     100  +25  +70  +60
%D 2D  100      A1   B0 - B1
%D 2D            |   |     |
%D 2D  +25 A2 - A3   B2 - B3
%D 2D
%D 2D  +20      C1   D0 - D1
%D 2D            |   |     |
%D 2D  +25 C2 - C3   D2 - D3
%D 2D
%D ren A1 A2 A3    ==> B A C
%D ren B1 B2 B3 B0 ==> B A C \setofst{(a,b){∈}A{×}B}{f(a){=}g(b)}
%D ren C1 C2 C3    ==> B A C
%D ren D1 D2 D3 D0 ==> B A C \setofst{a{∈}A}{f(a)∈B}
%D
%D (( A1 A3 -> .plabel= r g
%D    A2 A3 -> .plabel= b f
%D
%D    B0 B1 -> .plabel= a π'
%D    B0 B2 -> .plabel= l π
%D    B1 B3 -> .plabel= r g
%D    B2 B3 -> .plabel= b f
%D
%D    C1 C3 -> .plabel= r g
%D    C2 C3 -> .plabel= b f
%D
%D    D0 D1  ->
%D    D0 D2 `->
%D    D1 D3 `-> .plabel= r g
%D    D2 D3  -> .plabel= b f
%D
%D    A1 B2 harrownodes nil 15 30 |->
%D    C1 D2 harrownodes nil 15 30 |->
%D ))
%D enddiagram
%D
$$\pu
  \diag{two-different-pullback-ops}
$$

The second one only works when the right wall is an inclusion, but it
produces pullbacks whose left walls are inclusions. In both cases we
will write the left wall as $f^{-1}(g): f^{-1}(B) \to A$,

%D diagram pullback-left-wall
%D 2Dx     100  +30
%D 2D  100 A0 - A1 
%D 2D      |     | 
%D 2D  +25 A2 - A3 
%D 2D
%D ren A0 A1 A2 A3 ==> f^{-1}(B) B A C
%D
%D (( A0 A1 ->
%D    A0 A2 -> .plabel= l f^{-1}(g)
%D    A1 A3 -> .plabel= r g
%D    A2 A3 -> .plabel= b f
%D ))
%D enddiagram
%D
$$\pu
  \diag{pullback-left-wall}
$$
and there will be no default name for the top wall. When the right
wall is marked as an inclusion we will use the second pullback
operation, otherwise the first one.

In $\Set$ the classifying map of a monic $m: A \monicto B$ is defined
as:
%
%D diagram classifying-map-def-in-Set
%D 2Dx     100  +60
%D 2D  100 A0 - A1 
%D 2D      |     | 
%D 2D  +25 A2 - A3 
%D 2D
%D ren A0 A1 A2 A3 ==> A 1 B Ω
%D
%D (( A0 A1  ->
%D    A0 A2 >-> .plabel= l m
%D    A1 A3 `-> .plabel= r ⊤
%D    A2 A3  -> .plabel= b \sm{χ_f:=\\(λb:B.∃a∈A.m(a)=b)}
%D ))
%D enddiagram
%D
$$\pu
  \diag{classifying-map-def-in-Set}
$$
and the ``true'' map $⊤:1 \ito Ω$ is the inclusion $\{1\} \ito
\{0,1\}$.

The {\sl inclusion classified by a map $f: B \to Ω$} is the map
$f^{-1}(⊤)$; we will sometimes write it as $σ(f)$. Note that for any
monic $m: A \monicto B$ we have $m ≡ σ(χ_m)$, and we have $m = σ(χ_m)$
if $m$ is an inclusion; and for any $f: B \to Ω$ we have $χ_{σ(f)} =
f$.

\newpage

%  ___                                      _          _       
% |_ _|_ __   ___ ___   _ __  _ __ ___  ___(_)___  ___| |_   _ 
%  | || '_ \ / __/ __| | '_ \| '__/ _ \/ __| / __|/ _ \ | | | |
%  | || | | | (__\__ \ | |_) | | |  __/ (__| \__ \  __/ | |_| |
% |___|_| |_|\___|___/ | .__/|_|  \___|\___|_|___/\___|_|\__, |
%                      |_|                               |___/ 
%
% «inclusions-precisely»  (to ".inclusions-precisely")
% (cltp 5 "inclusions-precisely")
% (clta   "inclusions-precisely")
\subsection{Inclusions, precisely}
\label{inclusions-precisely}

A {\sl topos with inclusions} is a topos $\catE$ endowed with a class
of monics $\Incs(\catE)$, called the {\sl inclusions}, and two
pullback operations, as in the previous section, obeying the
properties below:

\begin{itemize}

\item[Inc1)] For any two object $C$ and $D$ of $\catE$ there is at
  most one inclusion from $C$ to $D$. When that inclusion map exists
  we write it as $C \ito D$ --- we don't need to name it --- and we
  say that {\sl $C$ is a subset of $D$} (notation: $C⊆D$).

\item[Inc2)] Each $[f]∈\ovl{\Sub}(D)$ contains exactly one inclusion
  map. This can be expressed as
%
%D diagram exactly-one-inclusion
%D 2Dx     100  +25
%D 2D  100 A0 - A1
%D 2D         \  |
%D 2D  +25      A2
%D 2D
%D ren A0 A1 A2 ==> ∀A ∃!B ∀D
%D
%D (( A0 A1 <-> .plabel= a ∃!
%D    A0 A2 >-> .plabel= l ∀f
%D    A1 A2 `-> .plabel= r ∃!g
%D ))
%D enddiagram
%D
$$\pu
  \diag{exactly-one-inclusion}
$$
in the variant of Freyd's diagrammatic language defined in
\cite[section 4.1]{FavC}. We will say that this $g$ is {\sl the
  inclusion associated} (or: {\sl equivalent}) {\sl to $f$}, and write
this as $\can(f)=g$.

\item[Inc3)] The composite of two inclusions is an inclusion. Or, in
  the language of Inc1: if $B⊆C$ and $C⊆D$ then $B⊆D$, with $B \ito D
  = B \ito C \ito D$.

\item[Inc4)] If $f: B \ito D$ and $g: C \ito D$ are inclusions with
  $f⊆g$ then the mediating map $m: B \monicto C$ is an inclusion. In
  the language of Inc1: $f⊆g$ implies $B⊆C$. We can visualize this as:
%
%D diagram mediating-map-incl
%D 2Dx     100 +15 +15
%D 2D  100 A0      A1
%D 2D
%D 2D  +25     A2
%D 2D
%D ren A0 A1 A2 ==> B C D
%D
%D (( A0 A1 `-> .plabel= a m
%D    A0 A2 `-> .plabel= l f
%D    A1 A2 `-> .plabel= r g
%D ))
%D enddiagram
%D
$$\pu
  \diag{mediating-map-incl}
$$

\item[Inc5)] The ``true'' map $⊤: 1 \ito Ω$ is an inclusion.

\item[Inc6)] The second operation that produces pullbacks in $\catE$
  receives maps $f: A \to C$ and $g: B \ito C$ and returns pullbacks
  whose left walls are inclusions. In a diagram:
%
%D diagram second-pullback-op
%D 2Dx     100  +25  +45  +25
%D 2D  100      A1   B0 - B1
%D 2D            |   |     |
%D 2D  +25 A2 - A3   B2 - B3
%D 2D
%D ren    A1 A2 A3 ==> B A C
%D ren B0 B1 B2 B3 ==> f^{-1}(B) B A C
%D
%D (( A1 A3 `-> .plabel= r g
%D    A2 A3  -> .plabel= b f
%D
%D    B0 B1  ->
%D    B0 B2 `-> .plabel= l f^{-1}(g)
%D    B1 B3 `-> .plabel= r g
%D    B2 B3  -> .plabel= b f
%D
%D    A1 B2 harrownodes nil 15 15 |->
%D ))
%D enddiagram
%D
$$\pu
  \diag{second-pullback-op}
$$

\item[Inc7)] The {\sl intersection} of two inclusions $B \ito D$ and
  $C \ito D$ is defined as their pullback:
%
%D diagram intersection-of-inclusions
%D 2Dx     100  +25
%D 2D  100 A0 - A1 
%D 2D      |     | 
%D 2D  +25 A2 - A3 
%D 2D
%D ren A0 A1 A2 A3 ==> B{∩}C C B D
%D
%D (( A0 A1 `->
%D    A0 A2 `->
%D    A1 A3 `->
%D    A2 A3 `->
%D ))
%D enddiagram
%D
$$\pu
  \diag{intersection-of-inclusions}
$$
Note that its upper wall is the mediating map from the composite
$B{∩}C \ito C \ito D$ to $C \ito D$, so it is an inclusion.

Using Inc2 we can see that $B{∩}C$ and $C{∩}B$ are the {\sl same}
subset of $D$, not just isomorphic subobjects.

\end{itemize}

We write $\Incs(D)$ for the class of inclusions with codomain $D$ and
$\Subsets(D)$ for the class of subsets of $D$. In a topos with
inclusions we have:
$$\Subsets(D)
  ≅ \Incs(D)
  ≅ \ovl{\Sub}(D)
  ≅ \Sub(D),
$$
where the first two `$≅$'s are isomorphisms and the last one is just
an ``equivalence of categories'': if we start with a monic $f$ in
$\Sub(D)$, take it to its equivalence class $[f]$ in $\ovl{\Sub}(D)$,
and then go back to $\Sub(D)$, what we get is $\can(f)$, and we have
$f≡\can(f)$ but not necessarily $f=\can(f)$.

%  ___                  ____       _    ____ 
% |_ _|_ __   ___ ___  / ___|  ___| |_ / ___|
%  | || '_ \ / __/ __| \___ \ / _ \ __| |    
%  | || | | | (__\__ \  ___) |  __/ |_| |___ 
% |___|_| |_|\___|___/ |____/ \___|\__|\____|
%                                            
% «inclusions-SetC»  (to ".inclusions-SetC")
% (cltp 6 "inclusions-SetC")
% (clta   "inclusions-SetC")
\subsection{Inclusions in a topos of the form $\SetC$}
\label{inclusions-SetC}

From here onwards $\catC$ will always denote a small category.

All toposes of the form $\SetC$ are toposes with inclusions. We define
the class of inclusions of $\SetC$, $\Incs(\SetC)$, as follows:

\begin{itemize}

\item[IncSC)] A morphism $i:A \ito B$ in $\SetC$ is an inclusion iff
  for every object $u$ of $\catC$ the map $i_u:A(u) \ito B(u)$ is an
  inclusion in $\Set$; that is, if $∀u∈\catC. ∀a∈A(u). i_u(a)=a$.

\end{itemize}

% TODO: fix the indentation above. See:
% (find-es "tex" "quote")

Take an inclusion $i:A \ito B$ and a morphism $v:u→w$ in $\catC$. As
$i$ is a natural transformation, these squares commute:
%
% (favp 22 "internal-view-NT" "5.3.")
% (fav     "internal-view-NT")
%
%D diagram IncSC
%D 2Dx     100 +20  +30 +25  +35
%D 2D  100 A0  B0 - B1  D0 - D1
%D 2D      |   |     |  |     |
%D 2D  +25 A1  B2 - B3  D2 - D3
%D 2D
%D 2D  +15     C0 - C1
%D 2D
%D ren A0 A1       ==> u w
%D ren B0 B1 B2 B3 ==> A(u) B(u) A(w) B(w)
%D ren C0 C1       ==> A B
%D ren D0 D1 D2 D3 ==> a a A(v)(a) A(v)(a)
%D
%D (( A0 A1  -> .plabel= l v
%D    B0 B1 `-> .plabel= a i_u
%D    B0 B2  -> .plabel= l A(v)
%D    B1 B3  -> .plabel= r B(v)
%D    B2 B3 `-> .plabel= a i_v
%D    C0 C1 `-> .plabel= a f
%D
%D    newnode: D3' at: @D3+v(0,-7) .TeX= B(v)(a) place
%D    D0 D1 |->
%D    D0 D2 |->
%D    D1 D3' |->
%D    D2 D3 |->
%D ))
%D enddiagram
%D
$$\pu
  \diag{IncSC}
$$

This means that $∀a∈A(u).A(v)(a)=B(v)(a)$ -- so $A(v):A(u)→A(w)$ is a
restriction of the function $B(v):B(u)→B(w)$ to $A(u)$.

%     _              _                   _   _                 _ _           
%    / \   _ __   __| |   __ _ _ __   __| | (_)_ __ ___  _ __ | (_) ___  ___ 
%   / _ \ | '_ \ / _` |  / _` | '_ \ / _` | | | '_ ` _ \| '_ \| | |/ _ \/ __|
%  / ___ \| | | | (_| | | (_| | | | | (_| | | | | | | | | |_) | | |  __/\__ \
% /_/   \_\_| |_|\__,_|  \__,_|_| |_|\__,_| |_|_| |_| |_| .__/|_|_|\___||___/
%                                                       |_|                  
%
% «and-and-implies»  (to ".and-and-implies")
% (cltp 6 "and-and-implies")
% (clta   "and-and-implies")
\subsection{`And' and `implies'}

\def\setofPQst#1{\setofst{(P,Q)∈Ω×Ω}{#1}}

In section \ref{top-to-clop} we will need the ``internal conjunction
map'', $(∧):Ω×Ω→Ω$, whose internal view is $(P,Q) \mapsto P∧Q$, and
the ``internal implication map'', $({→}):Ω×Ω→Ω$, that works as $(P,Q)
\mapsto (P{→}Q)$. They are well explained in sections 13.3 and 13.4 of
\cite{McLarty}, but only in their forms ``for adults'', that work in
arbitrary toposes. In this section I will just complement
\cite{McLarty} by showing briefly how those definitions that hold in
any topos are translations of definitions that make sense in $\Set$.

% (find-books "__cats/__cats.el" "mclarty")
% (find-mclartypage (+ 4 118) "13.3. Conjunction and intersection")
% (find-mclartypage (+ 4 119) "13.4. Order and implicates")

The arrow $(∧)$ is built as the classifying map of the inclusion
$σ(∧)$ in this diagram,
%
%D diagram definition-of-and
%D 2Dx     100  +60
%D 2D  100 A0 - A1
%D 2D      |     |
%D 2D  +25 A2 - A3
%D 2D
%D ren A0 A1 A2 A3 ==> \setofPQst{P∧Q} 1 Ω×Ω Ω
%D
%D (( A0 A1  ->
%D    A0 A2 `-> .plabel= l σ(∧)
%D    A1 A3 `->
%D    A2 A3  -> .plabel= b ∧:=χ_{σ(∧)}
%D
%D ))
%D enddiagram
%D
$$\pu
  \diag{definition-of-and}
$$
and the inclusion $σ(∧)$ is built as an equalizer. We have:
$$\begin{array}{l}
  % \setofst{P,Q}{{∧}} \\
  % := \;\;
  \setofPQst{P∧Q} \\
  = \;\; \setofPQst{P=⊤∧Q=⊤} \\
  = \;\; \setofPQst{\id_Ω(P)=⊤_Ω(P) ∧ \id_Ω(Q)=⊤_Ω(Q)} \\
  = \;\; \setofPQst{(\id_Ω×\id_Ω)(P,Q) = (⊤_Ω×⊤_Ω)(P,Q)} \\
  = \;\; \Eq((\id_Ω×\id_Ω),(⊤_Ω×⊤_Ω)) \\
  \end{array}
$$
%
%D diagram and-as-equalizer
%D 2Dx     100  +120  +40
%D 2D  100 B0 - B1 = B2
%D 2D
%D 2D  +20 A0 - A1 = A2
%D 2D
%D ren A0 A1 A2 ==> \Eq((\id_Ω×\id_Ω),(⊤_Ω×⊤_Ω)) Ω×Ω Ω
%D ren B0 B1 B2 ==> \setofPQst{P∧Q} Ω×Ω Ω
%D
%D (( A0 A1 `-> .plabel= a \eq((\id_Ω×\id_Ω),(⊤_Ω×⊤_Ω))
%D    A1 A2 -> sl^ .plabel= a \id_Ω×\id_Ω
%D    A1 A2 -> sl_ .plabel= b ⊤_Ω×⊤_Ω
%D
%D    B1 xy+= -15 0
%D    B0 B1 `-> .plabel= a σ(∧)
%D    B1 B2 -> sl^ .plabel= a (P,Q)↦(P,Q)
%D    B1 B2 -> sl_ .plabel= b (P,Q)↦(⊤,⊤)
%D ))
%D enddiagram
%D
$$\pu
  \diag{and-as-equalizer}
$$

Where the map $⊤_Ω$ is defined as:
%
%D diagram def-T_Omega
%D 2Dx     100  +20  +20  +20  +20  +20
%D 2D  100 A0 - A1 - A2   B0 - B1 - B2
%D 2D
%D ren A0 A1 A2 ==> A 1 Ω
%D ren B0 B1 B2 ==> Ω 1 Ω
%D
%D (( A0 A1 -> .plabel= a !_A
%D    A1 A2 -> .plabel= a ⊤
%D    A0 A2 -> .slide= -7.5pt .plabel= b ⊤_A:=⊤∘!_A
%D
%D    B0 B1 -> .plabel= a !_Ω
%D    B1 B2 -> .plabel= a ⊤
%D    B0 B2 -> .slide= -7.5pt .plabel= b ⊤_Ω:=⊤∘!_Ω
%D ))
%D enddiagram
%D
$$\pu
  \diag{def-T_Omega}
$$

The arrow $({→})$ is the classifier of the inclusion $σ({→})$, that is
built as another equalizer:
%
%D diagram def-P->Q
%D 2Dx     100  +60
%D 2D  100 A0 - A1
%D 2D      |     |
%D 2D  +25 A2 - A3
%D 2D
%D ren A0 A1 A2 A3 ==> \setofPQst{P→Q} 1 Ω×Ω Ω
%D
%D (( A0 A1  ->
%D    A0 A2 `-> .plabel= l σ({→})
%D    A1 A3 `->
%D    A2 A3  -> .plabel= b ({→}):=χ_{σ({→})}
%D ))
%D enddiagram
%D
$$\pu
  \diag{def-P->Q}
$$

$$\begin{array}{l}
  %\setofst{P,Q}{{→}} \\
  %:= \;\;
  \setofPQst{P→Q} \\
  = \;\; \setofPQst{⊤≤P→Q} \\
  = \;\; \setofPQst{⊤∧P≤Q} \\
  = \;\; \setofPQst{P≤Q} \\
  = \;\; \setofPQst{P=P∧Q} \\
  = \;\; \setofPQst{π(P,Q)=(∧)(P,Q)} \\
  = \;\; \Eq(π,∧) \\
  \end{array}
$$
%
%D diagram def-P->Q-2
%D 2Dx     100  +75  +40
%D 2D  100 B0 - B1 = B2
%D 2D
%D 2D  +20 A0 - A1 = A2
%D 2D
%D ren B0 B1 B2 ==> \setofPQst{P→Q} Ω×Ω Ω
%D ren A0 A1 A2 ==> \Eq(π,∧) Ω×Ω Ω
%D
%D (( A0 A1 `-> .plabel= a \eq(π,∧)
%D    A1 A2 -> sl^ .plabel= a π
%D    A1 A2 -> sl_ .plabel= b ∧
%D
%D  # B1 xy+= -15 0
%D    B0 B1 `-> .plabel= a σ({→})
%D    B1 B2 -> sl^ .plabel= a (P,Q)↦P
%D    B1 B2 -> sl_ .plabel= b (P,Q)↦(P∧Q)
%D ))
%D enddiagram
%D
$$\pu
  \diag{def-P->Q-2}
$$

Note that in $\Set$ we have:
$$\begin{array}{rcl}
  \setofPQst{P∧Q} &=& \{(1,1)\}, \\
  \setofPQst{P→Q} &=& \{(0,0),(0,1),(1,1)\}. \\
  \end{array}
$$

\newpage

% «clops»  (to ".clops")
% (cltp 8 "clops")
% (clta   "clops")
\section{Closure operators}
\label{clops}

A {\sl closure operator} $\clop$ on a topos with inclusions $\catE$ is
a family of operations like this,
$$\begin{array}{rrcl}
  \clop_E: & \Incs(E) &→& \Incs(E) \\
           & (d: D \ito E) &↦& (\ovl{d}: \ovl{D} \ito E), \\
  \end{array}
$$
where we have one $\clop_E$ for each object $E$ of the topos, and
these `$\clop_E$'s obey:

\msk

C1) $d ⊆ \ovl{d}$,

C2) $\ovl{d} = \ovl{\ovl{d}}$,

C3) $c⊆d$ implies $\ovl{c} ⊆ \ovl{d}$,

C4) $\ovl{c∩d} = \ovl{c} ∩ \ovl{d}$,

C5) $\ovl{f^{-1}(d)} = f^{-1}(\ovl{d})$,

\msk

\noindent for all inclusions $c: C \ito E$ and $d: D \ito E$ and for
all maps $f:B→E$.

We will draw the properties C1, C2, C3, C5 as:

%D diagram C1-C3
%D 2Dx     100  +25  +25  +25  +25  +20  +25  +25
%D 2D  100 A0        B0 - B1 - B2   C0 - C1
%D 2D      |  \      |              |  \    \
%D 2D  +25 |    A1   |              |    C2 - C3
%D 2D      |  /      |              |  /
%D 2D  +25 A2        B3             C4
%D 2D
%D ren A0 A1 A2       ==> D \ovl{D} E
%D ren B0 B1 B2 B3    ==> D \ovl{D} \ovl{\ovl{D}} E
%D ren C0 C1 C2 C3 C4 ==> C D \ovl{C} \ovl{D} E
%D
%D (( A0 A1 >->
%D    A0 A2 >-> .plabel= r d
%D    A1 A2 >-> .plabel= r \ovl{d}
%D
%D    B1 xy+= 0 5
%D    B2 xy+= -5 20
%D
%D    B0 B1 >->
%D    B1 B2 =
%D    B0 B3 >-> .plabel= r d
%D    B1 B3 >-> .plabel= r \ovl{d}
%D    B2 B3 >-> .plabel= r \ovl{\ovl{d}}
%D
%D    C0 C1 >->
%D    C0 C2 >->
%D    C1 C3 >->
%D    C2 C3 >->
%D    C0 C4 >-> .PLABEL= ^(0.22) c
%D    C1 C4 >-> .PLABEL= ^(0.22) d
%D    C2 C4 >-> .PLABEL= ^(0.22) \ovl{c}
%D    C3 C4 >-> .PLABEL= ^(0.22) \ovl{d}
%D ))
%D enddiagram
%D
$$\pu
  \diag{C1-C3}
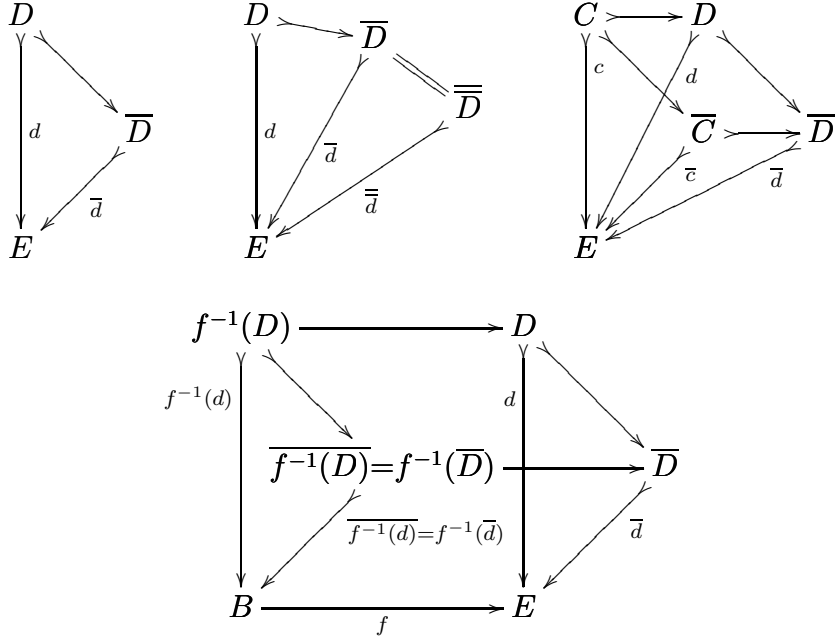
$$

%D diagram C5
%D 2Dx     100  +30 +30  +30
%D 2D  100 A0 ----- B0
%D 2D      |  \     |  \
%D 2D  +30 |    A1 -|--- B1
%D 2D      |  /     |  /
%D 2D  +30 A2 ----- B2
%D 2D
%D ren A0 A1 A2 ==> f^{-1}(D) \ovl{f^{-1}(D)}{=}f^{-1}(\ovl{D}) B
%D ren B0 B1 B2 ==> D \ovl{D} E
%D
%D (( A0 B0 ->
%D    A1 B1 ->
%D    A2 B2 -> .plabel= b f
%D    
%D    A0 A2 >-> .PLABEL= _(0.25) f^{-1}(d)
%D    A0 A1 >->
%D    A1 A2 >-> .PLABEL= ^(0.30) \ovl{f^{-1}(d)}{=}f^{-1}(\ovl{d})
%D    
%D    B0 B2 >-> .PLABEL= _(0.25) d
%D    B0 B1 >->
%D    B1 B2 >-> .PLABEL= ^(0.30) \ovl{d}
%D    
%D ))
%D enddiagram
%D
$$\pu
  \diag{C5}
$$

Where all the `$\monicto$'s in the diagrams are inclusions.

{\bf Important:} in all diagrams from this section to section to
\ref{clop-top-bij} all the `$\monicto$'s will stand for inclusions.
This is for typographical reasons, to make the diagrams a bit lighter.
The distinction between `$\monicto$'s and `$\ito$'s will reappear in
section \ref{without-inclusions}.

\newpage

% «topologies»  (to ".topologies")
% (cltp 9 "topologies")
% (clta   "topologies")

\subsection{Topologies}
\label{topologies}

A {\sl (Lawvere-Tierney) Topology} on a topos $\catE$ is a map $j:Ω→Ω$
obeying:

\ssk

LT1) $j∘j=j$,

LT2) $j∘⊤=⊤$,

LT3) $j∘∧ = ∧∘(j×j)$.

\ssk

We draw LT1, LT2, and LT3 as:
%
% (larp 3 "21._topologies")
% (lar    "21._topologies")
%
%D diagram topologies-1
%D 2Dx     100 +20 +20 +20 +30 +25
%D 2D  100 A0  A1  B0  B1  C0  C1
%D 2D
%D 2D  +20     A2      B2  C2  C3
%D 2D
%D ren A0 A1 A2 ==> Ω Ω Ω
%D ren B0 B1 B2 ==> Ω Ω Ω
%D ren C0 C1 C2 C3 ==> Ω×Ω Ω Ω×Ω Ω
%D
%D (( A0 A1 -> .plabel= a ⊤
%D    A0 A2 -> .plabel= l ⊤
%D    A1 A2 -> .plabel= r j
%D
%D    B0 B1 -> .plabel= a j
%D    B0 B2 -> .plabel= l j
%D    B1 B2 -> .plabel= r j
%D
%D    C0 C1 -> .plabel= a ∧
%D    C0 C2 -> .plabel= l j×j
%D    C1 C3 -> .plabel= r j
%D    C2 C3 -> .plabel= a ∧
%D ))
%D enddiagram
%D
$$\pu
  \diag{topologies-1}
$$

One way to grasp the intuitive meaning of LT1, LT2, and LT3 is to look
at their internal views. If we have maps $p,q:A→Ω$, the internal views
of
%
%D diagram topologies-2
%D 2Dx     100 +20 +20 +15 +20 +20 +15 +30 +25
%D 2D  100 A0  A1  A2  B0  B1  B2  C0  C1  C2
%D 2D
%D 2D  +20         A3          B3      C3  C4
%D 2D
%D ren A0 A1 A2 A3 ==> A Ω Ω Ω
%D ren B0 B1 B2 B3 ==> A Ω Ω Ω
%D ren C0 C1 C2 C3 C4 ==> A Ω×Ω Ω Ω×Ω Ω
%D
%D (( A0 A1 -> .plabel= a p
%D    A1 A2 -> .plabel= a ⊤
%D    A1 A3 -> .plabel= l ⊤
%D    A2 A3 -> .plabel= r j
%D
%D    B0 B1 -> .plabel= a p
%D    B1 B2 -> .plabel= a j
%D    B1 B3 -> .plabel= l j
%D    B2 B3 -> .plabel= r j
%D
%D    C0 C1 -> .plabel= a 〈p,q〉
%D    C1 C2 -> .plabel= a ∧
%D    C1 C3 -> .plabel= l j×j
%D    C2 C4 -> .plabel= r j
%D    C3 C4 -> .plabel= a ∧
%D ))
%D enddiagram
%D
$$\pu
  \diag{topologies-2}
$$
are:

%D diagram topologies-3a
%D 2Dx     100  +25  +20 +20  +25  +25
%D 2D  100 A0 - A1 - A2  B0 - B1 - B2 
%D 2D             \  |          \   | 
%D 2D  +18         \ A3'         \ B3'
%D 2D   +7           A3            B3 
%D 2D
%D ren A0 A1 A2 A3' A3 ==> a P(a) ⊤ \;⊤^* ⊤
%D ren B0 B1 B2 B3' B3 ==> a P(a) P(a)^* \;P(a)^{**} P(a)^*
%D
%D (( A0 A1  |-> # .plabel= a p
%D    A1 A2  |-> # .plabel= a ⊤
%D    A1 A3  |-> # .plabel= l ⊤
%D    A2 A3' |-> # .plabel= r j
%D
%D    B0 B1  |-> # .plabel= a p
%D    B1 B2  |-> # .plabel= a j
%D    B1 B3  |-> # .plabel= l j
%D                 .curve= _10pt
%D    B2 B3' |-> # .plabel= r j
%D ))
%D enddiagram
%D
%D diagram topologies-3b
%D 2Dx     100  +40  +55
%D 2D  100 C0 - C1 - C2
%D 2D           |     |
%D 2D  +18      |    C4'
%D 2D   +7      C3 - C4
%D 2D
%D ren C0 C1 C2        ==> a (P(a),Q(a)) P(a){∧}Q(a)
%D ren       C4' C3 C4 ==> (P(a){∧}Q(a))^* (P(a)^*,Q(a)^*) P(a)^*{∧}Q(a)^*
%D
%D (( C0 C1  |-> # .plabel= a 〈p,q〉
%D    C1 C2  |-> # .plabel= a ∧
%D    C1 C3  |-> # .plabel= l j×j
%D    C2 C4' |-> # .plabel= r j
%D    C3 C4  |-> # .plabel= a ∧
%D ))
%D enddiagram
%D
$$\pu
  \begin{array}{c}
    \diag{topologies-3a} \\[35pt]
    \diag{topologies-3b} \\
  \end{array}
$$

We are writing $j(P(a))$ as $P(a)^*$ to suggest a connection between
topologies and the J-operators of \cite{PH2}; we will develop this
idea in section .... .

\newpage

% «top-to-clop»  (to ".top-to-clop")
% (cltp 10 "top-to-clop")
% (clta    "top-to-clop")
\subsection{Topologies induce closure operators}
\label{top-to-clop}

\Theoremsubsection
\label{thm:top-to-clop}
Let $\catE$ be a topos with inclusions, and let $j$ be a topology on
it. For each inclusion $d: D \ito E$ let $\ovl{d}: \ovl{D} \ito E$ be
the inclusion that is classified by $j∘χ_d$, as in the diagram below:
%
%D diagram j-to-clo-def
%D 2Dx     100  +20  +20   +20  +20
%D 2D  100 A0 ------ A1
%D 2D      |          |
%D 2D  +20 |    A2 ---|-------> A3
%D 2D      v  /       v       /
%D 2D  +20 A4 -----> A5 -> A6
%D 2D
%D ren A0 A1 ==> D 1
%D ren A2 A3 ==> \ovl{D} 1
%D ren A4 A5 A6 ==> E Ω Ω
%D
%D (( A0 A1 ->
%D    A0 A4 >-> .PLABEL= _(0.30) d
%D    A1 A5 >->
%D    A2 A3 ->
%D    A2 A4 >-> .PLABEL= ^(0.30) \ovl{d}
%D    A3 A6 >->
%D    A4 A5 -> .plabel= b χ_d
%D    A5 A6 -> .plabel= b j
%D
%D ))
%D enddiagram
%D
$$\pu
  \diag{j-to-clo-def}
$$

Then this operation $d \mapsto \ovl{d}$ is a closure operator ---
i.e., it obeys C1, C2, C3, C4, C5.

\msk

{\bf Proof.}

% «top-to-clop-C1»  (to ".top-to-clop-C1")
% (cltp 10 "top-to-clop-C1")
% (clta    "top-to-clop-C1")

For C1, rename the second 1 to $1'$ in the diagram above and draw the
identity map $1 → 1'$. The slanted rectangle with $\ovl{D}$ in its
upper left corner is a pullback. We can factor the maps $d: D \monicto
E$ and $!: D → 1'$ through it,
%
%D diagram j-to-clo-C1
%D 2Dx     100  +20  +20   +20  +20
%D 2D  100 A0 ------ A1
%D 2D      |          |
%D 2D  +20 |    A2 ---|-------> A3
%D 2D      v  /       v       /
%D 2D  +20 A4 -----> A5 -> A6
%D 2D
%D ren A0 A1 ==> D 1
%D ren A2 A3 ==> \ovl{D} 1'
%D ren A4 A5 A6 ==> E Ω Ω
%D
%D (( A0 A1 ->
%D    A0 A4 >-> .PLABEL= _(0.30) d
%D    A1 A5 >->
%D    A2 A3 ->
%D    A2 A4 >-> .PLABEL= ^(0.30) \ovl{d}
%D    A3 A6 >->
%D    A4 A5 -> .plabel= b χ_d
%D    A5 A6 -> .plabel= b j
%D
%D    A0 A2 >-> .plabel= a m
%D    A1 A3 ->
%D ))
%D enddiagram
%D
$$\pu
  \diag{j-to-clo-C1}
$$
and this gives us a mediating map $m: D \monicto \ovl{D}$. It is easy
to check that this $m$ is a monic and an inclusion.\footnote{I thank
  David Michael Roberts for helping me with this.}

\msk

% «top-to-clop-C2»  (to ".top-to-clop-C2")
% (cltp 10 "top-to-clop-C2")
% (clta    "top-to-clop-C2")

For C2, draw the diagram below:
%
%D diagram j-to-clo-C2
%D 2Dx     100 +20 +15 +10 +20 +20 +35
%D 2D  100 A0 -------- A1
%D 2D  +15     A2 ------------ A3
%D 2D  +20         A4 ------------ A5
%D 2D
%D 2D  +15 A6 -------- A7  A8  A9
%D 2D
%D ren A0 A1 ==> D 1
%D ren A2 A3 ==> \ovl{D} 1
%D ren A4 A5 ==> \ovl{\ovl{D}} 1
%D ren A6 A7 A8 A9 ==> E Ω Ω Ω
%D
%D (( A0 A1 -> A0 A6 >-> .plabel= l d             A1 A7 >->
%D    A2 A3 -> A2 A6 >-> .plabel= l \ovl{d}       A3 A8 >->
%D    A4 A5 -> A4 A6 >-> .plabel= l \ovl{\ovl{d}} A5 A9 >->
%D    A0 A2 >-> A1 A3 >->
%D    A2 A4 >-> A3 A5 >->
%D    A6 A7 -> .plabel= b χ_d
%D    A7 A8 -> .plabel= b j
%D    A8 A9 -> .plabel= b j
%D ))
%D enddiagram
%D
$$\pu
  \diag{j-to-clo-C2}
$$

The inclusion $\ovl{d}$ is classified by $j∘χ_d$ and $\ovl{\ovl{d}}$
is classified by $j∘j∘χ_d$. By LT1 we have $j∘j=j$, and so $j∘χ_d =
j∘j∘χ_d$. This means that $\ovl{d}$ and $\ovl{\ovl{d}}$ are two
inclusions classified by the same map --- so $\ovl{d} =
\ovl{\ovl{d}}$, and the inclusion $\ovl{D} \ito \ovl{\ovl{D}}$ is the
identity.

\msk

% «top-to-clop-C4»  (to ".top-to-clop-C4")
% (cltp 11 "top-to-clop-C4")
% (clta    "top-to-clop-C4")

To prove C4 we use the diagram below and the series of equalities at
the right of it:
%
%D diagram j-to-clo-C4
%D 2Dx     100 +30 +25
%D 2D  100 C0  C1  C2
%D 2D
%D 2D  +25     C3  C4
%D 2D
%D ren C0 C1 C2 C3 C4 ==> E Ω{×}Ω Ω Ω{×}Ω Ω
%D
%D (( C0 C1 -> .plabel= a 〈χ_c,χ_d〉
%D    C1 C2 -> .plabel= a ∧
%D    C1 C3 -> .plabel= l j×j
%D    C2 C4 -> .plabel= r j
%D    C3 C4 -> .plabel= a ∧
%D ))
%D enddiagram
%D
$$\pu
  \diag{j-to-clo-C4}
  \quad
  \begin{array}{rcl}
    χ_{(\ovl{c∩d})} &=& j∘χ_{c∩d} \\
                    &=& j∘∧∘〈χ_c,χ_d〉 \\
                    &=& ∧∘(j×j)∘〈χ_c,χ_d〉 \\
                    &=& ∧∘〈j∘χ_c,j∘χ_d〉 \\
                    &=& ∧∘〈χ_{\ovl{c}},χ_{\ovl{d}}〉 \\
                    &=& χ_{(\ovl{c}∩\ovl{d})} \\
  \end{array}
$$

The inclusions $\ovl{c∩d}$ and $\ovl{c}∩\ovl{d}$ are classified by the
same map, so they are equal.

\msk

% «top-to-clop-C3»  (to ".top-to-clop-C3")
% (cltp 11 "top-to-clop-C3")
% (clta    "top-to-clop-C3")

The proof of C3 is this series of inferences:
%:
%:  c⊆d
%:  -----
%:  c=c∧d
%:  -----------------  -------------------------
%:  \ovl{c}=\ovl{c∧d}  \ovl{c∧d}=\ovl{c}∧\ovl{d}
%:  --------------------------------------------
%:         \ovl{c}=\ovl{c}∧\ovl{d}
%:         -----------------------
%:           \ovl{c}⊆\ovl{d}
%:
%:           ^j-to-clo-C3
%:
\pu
$$\ded{j-to-clo-C3}$$

% «top-to-clop-C5»  (to ".top-to-clop-C5")
% (cltp 11 "top-to-clop-C5")
% (clta    "top-to-clop-C5")

The proof of C5 is this diagram
%
%D diagram j-to-clo-C5
%D 2Dx     100  +20  +20   +20 +20  +20  +20
%D 2D  100 A0 ------ A1 ------ A2 _____
%D 2D      |  \      |  \      |       \
%D 2D  +20 |    A3 --|--- A4 --|-------- A5
%D 2D      |  /      |  /      |       /
%D 2D  +20 A6 ------ A7 ------ A8 - A9
%D 2D
%D ren A0  A1  A2      ==> f^{-1}(D) D 1
%D ren   A3  A4     A5 ==> \ovl{f^{-1}(D)}{=}f^{-1}(\ovl{D}) \ovl{D} 1
%D ren   A3  A4     A5 ==> \ovl{f^{-1}(D)}                   \ovl{D} 1
%D ren A6  A7  A8 A9   ==> B E Ω Ω
%D
%D (( # Horizontal arrows:
%D    A0 A1 ->
%D    A1 A2 ->
%D    A3 A4 ->
%D    A4 A5 ->
%D    A6 A7 -> .plabel= b f
%D    A7 A8 -> .plabel= b χ_d
%D    A8 A9 -> .plabel= b j
%D
%D    # Vertical arrows:
%D    A0 A6 >-> .PLABEL= _(0.30) f^{-1}(d)
%D    A1 A7 >-> .PLABEL= _(0.30) d
%D    A2 A8 >->
%D
%D    # Diagonal arrows:
%D    A0 A3 >->
%D    A1 A4 >->
%D    A2 A5  ->
%D    A3 A6 >-> .PLABEL= ^(0.20) \ovl{f^{-1}(d)} # \sm{\ovl{f^{-1}(s)}=\\f^{-1}(\ovl{s})}
%D    A4 A7 >-> .plabel= r \ovl{d}
%D    A5 A9 >->
%D
%D ))
%D enddiagram
%D
$$\pu
  \diag{j-to-clo-C5}
$$
plus these equalities:
$$\begin{array}{rcl}
    χ_{(\ovl{f^{-1}(d)})} &=& j∘χ_{f^{-1}(d)} \\
                          &=& j∘χ_d∘f \\
                          &=& χ_{\ovl{d}}∘f \\
                          &=& χ_{f^{-1}(\ovl{d})} \\
  \end{array}
$$
The inclusions $f^{-1}(d)$ and $f^{-1}(\ovl{d})$ are classified by the
same map, so they are the same inclusion.

% «restricting-a-clop»  (to ".restricting-a-clop")
% (cltp 13 "restricting-a-clop")
% (clta    "restricting-a-clop")

\subsection{Restricting a $\clop_E$}

In this section we will see how a closure operation $\clop_E$ can be
``restricted'' to a subset $D⊆E$.

\msk

\Theoremsubsection
\label{thm:restr-clop-1}
Let $\catE$ be a topos with inclusions, with a
closure operator $\clop$. If $C⊆D⊆E$ in it, then the closure of $m: C
\ito D$ can be calculated from the closures of $c: C \ito D$ and $d: D
\ito E$ --- and we have $\ovl{m} = d^{-1}(\ovl{c})$ and $\dom(m) =
\ovl{C}∩D$.

\ssk

{\bf Proof.} draw the diagram at the left below, that is the diagram
for C5 with some things renamed. The pullback of $c: C \ito E$ and $d:
D \ito E$ is $C∩D$, which is $C$; so $\ovl{d}^{-1}(c) = m$, and we
have the diagram at the right.
%
%D diagram restr-clop-1
%D 2Dx     100  +20 +20  +20 +20  +20 +20  +20 
%D 2D  100 A0 ----- B0       C0 ----- D0       
%D 2D      |  \     |  \     |  \     |  \     
%D 2D  +20 |    A1 -|--- B1  |    C1 -|--- D1  
%D 2D      |  /     |  /     |  /     |  /     
%D 2D  +20 A2 ----- B2       C2 ----- D2       
%D 2D
%D ren A0 A1 A2 ==> d^{-1}(C) \ovl{d^{-1}(C)} D
%D ren B0 B1 B2 ==> C \ovl{C} E
%D ren C0 C1 C2 ==> C C^D{=}\ovl{C}{∩}D D
%D ren D0 D1 D2 ==> C \ovl{C} E
%D
%D (( A0 B0 >->
%D    A1 B1 >->
%D    A2 B2 >-> .plabel= b d
%D    
%D    A0 A2 >-> .PLABEL= _(0.25) d^{-1}(c)
%D    A0 A1 >->
%D    A1 A2 >-> .PLABEL= ^(0.30) \ovl{d^{-1}(c)}
%D    
%D    B0 B2 >-> .PLABEL= _(0.25) c
%D    B0 B1 >->
%D    B1 B2 >-> .PLABEL= ^(0.30) \ovl{c}
%D ))
%D (( C0 D0 >->
%D    C1 D1 >->
%D    C2 D2 >-> .plabel= b d
%D    
%D    C0 C2 >-> .PLABEL= _(0.25) m
%D    C0 C1 >->
%D    C1 C2 >-> .PLABEL= ^(0.30) \!\!\ovl{m}=d^{-1}(\ovl{c})
%D    
%D    D0 D2 >-> .PLABEL= _(0.25) c
%D    D0 D1 >->
%D    D1 D2 >-> .PLABEL= ^(0.30) \ovl{c}
%D ))
%D enddiagram
%D
$$\pu
  \diag{restr-clop-1}
$$

Our notation for the domain of the closure of an $m: C \ito D$ when
the name $\ovl{C}$ is taken will be $C^D$, for ``the closure of $C$ in
$D$''; the operation `$·^D$' will generalize the `$·^*$' of
\cite{PH2}. As $C^D$ is the pullback of $\ovl{c}: \ovl{C} \ito E$ and
$d: D \ito E$ we have $C^D = \ovl{C}∩D = C^E ∩ D$.

\msk

\Theoremsubsection
\label{thm:restr-clop-2}
Let $\catE$ be a topos with inclusions with closure
operator $\clop$. If $D⊆E$ in $\catE$, then $\clop_D$ can be obtained
from $\clop_E$ in the following way:
$$\begin{array}{rrcl}
  \clop_D: &     \Incs(D) &→& \Incs(D) \\
           & (m:C \ito D) &→& (\ovl{m}: C^D → D) \\
           &              & & := (d^{-1}(\ovl{c}): \ovl{C}∩D \ito D) \\
  \end{array}
$$
where $c$ is $c:C \ito E$ and $\ovl{c}$ is its closure, $\ovl{c}:
\ovl{C} \ito E$.

\ssk

{\bf Proof.} This is an easy corollary of Theorem
\ref{thm:restr-clop-1}.

\newpage

% «dense-and-closed»  (to ".dense-and-closed")
% (cltp 13 "dense-and-closed")
% (clta    "dense-and-closed")

\subsection{Dense and closed}

For the next theorems we need some definitions:

An inclusion $c:C → D$ is {\sl dense} iff $\ovl{c} = \id_D$.

An inclusion $d:D → E$ is {\sl closed} iff $\ovl{d} = d$.

\ssk

\def\TE#1{\text{ #1}}

\Theoremsubsection
\label{thm:dense-and-closed-1}
If an inclusion $a: A \ito B$ is dense and closed then it is the
identity.

{\bf Proof:}
%:
%:   a\TE{dense}   a\TE{closed}
%:  -------------  ------------
%:  \ovl{a}=\id_B  \ovl{a}=a
%:  ------------------------
%:      a=\id_B
%:
%:      ^dense-and-closed
%:
\pu
$$\ded{dense-and-closed}
$$

\msk

\Theoremsubsection
\label{thm:dense-and-closed-2}
In a topos with inclusions $\catE$ with closure operator $\clop$, for
any inclusion $d:D \ito E$ we have:
%
%D diagram dense-and-closed-1
%D 2Dx     100  +20
%D 2D  100 A0
%D 2D      |  \
%D 2D  +20 |    A1
%D 2D      |  /
%D 2D  +20 A2
%D 2D
%D ren A0 A1 A2 ==> D \ovl{D} E
%D
%D (( A0 A1 >-> .plabel= r m\TE{(dense)}
%D    A1 A2 >-> .plabel= r \ovl{d}\TE{(closed)}
%D    A0 A2 >-> .plabel= l d
%D ))
%D enddiagram
%D
$$\pu
  \diag{dense-and-closed-1}
$$

{\bf Proof.} $\ovl{\ovl{d}} = \ovl{d}$, so $\ovl{d}$ is closed. To see
that $m: D \ito \ovl{D}$ is dense, we build the diagram at the left
below:
%
%D diagram dense-and-closed-2
%D 2Dx     100  +20 +20  +20 +20  +20 +20  +20 
%D 2D  100 A0 ----- B0       C0 ----- D0       
%D 2D      |  \     |  \     |  \     |  \     
%D 2D  +20 |    A1 -|--- B1  |    C1 -|--- D1  
%D 2D      |  /     |  /     |  /     |  /     
%D 2D  +20 A2 ----- B2       C2 ----- D2       
%D 2D
%D ren A0 A1 A2 ==> \ovl{d}^{-1}(D) \ovl{\ovl{d}^{-1}(\ovl{D})} \ovl{D}
%D ren B0 B1 B2 ==> D \ovl{D} E
%D ren C0 C1 C2 ==> D \ovl{D} \ovl{D}
%D ren D0 D1 D2 ==> D \ovl{D} E
%D
%D (( A0 B0 >->
%D    A1 B1 >->
%D    A2 B2 >-> .plabel= b \ovl{d}
%D    
%D    A0 A2 >-> .PLABEL= _(0.25) \ovl{d}^{-1}(d)
%D    A0 A1 >->
%D    A1 A2 >-> .PLABEL= ^(0.30) \ovl{\ovl{d}^{-1}(d)}
%D    
%D    B0 B2 >-> .PLABEL= _(0.25) d
%D    B0 B1 >->
%D    B1 B2 >-> .PLABEL= ^(0.30) \ovl{d}
%D ))
%D (( C0 D0 ->
%D    C1 D1 ->
%D    C2 D2 -> .plabel= b d
%D    
%D    C0 C2 >-> .PLABEL= _(0.25) m
%D    C0 C1 >->
%D    C1 C2 >-> .PLABEL= ^(0.30) \ovl{m}=\id
%D    
%D    D0 D2 >-> .PLABEL= _(0.25) d
%D    D0 D1 >->
%D    D1 D2 >-> .PLABEL= ^(0.30) \ovl{d}
%D ))
%D enddiagram
%D
$$\pu
  \diag{dense-and-closed-2}
$$
we have $\ovl{d}^{-1}(D) = \ovl{D}∩D = D$ and $\ovl{d}^{-1}(\ovl{D}) =
\ovl{D}∩\ovl{D} = \ovl{D}$, so we can rewrite it as the diagram at the
right above, and we get that $\ovl{m} = \id$.

\newpage

% «clop-to-top»  (to ".clop-to-top")
% (cltp 14 "clop-to-top")
% (clta    "clop-to-top")

\subsection{Closure operators induce topologies}

Let $\catE$ be a topos with inclusions, and $\clop$ a
closure operator on it. Build this diagram on it:
%
%D diagram def-J
%D 2Dx     100  +20  +20  +20
%D 2D  100 A0
%D 2D      |  \
%D 2D  +20 |    A1 ------ A2
%D 2D      |  /         /
%D 2D  +20 A3 ------ A4
%D 2D
%D ren A0 A1 A2 A3 A4 ==> 1 J{:=}\ovl{1} 1 Ω Ω
%D
%D (( A0 A1 >->
%D    A0 A3 >-> .plabel= l ⊤
%D    A1 A2  ->
%D    A1 A3 >-> .plabel= r \ovl{⊤}
%D    A2 A4 >->
%D    A3 A4 -> .plabel= b j:=χ_{\ovl{⊤}}
%D ))
%D enddiagram
%D
$$\pu
  \diag{def-J}
$$

Here the closure of $⊤:1→Ω$ is $\ovl{⊤}:\ovl{1}→Ω$, and $J$ is an
alternate name for this $\ovl{1}$; and $j:=χ_{\ovl{⊤}}$ is the map
that classifies $\ovl{⊤}$.

\msk

% «clop-to-j»  (to ".clop-to-j")
% (cltp 14 "clop-to-j")
% (clta    "clop-to-j")

\Theoremsubsection
\label{thm:clop-to-top-1}
For every inclusion $d: D \ito E$ we have $χ_{\ovl{d}} = j∘χ_d$, where
$j$ is the map above.

\ssk

{\bf Proof.} Take a map $f:E→Ω$, and add to the diagram above the
diagram for $\ovl{f^{-1}(⊤)} = f^{-1}(\ovl{⊤})$. We get this:
%
%D diagram chi(ovl(d))-1
%D 2Dx     100  +20  +20  +20  +20  +20
%D 2D  100 A0 ------ A1
%D 2D      |  \      |  \
%D 2D  +20 |    A2 --|--- A3 ------ A4
%D 2D      |  /      |  /         /
%D 2D  +20 A5 ------ A6 ------ A7
%D 2D
%D ren A0 A1    ==> f^{-1}(1) 1
%D ren A2 A3 A4 ==> \ovl{f^{-1}(1)} J 1
%D ren A5 A6 A7 ==> E Ω Ω
%D
%D (( A0 A1 ->
%D    A2 A3 ->
%D    A3 A4 ->
%D    A5 A6 -> .plabel= b f
%D    A6 A7 -> .plabel= b j
%D
%D    A0 A5 >-> .PLABEL= _(0.30) f^{-1}(⊤)
%D    A0 A2 >->
%D    A2 A5 >-> .PLABEL= ^(0.30) \ovl{f^{-1}(⊤)}
%D
%D    A1 A6 >-> .PLABEL= _(0.30) ⊤
%D    A1 A3 >->
%D    A3 A6 >-> .PLABEL= ^(0.30) \ovl{⊤}
%D
%D    A4 A7 >-> .PLABEL= ^(0.30) ⊤
%D ))
%D enddiagram
%D
$$\pu
  \diag{chi(ovl(d))-1}
$$

This map $f$ is the classifying map for some inclusion; let's call it
$d:D \ito E$, and rewrite $f$ as $χ_d$. We get:
%
%D diagram chi(ovl(d))-2
%D 2Dx     100  +20  +20  +20  +20  +20
%D 2D  100 A0 ------ A1
%D 2D      |  \      |  \
%D 2D  +20 |    A2 --|--- A3 ------ A4
%D 2D      |  /      |  /         /
%D 2D  +20 A5 ------ A6 ------ A7
%D 2D
%D ren A0 A1    ==> D 1
%D ren A2 A3 A4 ==> \ovl{D} J 1
%D ren A5 A6 A7 ==> E Ω Ω
%D
%D (( A0 A1 ->
%D    A2 A3 ->
%D    A3 A4 ->
%D    A5 A6 -> .plabel= b χ_d
%D    A6 A7 -> .plabel= b j
%D
%D    A0 A5 >-> .PLABEL= _(0.30) d
%D    A0 A2 >->
%D    A2 A5 >-> .PLABEL= ^(0.30) \ovl{d}
%D
%D    A1 A6 >-> .PLABEL= _(0.30) ⊤
%D    A1 A3 >->
%D    A3 A6 >-> .PLABEL= ^(0.30) \ovl{⊤}
%D
%D    A4 A7 >-> .PLABEL= ^(0.30) ⊤
%D ))
%D enddiagram
%D
$$\pu
  \diag{chi(ovl(d))-2}
$$

We have $\ovl{d} = \ovl{f^{-1}(⊤)} = f^{-1}(\ovl{⊤}) =
{χ_d}^{-1}(\ovl{⊤}) = {χ_d}^{-1}(j^{-1}(⊤)) = (j∘χ_d)^{-1}(⊤)$, and so
$χ_{\ovl{d}} = χ_{((j∘χ_d)^{-1}(⊤))} = j∘χ_d$.

\bsk

\newpage

% «clop-to-top-LT1»  (to ".clop-to-top-LT1")
% (cltp 15 "clop-to-top-LT1")
% (clta    "clop-to-top-LT1")

\Theoremsubsection
\label{thm:clop-to-top-2}
The map $j$ defined above is a topology.

{\bf Proof.} To prove LT1 we have to see that $j=j∘j$. We have
$\ovl{d} = \ovl{\ovl{d}}$ for all inclusions $d$; so $χ_{\ovl{d}} =
χ_{\ovl{\ovl{d}}}$ always. We have $χ_{\ovl{d}} = j∘χ_d$ and
$χ_{\ovl{\ovl{d}}} = j∘j∘χ_d$, so $j∘χ_d = j∘j∘χ_d$ always holds.
There is a way to make $χ_d = \id$ here --- which is when $d:D \ito E$
is $⊤:1 \ito Ω$ --- and so a particular case of $j∘χ_d = j∘j∘χ_d$ is
$j∘\id = j∘j∘\id$, which gives us $j=j∘j$.

\msk

% «clop-to-top-LT2»  (to ".clop-to-top-LT2")
% (cltp 15 "clop-to-top-LT2")
% (clta    "clop-to-top-LT2")

To prove LT2 we have to see that $⊤_Ω=j∘⊤_Ω$, i.e., that
$⊤∘!_Ω=j∘⊤∘!_Ω$. To do this we draw this diagram,
%
%D diagram LT2
%D 2Dx     100  +20  +20  +20  +20  +20  +20  +20
%D 2D  100 A0 ------ A1 ------ A2
%D 2D      |  \      |  \      |  \
%D 2D  +20 |    A3 --|--- A4 --|--- A5 ------ A6
%D 2D      |  /      |  /      |  /         /
%D 2D  +20 A7 ------ A8 ------ A9 ------ A10
%D 2D
%D ren A0 A1 A2     ==> Ω 1 1
%D ren A3 A4 A5 A6  ==> Ω 1 J 1
%D ren A7 A8 A9 A10 ==> Ω 1 Ω Ω
%D
%D (( A0 A1 -> A1 A2 ->
%D    A3 A4 -> A4 A5 -> A5 A6 ->
%D    A7 A8 -> .plabel= b !_Ω
%D    A8 A9 -> .plabel= b ⊤
%D    A9 A10 -> .plabel= b j
%D
%D    A0 A7 >-> .PLABEL= _(0.30) \id_Ω
%D    A0 A3 >->
%D    A3 A7 >-> .PLABEL= ^(0.30) \ovl{\id_Ω}=\id_Ω
%D
%D    A1 A8 >-> .PLABEL= _(0.30) !_1
%D    A1 A4 >->
%D    A4 A8 >-> .PLABEL= ^(0.30) \ovl{!_1}=!_1
%D
%D    A2 A9 >-> .PLABEL= _(0.30) ⊤
%D    A2 A5 >->
%D    A5 A9 >-> .PLABEL= ^(0.30) \ovl{⊤}
%D
%D    A6 A10 >-> .PLABEL= ^(0.30) ⊤
%D ))
%D enddiagram
%D
$$\pu
  \diag{LT2}
$$
and check that its two upright squares and its three lower slanted
squared are pullbacks. With this we get that both $⊤∘!_Ω$ and
$j∘⊤∘!_Ω$ classify $\id_Ω$, so $⊤∘!_Ω = j∘⊤∘!_Ω$.

\msk

% «clop-to-top-LT3»  (to ".clop-to-top-LT3")
% (cltp 15 "clop-to-top-LT3")
% (clta    "clop-to-top-LT3")

To prove LT3 we start by choosing any two inclusions with the same
codomain, $c:C \ito E$ and $d:D \ito E$. From the maps $χ_c,χ_d:E→Ω$
we build a map $〈χ_c,χ_d〉:Ω→Ω×Ω$, and we plug it on the diagram for
LT3. We get:
%
%D diagram clo-to-LT3
%D 2Dx     100  +30  +30
%D 2D  100 C0 - C1 - C2
%D 2D           |     |
%D 2D           |     |
%D 2D  +25      C3 - C4
%D 2D
%D ren C0 C1 C2 ==> E Ω{×}Ω Ω
%D ren    C3 C4 ==>   Ω{×}Ω Ω
%D
%D (( C0 C1 -> .plabel= a 〈χ_c,χ_d〉
%D    C1 C2 -> .plabel= a ∧
%D    C1 C3 -> .plabel= l j×j
%D    C2 C4 -> .plabel= r j
%D    C3 C4 -> .plabel= a ∧
%D ))
%D enddiagram
%D
$$\pu
  \diag{clo-to-LT3}
$$

We have
$$\begin{array}{rcl}
  ∧∘(j×j)∘〈χ_c,χ_d〉 &=& χ_{\ovl{c}∩\ovl{d}} \\
  j∘∧∘〈χ_c,χ_d〉 &=& χ_{\ovl{c∩d}} \\
  \end{array}
$$
and C4 tells us that $\ovl{c}∩\ovl{d} = \ovl{c∩d}$; so it is always
true that $∧∘(j×j)∘〈χ_c,χ_d〉 = j∘∧∘〈χ_c,χ_d〉$. We can make $〈χ_c,χ_d〉$
be the identity map if we take $E:=Ω×Ω$, $〈χ_c,χ_d〉 = \id_{Ω×Ω} =
〈π,π'〉$. The internal views of $χ_c$ and $χ_d$ are:
%
%D diagram chi_c-and-chi_d
%D 2Dx     100 +30 +25 +25
%D 2D  100 A0  A1  B0  B1
%D 2D
%D 2D  +25 A2  A3  B2  B3
%D 2D
%D 2D  +20 C0  C1  D0  D1
%D 2D
%D 2D  +25 C2  C3  D2  D3
%D 2D
%D ren A0 A1 A2 A3 ==> C 1 Ω{×}Ω Ω
%D ren C0 C1 C2 C3 ==> D 1 Ω{×}Ω Ω
%D ren B2 B3       ==> (P,Q) P
%D ren D2 D3       ==> (P,Q) Q
%D
%D (( A0 A1 ->
%D    A0 A2 >-> .plabel= l c
%D    A1 A3 >-> .plabel= r ⊤
%D    A2 A3 -> .plabel= b χ_c=π
%D
%D    C0 C1 ->
%D    C0 C2 >-> .plabel= l d
%D    C1 C3 >-> .plabel= r ⊤
%D    C2 C3 -> .plabel= b χ_d=π
%D
%D    B2 B3 |->
%D
%D    D2 D3 |->
%D ))
%D enddiagram
%D
$$\pu
  \diag{chi_c-and-chi_d}
$$

In $\Set$ we can construct the subsets $C$ and $D$ as:
$$\begin{array}{rcl}
  C &=& \setofst{(P,Q)∈Ω×Ω}{P=⊤} \\
    &=& \{⊤\}×Ω \\
  D &=& \setofst{(P,Q)∈Ω×Ω}{Q=⊤} \\
    &=& Ω×\{⊤\} \\
  \end{array}
$$

This {\sl suggests} that we can generalize that construction to any
topos as:
%
%D diagram LT3-last-construction
%D 2Dx     100 +25 +30 +25
%D 2D  100 A0  A1  B0  B1
%D 2D
%D 2D  +25 A2  A3  B2  B3
%D 2D
%D 2D  +20 C0  C1  D0  D1
%D 2D
%D 2D  +25 C2  C3  D2  D3
%D 2D
%D ren A0 A1 A2 A3 ==> C 1 E Ω
%D ren B0 B1 B2 B3 ==> 1{×}Ω 1 Ω{×}Ω Ω
%D ren C0 C1 C2 C3 ==> D 1 E Ω
%D ren D0 D1 D2 D3 ==> Ω{×}1 1 Ω{×}Ω Ω
%D
%D (( A0 A1 ->
%D    A0 A2 >-> .plabel= l c
%D    A1 A3 >-> .plabel= r ⊤
%D    A2 A3 -> .plabel= b χ_c
%D
%D    B0 B1 ->
%D    B0 B2 >-> .plabel= l (⊤×\id)
%D    B1 B3 >-> .plabel= r ⊤
%D    B2 B3 -> .plabel= b χ_c=π
%D
%D    C0 C1 ->
%D    C0 C2 >-> .plabel= l d
%D    C1 C3 >-> .plabel= r ⊤
%D    C2 C3 -> .plabel= b χ_d
%D
%D    D0 D1 ->
%D    D0 D2 >-> .plabel= l (\id×⊤)
%D    D1 D3 >-> .plabel= r ⊤
%D    D2 D3 -> .plabel= b χ_d=π
%D ))
%D enddiagram
%D
$$\pu
  \diag{LT3-last-construction}
$$

These constructions do work, but I will skip the details of the proof.
So: with $c=(⊤×\id)$ and $d=(\id×⊤)$ we have $〈χ_c,χ_d〉=\id_{Ω×Ω}$,
and in this particular case our equality $∧∘(j×j)∘〈χ_c,χ_d〉 =
j∘∧∘〈χ_c,χ_d〉$ reduces to $∧∘(j×j)=j∘∧$ --- and this proves LT3.

\newpage

%  ____  _  _           _   _             
% | __ )(_)(_) ___  ___| |_(_) ___  _ __  
% |  _ \| || |/ _ \/ __| __| |/ _ \| '_ \ 
% | |_) | || |  __/ (__| |_| | (_) | | | |
% |____/|_|/ |\___|\___|\__|_|\___/|_| |_|
%        |__/                             
%
% «clop-top-bij»  (to ".clop-top-bij")
% (cltp 17 "clop-top-bij")
% (clta    "clop-top-bij")
\subsection{A bijection}
\label{clop-top-bij}

% (favp 28 "basic-example-bij")
% (fav     "basic-example-bij")

We saw that a closure operator induces a topology and that a topology
induces a closure operator. Now we need to check that these two
operations, that we can abbreviate as $\clop \mapsto j$ and $j \mapsto
\clop$, as below,
%
%D diagram clop-top-bij-1
%D 2Dx     100    +50
%D 2D  100 A0 |-> A1
%D 2D   +8 A3 <-| A2
%D 2D
%D ren A0 A1 ==> \clop  j
%D ren A3 A2 ==> \clop  j
%D
%D (( A0 A1 |-> .plabel= a j:=χ_{\ovl{⊤}}
%D    A3 A2 <-| .plabel= b \clop:=(λd.σ(j∘χ_d))
%D ))
%D enddiagram
%D
$$\pu
  \diag{clop-top-bij-1}
$$
are inverses to one another --- i.e., that the composites
$\clop \mapsto j \mapsto \clop$ and $j \mapsto \clop \mapsto j$ are
identity maps. We will organize all this visually as:
%
%
%D diagram clop-top-bij-2
%D 2Dx     100    +45   +50    +55
%D 2D  100 A0'                 B1'
%D 2D  +15 A0 |-> A1    B0 |-> B1 
%D 2D   +8 A3 <-| A2    B3 <-| B2 
%D 2D
%D ren A0'   ==> (λd.\ovl{d})
%D ren A0 A1 ==> \clop                   χ_{\ovl{⊤}}
%D ren A3 A2 ==> (λd.σ(χ_{\ovl{⊤}}∘χ_d)) χ_{\ovl{⊤}}
%D 
%D ren    B1' ==>               j∘χ_{⊤}
%D ren B0 B1  ==> (λd.σ(j∘χ_d)) χ_{((λd.σ(j∘χ_d))(⊤))}
%D ren B3 B2  ==> (λd.σ(j∘χ_d)) j
%D
%D (( A0' A0 =
%D    A0 A1 |->
%D    A3 A2 <-|
%D ))
%D (( B1' B1 =
%D    B0 B1 |->
%D    B3 B2 <-|
%D ))
%D enddiagram
%D
$$\pu
  \diag{clop-top-bij-2}
$$

To prove that $\clop \mapsto j \mapsto \clop$ is the identity we need
to check that in any topos with inclusions with a closure operator
$\clop$ we have that $\clop$, i.e., $(λd.\ovl{d})$, is equal to
$(λd.σ(χ_{\ovl{⊤}}∘χ_d))$. It is enough that check that we have
$\ovl{d} = σ(χ_{\ovl{⊤}}∘χ_d)$ for any inclusion $d$. Look at the
diagram below...

%D diagram clop-j-clop
%D 2Dx     100  +20  +20  +20  +20  +20
%D 2D  100 A0 ------ A1
%D 2D      |  \      |  \
%D 2D  +20 |    A2 --|--- A3 ------ A4
%D 2D      |  /      |  /         /
%D 2D  +20 A5 ------ A6 ------ A7
%D 2D
%D ren A0 A1    ==> D 1
%D ren A2 A3 A4 ==> \ovl{D} J 1
%D ren A5 A6 A7 ==> E Ω Ω
%D
%D (( A0 A1 ->
%D    A2 A3 ->
%D    A3 A4 ->
%D    A5 A6 -> .plabel= b χ_d
%D    A6 A7 -> .plabel= b χ_{\ovl{⊤}}
%D
%D    A0 A5 >-> .PLABEL= _(0.30) d
%D    A0 A2 >->
%D    A2 A5 >-> .PLABEL= ^(0.30) \ovl{d}
%D
%D    A1 A6 >-> .PLABEL= _(0.30) ⊤
%D    A1 A3 >->
%D    A3 A6 >-> .PLABEL= ^(0.30) \ovl{⊤}
%D
%D    A4 A7 >-> .PLABEL= ^(0.30) ⊤
%D ))
%D enddiagram
%D
$$\pu
  \diag{clop-j-clop}
$$

To prove that $j \mapsto \clop \mapsto j$ is the identity we need to
check that in any topos with inclusions with a topology $j$ we have
$j=j∘χ_T$. Look at the diagram below:
%
%D diagram j-clop-j
%D 2Dx     100 +25 +25
%D 2D  100 A0  A1
%D 2D
%D 2D  +25 A2  A3  A4
%D 2D
%D ren A0 A1 A2 A3 A4 ==> 1 1 Ω Ω Ω
%D
%D (( A0 A1 ->
%D    A0 A2 >-> .plabel= l ⊤
%D    A1 A3 >-> .plabel= r ⊤
%D    A2 A3 -> .plabel= l χ_⊤
%D    A3 A4 -> .plabel= b j
%D ))
%D enddiagram
%D
$$\pu
  \diag{j-clop-j}
$$
We have $χ_⊤=\id_Ω$, and so $j∘χ_T = j∘\id = j$.

\newpage

% __        ___ _   _                 _     _                
% \ \      / (_) |_| |__   ___  _   _| |_  (_)_ __   ___ ___ 
%  \ \ /\ / /| | __| '_ \ / _ \| | | | __| | | '_ \ / __/ __|
%   \ V  V / | | |_| | | | (_) | |_| | |_  | | | | | (__\__ \
%    \_/\_/  |_|\__|_| |_|\___/ \__,_|\__| |_|_| |_|\___|___/
%                                                            
% «without-inclusions»  (to ".without-inclusions")
% (cltp 17 "without-inclusions")
% (clta    "without-inclusions")
\section{Translating all this to toposes without inclusions}
\label{without-inclusions}

Let's start by an example -- we will translate Theorem
% (cltp 12 "restricting-a-clop")
% (clta    "restricting-a-clop")
\ref{thm:restr-clop-1}.
This diagram condenses the two diagrams of the original proof into a
single one:
%
%D diagram before-translation-1
%D 2Dx     100  +70 +60  +20 
%D 2D  100 A0 ----- B0       
%D 2D      |  \     |  \     
%D 2D  +20 |    A1 -|--- B1  
%D 2D      |  /     |  /     
%D 2D  +20 A2 ----- B2       
%D 2D
%D ren A0    ==> d^{-1}(C){=}D{∩}C{=}C
%D ren A1 A2 ==> C^D{=}\ovl{d¹(C)}{=}d¹(\ovl{C}){=}D{∩}\ovl{C} D
%D ren B0 B1 B2 ==> C \ovl{C} E
%D
%D (( A0 B0 `->
%D    A1 B1 `->
%D    A2 B2 `-> .plabel= b d
%D    
%D  # A0 A2 `-> .PLABEL= _(0.25) \sm{d¹(c)\\=(d∩c)'\\m}
%D    A0 A2 `-> .plabel= l       \sm{d¹(c)\\=(d∩c)'\\m}
%D    A0 A1 `->
%D    A1 A2 `-> .PLABEL= ^(0.40) \ph{mmmmm}\ovl{m}=\ovl{d¹(c)}=d¹(\ovl{c})=(d∩\ovl{c})'
%D    
%D    B0 B2 `-> .PLABEL= _(0.25) c
%D    B0 B1 `->
%D    B1 B2 `-> .PLABEL= ^(0.30) \ovl{c}
%D ))
%D enddiagram
%D
$$\pu
  \diag{before-translation-1}
$$
where $(d∩c)'$ and $(d∩\ovl{c})'$ are maps in intersection pullbacks.
The convention is that if $a:A \monicto C$ and $b: B \monicto C$ are
monics then the components of the diagram for $A∩B$ are named like
this:
%
%D diagram intersection-pullback
%D 2Dx     100  +30
%D 2D  100 A0 - A1
%D 2D
%D 2D  +25 A2 - A3
%D 2D
%D ren A0 A1 A2 A3 ==> A∩B B A C
%D
%D (( A0 A1 >-> .plabel= a (a∩b)''
%D    A0 A2 >-> .plabel= l (a∩b)'
%D    A0 A3 >-> .plabel= m (a∩b)
%D    A1 A3 >-> .plabel= r b
%D    A2 A3 >-> .plabel= b a
%D ))
%D enddiagram
%D
$$\pu
  \diag{intersection-pullback}
$$

This is how I would start to structure the proof above to implement it
in a proof assistant. Most nodes in this tree
%:
%:                      C⊆D
%:                     -----
%:                     C=C∩D
%:                     -----
%:    d\isinclusion    C=D∩C
%:    -------------    -----
%:    d¹(C)=D∩C        D∩C=C
%:    ----------------------
%:       C=d¹(C)                                      d\isinclusion
%:     ---------------   -----------------------    ---------------------
%:     C^D=\ovl{d¹(C)}   \ovl{d¹(C)}=d¹(\ovl{C})    d¹(\ovl{C})=D∩\ovl{C}
%:     ------------------------------------------------------------------
%:                           C^D=D∩\ovl{C}
%:
%:                           ^without-incs-1
%:
%:
\pu
\def\isinclusion{\text{ inclusion}}
$$\ded{without-incs-1}$$
state that two inclusions with different constructions are isomorphic,
and so they are the same morphism. For example, ``$C=d¹(C)$'' is an
abbreviation for this:
$$(m:C \ito D) = (d¹(c): D¹(C) \ito D)$$

The properties of inclusions let us omit the codomains and the names
of the arrows in many cases, and write only their domains.

\newpage

We can regard the tree above as a {\sl proof} of this {\sl equality}
of {\sl inclusions} that appears at the root node:
$$(\ovl{m}:C^D \ito D) = ((d∩\ovl{c})': D∩\ovl{C} \ito D)$$

We can translate it to a {\sl construction} of this {\sl isomorphism}
of {\sl monics}:
$$(\ovl{m}:C^D \monicto D) ≡ ((d∩\ovl{c})': D∩\ovl{C} \monicto D)$$

Now the names of the morphisms are primary and the names of the
objects secondary. I prefer write both, otherwise I feel that the
translated tree becomes unreadable. Here is the translation of the
upper left part of the previous tree:
%:
%L addabbrevs("`->", "\\ito ")
%L addabbrevs(">->", "\\monicto ")
%:
%:                                          (m:C>->D)⊆(\id:D>->D)
%:                                        -------------------------
%:                                        (m:C>->D)≡(m∩\id:C∩D>->D)
%:                                        -------------------------
%:                d\ismonic               (m:C>->D)≡(\id∩m:D∩C>->D)
%:    ---------------------------------   -------------------------
%:    (d¹(c):d¹(C)>->D)≡(\id∩m:D∩C>->D)   (\id∩m:D∩C>->D)≡(m:C>->D)
%:    -------------------------------------------------------------
%:          (d¹(c):d¹(C)>->D)≡(m:C>->D)
%:          ---------------------------
%:          (m:C>->D)≡(d¹(c):d¹(C)>->D)
%:     -------------------------------------------------
%:     (\ovl{m}:C^D>->D)≡({\ovl{d¹(c)}:\ovl{d¹(C)}>->D})
%:
%:     ^without-incs-2
%:
%:
\pu
\def\ismonic{\text{ monic}}
$$\scalebox{0.85}{$
  \ded{without-incs-2}
  $}
$$

I tried to draw a diagram with all the morphisms in the tree above
following my usual conventions, and I found the result too messy. But
if we translate the original diagram to this,
%
%D diagram after-translation-1
%D 2Dx     100  +70 +60  +20 
%D 2D  100 A0 ----- B0       
%D 2D      |  \     |  \     
%D 2D  +20 |    A1 -|--- B1  
%D 2D      |  /     |  /     
%D 2D  +20 A2 ----- B2       
%D 2D
%D ren A0    ==> d^{-1}(C){≡}D{∩}C{≡}C
%D ren A1 A2 ==> C^D{≡}\ovl{d¹(C)}{≡}d¹(\ovl{C}){≡}D{∩}\ovl{C} D
%D ren B0 B1 B2 ==> C \ovl{C} E
%D
%D (( A0 B0 >->
%D    A1 B1 >->
%D    A2 B2 >-> .plabel= b d
%D    
%D  # A0 A2 >-> .PLABEL= _(0.25) \sm{d¹(c)\\≡(d∩c)'\\m}
%D    A0 A2 >-> .plabel= l       \sm{d¹(c)\\≡(d∩c)'\\m}
%D    A0 A1 >->
%D    A1 A2 >-> .PLABEL= ^(0.40) \ph{mmmmm}\ovl{m}≡\ovl{d¹(c)}≡d¹(\ovl{c})≡(d∩\ovl{c})'
%D    
%D    B0 B2 >-> .PLABEL= _(0.25) c
%D    B0 B1 >->
%D    B1 B2 >-> .PLABEL= ^(0.30) \ovl{c}
%D ))
%D enddiagram
%D
$$\pu
  \diag{after-translation-1}
$$
and we define in the right way how to interpret the `$≡$'s in it, then
everything works. In
%
%D diagram how-to-interpret-the-isos
%D 2Dx     100
%D 2D  100 A0
%D 2D
%D 2D  +25 A1
%D 2D
%D ren A0 ==> A_1{≡}A_2{≡}A_3
%D ren A1 ==> B_1{≡}B_2{≡}B_3{≡}B_4
%D
%D (( A0 A1 -> .plabel= r f_1{≡}f_2{≡}f_3{≡}f_4{≡}f_5
%D
%D ))
%D enddiagram
%D
$$\pu
  \diag{how-to-interpret-the-isos}
$$
the ``object'' $A_1{≡}A_2{≡}A_3$ means that we have three objects with
known, but unnamed, isos between each one and the next, like this:
$A_1 ↔ A_2 ↔ A_3$, and the ``arrow'' $f_1{≡}f_2{≡}f_3{≡}f_4{≡}f_5$ is
in fact five ``isomorphic'' arrows, and each $f_i$ goes from some
$A_j$ to some $B_k$, but the diagram does not say what are these
`$j$'s and `$k$'s; in this context ``the `$f_i$'s are isomorphic''
means that the diagram made by $A_1 ↔ A_2 ↔ A_3$, $B_1 ↔ B_2 ↔ B_3 ↔
B_4$, and all the `$f_i$'s commutes.

\msk

The translation sketched above works for all constructions and proofs
in sections \ref{clops}--\ref{clop-top-bij}. It may be possible to
characterize the class of constructions and proofs on which it works,
but this is far beyond the scope of these notes.

\newpage

\pu

\makeatletter
\def\LittleNSetArgs#1{\LittleNArgs@#1}
\def\LittleNArgs@#1#2#3#4{%
  \sa{L2}{#1}\sa{R2}{#2}%
  \sa{L1}{#3}\sa{R1}{#4}%
  }
\makeatother

\makeatletter
\def\OLittleNSetArgs#1{\OLittleNArgs@#1}
\def\OLittleNArgs@#1#2#3#4#5#6#7#8{%
             \sa{21}{#1}\sa{22}{#2}%
  \sa{10}{#3}\sa{11}{#4}\sa{12}{#5}%
  \sa{00}{#6}\sa{01}{#7}\sa{02}{#8}%
  }
\makeatother
\def\littlena  #1{{        \LittleNSetArgs{#1}\tcg{LittleNSmall}         }}
\def\littlenap #1{{ \left( \LittleNSetArgs{#1}\tcg{LittleNSmall} \right) }}
\def\littlen   #1{{        \LittleNSetArgs{#1}\tcg{LittleNMedium}        }}
\def\littlenw  #1{{        \LittleNSetArgs{#1}\tcg{LittleN_chi_f}        }}
\def\littlenbig#1{{        \LittleNSetArgs{#1}\tcg{LittleNBig}           }}
\def\olittlen  #1{{       \OLittleNSetArgs{#1}\zha{OLittleNSmall}        }}

\def\littlenbullets {      \littlena{{•}{•}{•}{•}}       }
\def\littlenbulletsp{\left(\littlena{{•}{•}{•}{•}}\right)}

% ----------------------------------------/

%  ____       _   ____      
% / ___|  ___| |_|  _ \ ___ 
% \___ \ / _ \ __| | | / __|
%  ___) |  __/ |_| |_| \__ \
% |____/ \___|\__|____/|___/
%                           
% «SetDs»  (to ".SetDs")
% (cltp 22 "SetDs")
% (clta    "SetDs")
\section{Toposes of the form $\SetD$}
\label{SetDs}

In sec.\ref{inclusions-SetC} we conventioned that $\catC$ would always
denote a small category. In this section we will need many other
conventions like that, especially to define Grothendieck topologies in
sections \ref{SetD-Jcan} and \ref{SetD-J}.

\msk

From here onwards $\catD$ will always denote a finite DAG regarded as
posetal category. We will consider that $\catD$ is downward directed:
if $u,v∈\catD$ then we say that $u$ is {\sl above} $v$, or $v$ is {\sl
  below} $u$, when we have an arrow $u→v$. We will say that $u$ is
{\sl strictly above} $v$ when we have an arrow $u→v$ and $u≠v$. Note
that every $u∈\catD$ is above itself.

In contexts in which we have a category $\catD$ our (default) topos
will be the topos $\catE = \SetD$, with the inclusions given by the
definitions in sec.\ref{inclusions-SetC}. We will use the conventions
from \cite[section 7.12]{FavC} to draw its objects; for example, if
%
% (favp 50 "functors-as-objects")
% (fava    "functors-as-objects")
%
$$\catD = \littlenbulletsp = \tcg{LittleN 2CG}$$
then these are two objects of $\SetD$:
%
%D diagram an-object-A-of-SetLittleN
%D 2Dx     100 +20
%D 2D  100 A0  A1
%D 2D
%D 2D  +20 A2  A3
%D 2D
%D ren A0 A1 A2 A3 ==> ∅ \{4\} \{5\} \{7\}
%D
%D (( A0 A2 ->
%D    A0 A3 ->
%D    A1 A3 ->
%D    
%D ))
%D enddiagram
%D
%D diagram an-object-B-of-SetLittleN
%D 2Dx     100 +25
%D 2D  100 A0  A1
%D 2D
%D 2D  +20 A2  A3
%D 2D
%D ren A0 A1 A2 A3 ==> \{1,2\} \{3,4\} \{5,6\} \{7\}
%D
%D (( A0 A2 -> .plabel= l \sm{1↦5\\2↦6}
%D    A0 A3 ->
%D    A1 A3 ->
%D    
%D ))
%D enddiagram
%D
\pu
$$A =
  \pdiag{an-object-A-of-SetLittleN} \;,
  \qquad
  %\text{and}
  %\quad
  B =
  \pdiag{an-object-B-of-SetLittleN} \;.
$$

There is an inclusion $f: A \ito B$ between them.

We will use a notation with expressions like `$(a∈A(u))$'s to make
certain calculations easier to visualize. For example,
$$\begin{array}{rrl}
  A \psm{▁2 \\ ↓ \\ ▁1} (4∈A(▁2)) = (7∈A(▁1))
    & \text{will mean}
    & A \psm{▁2 \\ ↓ \\ ▁1}: A(▁2) → A(▁1)
    \\
    [7pt]
    & \text{and}
    & A \psm{▁2 \\ ↓ \\ ▁1}(4) = 7,
      \quad \text{and}
    \\
    [10pt]
  f (4∈A(▁2)) = (4∈B(▁2))
    & \text{will mean}
    & f(▁2): A(▁2) → B(▁2)
    \\
    & \text{and}
    & f(▁2)(4) = 4.
  \end{array}
$$

Note that the way to expand the expressions at the left is different
for functors (first case) and for natural transformations (second
case).

\newpage

%  ____       _   ____      _   _ 
% / ___|  ___| |_|  _ \ _  | | | |
% \___ \ / _ \ __| | | (_) | |_| |
%  ___) |  __/ |_| |_| |_  |  _  |
% |____/ \___|\__|____/(_) |_| |_|
%                                 
% «SetD-H»  (to ".SetD-H")
\subsection{The Heyting Algebra $H$}
\label{SetD-H}

A $\SetD$ has several diffent terminal objects. The symbol 1, or
$1_\catE$, will by default denote the one that has $1(u) = \{*\}$ for
every $u∈\catD$. Its class of subsets always forms a set; for example,
when $D = \littlenbulletsp$ we have:
$$\Subsets(1_\catE) = \left\{
  \littlena{0000}, \;
  \littlena{0001}, \;
  \littlena{0010}, \;
  \littlena{0011}, \;
  \littlena{0101}, \;
  \littlena{0111}, \;
  \littlena{1011}, \;
  \littlena{1111}
  \right\}
$$
Note that here the small `0's mean `$∅$' and the small `1's mean
`$\{*\}$'. We will use the notations for piles from \cite[sec.15]{PH1}
to draw this as:
$$H \; = \; \Subsets(1_\catE) = \;\; \zha{LittleN ZHA}$$
but with a difference: in \cite[sec.15]{PH1} two-digit numbers are
interpreted as subsets of the set of points of the current default
two-column graph --- so, there we would have $21 = \csm{2▁, & \\ 1▁, &
  ▁2} ⊂ 𝐃_0$ --- while here they will be subsets of $1_𝐄$; so, here
this holds:
$$21
  = \littlenap{1011}
  = \littlen{ {\{*\}}   {∅}   {\{*\}} {\{*\}} }
  ⊂ \littlen{ {\{*\}} {\{*\}} {\{*\}} {\{*\}} }
  = 1_𝐄 \; .
$$

In contexts in which we have a category $\catD$ the letter $H$ will
denote this set $\Subsets(1_\catE)$, regarded as a Heyting Algebra,
and we will refer to its elements as {\sl truth-values}. In contexts
in which our $\catD$ is a 2-column graph this $H$ will be a {\sl
  Planar} Heyting Algebra --- in the sense of sections 4 and 17 of
\cite{PH1} --- and we will denote these truth-values by two-digit
numbers.

% (ph1p  6 "ZDAGs" "2")
% (ph1a    "ZDAGs")
% (ph1p 28 "2CGs" "14")
% (ph1a    "2CGs")
% (ph1p 29 "topologies-on-2CGs" "15")
% (ph1a    "topologies-on-2CGs")
% (ph1p 32 "converting-ZHAs-2CAGs" "17")
% (ph1a    "converting-ZHAs-2CAGs")

% (find-sh "dict caveat")

%  ____       _   ____                           _       
% / ___|  ___| |_|  _ \ _   _ __   ___  ___  ___| |_ ___ 
% \___ \ / _ \ __| | | (_) | '_ \ / _ \/ __|/ _ \ __/ __|
%  ___) |  __/ |_| |_| |_  | |_) | (_) \__ \  __/ |_\__ \
% |____/ \___|\__|____/(_) | .__/ \___/|___/\___|\__|___/
%                          |_|                           
%
% «SetD-defs-posets»  (to ".SetD-defs-posets")
% (cltp 23 "SetD-defs-posets")
% (clta    "SetD-defs-posets")
\subsection{Some definitions on posets}
\label{SetD-defs-posets}

Suppose that $𝐏$ is a downward-directed poset; we will denote its set
of points as $𝐏_0$. A subset $𝓐⊂𝐏_0$ is a {\sl down-set} iff it obeys
this:
$$∀u,v∈𝐏_0.\; \pmat{u \\ \text{above} \\ v} → \pmat{u∈𝓐 \\ ↓ \\ v∈𝓐}$$

\def\phu {\phantom{1▁}}
\def\phuc{\phantom{1▁,}}
\def\csq#1#2#3#4{\csm{#1&#2\\#3&#4\\}}

We will denote the set of all down-sets of $𝐏$ by $\Downs(𝐏)$, and if
$𝓐⊂𝐏_0$ then we will denote the down-set generated by $𝓐$ --- i.e.,
the smallest down-set of $𝐏$ containing $𝓐$ --- by ${↓}_𝐏 𝓐$, or by
${↓}𝓐$. When $u∈𝐏_0$ we will usually write ${↓}\{u\}$ as ${↓}u$; so,
if
$$\catD = \littlenbulletsp = \tcg{LittleN 2CG}$$
then:
$$\begin{array}{rr}
    & {↓}_𝐃  \{1▁,\,▁2\}
      = {↓}_𝐃 \csq{}{▁2,}{1▁}{}
      =        \csq{}{▁2,}{1▁,}{▁1} ,
      \\[5pt]
    % \text{and:} \phantom{m}
    & {↓}_𝐃 \, 2▁
      = {↓}_𝐃 \{2▁\}
      = {↓}_𝐃 \csq{2▁}{}{}{\phuc}
      =        \csq{2▁,}{\phuc}{1▁,}{▁1} .
  \end{array}
$$

\msk

We will use the poset $𝐃$ above in all examples in this section.

Every ${↓}_𝐏𝓐$ can be seen as a poset --- it inherits the order
from $𝐏$.

\msk

Every object $B$ of a $\SetD$ can be transformed into a poset $\Po(B)$
whose points are the pairs $(u,b)$ in which $u∈𝐃$ and $b∈B(u)$. For
example:
%
%D diagram B-in-SetLittleN-as-poset
%D 2Dx     100 +25  +40 +25
%D 2D  100 A1  A2   A3  A4
%D 2D
%D 2D  +40 A5  A6   A7
%D 2D
%D ren A1 A2  A3 A4 ==> (2▁,1) (2▁,2)   (▁2,3) (▁2,4)
%D ren A5 A6  A7    ==> (1▁,5) (1▁,6)   (▁1,7)
%D
%D (( A1 A5 -> A2 A6 ->
%D    A1 A7 -> A2 A7 ->
%D    A3 A7 -> A4 A7 ->
%D
%D ))
%D enddiagram
%D
$$\pu
  B = \pdiag{an-object-B-of-SetLittleN}
  \qquad
  \Po(B) \;\;=\; \scalebox{0.8}{$\diag{B-in-SetLittleN-as-poset}$}
$$

\def\mym   #1#2#3{\psm{#1 \\ ↓&↘ \\ #2&&#3 }}
\def\mymm#1#2#3#4{\psm{#1 && #2 \\ ↓&↘&↓ \\ #3&&#4 }}

The operation that does the inverse of `$\Po$' is called `$\Ob$'. For
every object $B$ in a $\SetD$ there is a bijection between the subsets
of $B$ in $\SetD$ and the down-sets of $\Po(B)$. Let's define:
%
% $${↓}(b∈B(u)) \;\; := \;\; \Ob({↓}_{\Po(B)} (u,b))$$

% For example,
%
$$\begin{array}{rcl}
  {↓}(b∈B(u)) & := & \Ob({↓}_{\Po(B)} (u,b)),
  \\
  [5pt]
  \text{so:} \phantom{mm}
  {↓}(1∈B(2▁))
    &=& \Ob({↓}_{\Po(B)} (2▁,1)) \\
    &=& \Ob \mym {(2▁,1)} {(1▁,5)} {(▁1,7)} \\
    &=& \mymm {\{1\}} {∅} {\{5\}} {\{7\}} .
        \\
  \end{array}
$$

We will also use this other shorthand: if $u∈𝐃$, then
$${↓}u \; = \; {↓}(*∈1_𝐄(u)).
$$

This looks ambiguous. We have both
$$\begin{array}{rl}
    & {↓}2▁
      \;=\; {↓}(*∈1_𝐄(2▁))
      \;=\; \littlenap{1011}
      \;=\; 21
      \;⊂\; 1_𝐄 \\
    \text{and} \phantom{mmm}
    & {↓}2▁
      \;=\; {↓}_𝐃 \{2▁\}
      \;=\; \csq{2▁,}{}{1▁,}{▁1} % \{2▁,1▁, \, 1▁\}
      \;⊂\; 𝐃\, ,
  \end{array}
$$

but will always use the first meaning.

%  ____       _   ____           _               _  __ _           
% / ___|  ___| |_|  _ \ _    ___| | __ _ ___ ___(_)/ _(_) ___ _ __ 
% \___ \ / _ \ __| | | (_)  / __| |/ _` / __/ __| | |_| |/ _ \ '__|
%  ___) |  __/ |_| |_| |_  | (__| | (_| \__ \__ \ |  _| |  __/ |   
% |____/ \___|\__|____/(_)  \___|_|\__,_|___/___/_|_| |_|\___|_|   
%                                                                  
% «SetD-classifier»  (to ".SetD-classifier")
% (cltp 25 "SetD-classifier")
% (clta    "SetD-classifier")

\subsection{The classifier $\Omega$}
\label{SetD-classifier}

The classifier object on a $\SetD$ is the object $Ω∈\SetD$ whose
action on objects is $Ω(u) = \Downs({↓}u)$ and whose action on
morphisms is $Ω \psm{u \\ ↓ \\ v} (𝓢) = 𝓢∩{↓}v$. Let's decypher this
--- it contains an abuse of language.

Let $𝐃=\littlenbulletsp$ and $u=2▁$. Then ${↓}2▁ = \littlenap{1011} ⊂
1_𝐄$. We will convert that $\littlenap{1011}$ to a poset
$\littlenap{{•}{·}{•}{•}}$, and $\Downs({↓}2▁)$ will be the set of
down-sets of this new poset, which is:
$$\begin{array}{rcl}
  \Downs({↓}2▁) &=&
    \left \{
    \littlena{{0}{·}{0}{0}}, \;
    \littlena{{0}{·}{0}{1}}, \;
    \littlena{{0}{·}{1}{0}}, \;
    \littlena{{0}{·}{1}{1}}, \;
    \littlena{{1}{·}{1}{1}}
    \right\} \\
    &=& \{ 00, 01, 10, 11, 21 \} \\
    &=& \left( \olittlen{ {21}{·} {10}{11}{·} {00}{01}{·} } \right) \\
        [10pt]
    &⊂& H \\
  \end{array}
$$

Each `$·$' means ``this is out of the domain'', but the precise
details vary according to the context.

\def\polittlen#1{\left( \olittlen{#1} \right)}

If we do the same as above for $1▁$, $▁1$, and $▁2$, we get:
$$\begin{array}{ccc}
  \Downs({↓}2▁) = {\polittlen{{21}{·} {10}{11}{·} {00}{01}{·}}} &&
  \Downs({↓}▁2) = {\polittlen{ {·}{·}  {·}{·}{·} {00}{01}{02}}} \\
  \\
  \Downs({↓}1▁) = {\polittlen{ {·}{·} {10}{·}{·}   {00}{·}{·}}} &&
  \Downs({↓}▁1) = {\polittlen{ {·}{·}  {·}{·}{·}  {00}{01}{·}}} \\
  \end{array}
$$

When $\catD = \littlenbulletsp$ we can draw the classifier
of $\SetD$ as:
\def\Db{\littlenbig{
  {\olittlen{{21}{·} {10}{11}{·} {00}{01}{·}}}
  {\olittlen{ {·}{·}  {·}{·}{·} {00}{01}{02}}}
  {\olittlen{ {·}{·} {10}{·}{·}   {00}{·}{·}}}
  {\olittlen{ {·}{·}  {·}{·}{·}  {00}{01}{·}}}
  }}
\def\Da{\littlen{ {\{21\}} {\{20\}} {\{10\}} {\{01\}} }}
$$Ω \;\; = \; \Db$$

When $𝐃$ is a 2-column graph the classifier $Ω∈\SetD$ can always be
drawn in this (nice) way. 

\newpage

% (favp 20 "internal-view-functor")
% (fava    "internal-view-functor")

We still need to understand what the action of $Ω$ on morphisms
``means''. We can use the definition to calculate a particular case,
$$\begin{array}{ccl}
  Ω \psm{u \\ ↓ \\ v} (𝓢) &=& 𝓢∩{↓}v \\[8pt]
  Ω \psm{2▁ \\ ↓ \\ 1▁} (11) &=& 11∩{↓}1▁ \\
                             &=& 11∩10 \;=\; 10 \\
  \end{array}
$$
and can use a diagram like the one in \cite[section 5.2]{FavC}, but
with the particular case in its right half, to put that in a more
categorical form:
%
%D diagram SetD-classifier-restriction-map
%D 2Dx     100  +30 +25    +40  +35 +40
%D 2D  100 A0 - A1  C0     a0 - a1  c0    
%D 2D      |     |         |     |        
%D 2D  +30 A2 - A3  C1     a2 - a3  c1    
%D 2D                                     
%D 2D  +20 B0 - B1         b0 - b1        
%D 2D
%D ren A0 A1 A2 A3 ==> u Ω(u) v Ω(v)
%D ren B0 B1       ==> \catD \Set
%D ren C0 C1       ==> 𝓢 𝓢∩{↓}v
%D
%D ren a0 a1 a2 a3 ==> 2▁ \DDLtwo 1▁ \DDLone
%D ren b0 b1       ==> \littlenbulletsp \Set
%D ren c0 c1       ==> 11 11∩{↓}1▁=10
%D
%D (( A0 A1 |->
%D    A0 A2  ->
%D    A1 A3  ->
%D    A0 A3 harrownodes nil 20 nil |->
%D    A2 A3 |->
%D    B0 B1  -> .plabel= a Ω
%D    C0 C1 |->
%D ))
%D (( a0 a1 |->
%D    a0 a2  ->
%D    a1 a3  ->
%D    a0 a3 harrownodes nil 20 nil |->
%D    a2 a3 |->
%D    b0 b1  -> .plabel= a Ω
%D    c0 c1 |->
%D ))
%D enddiagram
%D
$$\pu
  \def\DDLtwo{\polittlen{{21}{·} {10}{11}{·} {00}{01}{·}}}
  \def\DDRtwo{\polittlen{ {·}{·}  {·}{·}{·} {00}{01}{02}}}
  \def\DDLone{\polittlen{ {·}{·} {10}{·}{·}   {00}{·}{·}}}
  \def\DDRone{\polittlen{ {·}{·}  {·}{·}{·}  {00}{01}{·}}}
  \diag{SetD-classifier-restriction-map}
$$

%  ____       _   ____      _                   
% / ___|  ___| |_|  _ \ _  | |_ _ __ _   _  ___ 
% \___ \ / _ \ __| | | (_) | __| '__| | | |/ _ \
%  ___) |  __/ |_| |_| |_  | |_| |  | |_| |  __/
% |____/ \___|\__|____/(_)  \__|_|   \__,_|\___|
%                                               
% «SetD-true»  (to ".SetD-true")
% (cltp 26 "SetD-true")
% (clta    "SetD-true")

\subsection{The `true' map $⊤:1 \ito Ω$}
\label{SetD-true}

The `true' map $⊤:1 \ito Ω$ will not be an inclusion if we take $1$ as
the default terminal. We can fix this by defining $1_⊤$ as the
terminal that takes each $u∈\catD$ to $\{{↓}u\}$, and by using this as
our `true' map: $⊤:1_⊤ \ito Ω$. Our {\sl default} meaning for $1$ is
still the terminal that takes each $u$ to $\{*\}$, but in contexts
like ``$⊤:1 \ito Ω$'' the default meaning for the `1' will change to
`$1_⊤$'.

When $\catD = \littlenbulletsp$ the `true' map will be this one:

%D diagram SetD-true
%D 2Dx     100
%D 2D  100 \Da
%D 2D
%D 2D  +50 \Db
%D 2D
%D # ren ==>
%D
%D (( \Da \Db `-> .plabel= r ⊤
%D
%D ))
%D enddiagram
%D
$$\pu
  \diag{SetD-true}
$$

\newpage

%  ____       _   ____           _     _ 
% / ___|  ___| |_|  _ \ _    ___| |__ (_)
% \___ \ / _ \ __| | | (_)  / __| '_ \| |
%  ___) |  __/ |_| |_| |_  | (__| | | | |
% |____/ \___|\__|____/(_)  \___|_| |_|_|
%                                        
% «SetD-chi»  (to ".SetD-chi")
\subsection{The classifying map $χ_f$ of an inclusion $f$}
\label{SetD-chi}

\unitlength=1pt  % for the pullbacks

\def\chifu {(χ_f)_u}
\def\chifub{(χ_f)_u(b)}
\def\chifub{χ_f(b∈B(u))}
\def\dnbu  {{↓}(b∈B(u))}

Suppose that this diagram is a pullback in a $\SetD$:
%
%D diagram SetD-chi-1
%D 2Dx     100  +25
%D 2D  100 A0 - A1
%D 2D      |     |
%D 2D  +25 A2 - A3
%D 2D
%D ren A0 A1 A2 A3 ==> A 1_⊤ B Ω
%D
%D (( A0 A1  -> .plabel= a !
%D    A0 A2 `-> .plabel= l f
%D    A1 A3 `-> .plabel= r ⊤
%D    A2 A3  -> .plabel= a χ_f
%D    A0 relplace 7 7 \pbsymbol{7}
%D ))
%D enddiagram
%D
$$\pu
  \diag{SetD-chi-1}
$$
then for every $u∈𝐃$ this is also a pullback:
%
%D diagram SetD-chi-2
%D 2Dx     100  +30 +30  +30
%D 2D  100 A0 - A1  B0 - B1
%D 2D      |     |  |     |
%D 2D  +30 A2 - A3  B2 - B3
%D 2D
%D ren A0 A1 A2 A3 ==> A(u) 1_⊤(u) B(u) Ω(u)
%D ren B0 B1 B2 B3 ==> a {↓}u a (χ_f)_u(a)
%D
%D (( A0 A1  -> .plabel= a !
%D    A0 A2 `-> .plabel= l f_u
%D    A1 A3 `-> .plabel= r ⊤_u
%D    A2 A3  -> .plabel= a (χ_f)_u
%D    A0 relplace 7 7 \pbsymbol{7}
%D
%D    newnode: B3' at: @B3+v(0,-8) .TeX= {↓}u
%D    B0 B1  |->
%D    B1 B3' |->
%D    B0 B2  |->
%D    B2 B3  |->
%D ))
%D enddiagram
%D
$$\pu
  \diag{SetD-chi-2}
$$
so, for all $u∈𝐃$:
$$A(u) \;=\; \setofst{b∈B(u)}{(χ_f)_u(b)={↓}u}
$$

This gives us an elementary way to start from just $χ_f:B→Ω$ and
construct the pullback in a way that makes the left wall an inclusion.

The other direction is harder. Suppose that we have an inclusion $f:A
\ito B$ in a $\SetD$, and we want to construct the map $χ_f:B→Ω$ that
completes the pullback. For every $u∈𝐃$ and $b∈B(u)$ we will define:
$$\chifub \; = \; \Cst(A ∩ \dnbu) 
$$

The operation $\Cst$ is new: it ``canonicalizes subobjects of the
terminal''. If $C∈\SetD$, then $C$ is a subobject of the terminal iff
the (unique) map $!_C: C→1$ is a monic. $\Cst(C)$ is defined if an
only if $C$ is a subobject of the terminal, and is defined as $\Cst(C)
:= \dom(\can(!_C))$. Here are a diagram for the general case and an
example:
\def\CstA{\littlen{    {∅}   {\{4\}} {\{5\}} {\{6\}} }}
\def\CstB{\littlen{    {∅}   {\{*\}} {\{*\}} {\{*\}} }}
\def\CstC{\littlen{  {\{*\}} {\{*\}} {\{*\}} {\{*\}} }}

%D diagram CST-def-and-example
%D 2Dx     100 +45  +40 +45
%D 2D  100 A0  A1   B0  B1
%D 2D
%D 2D  +30     A2       B2
%D 2D
%D ren A0 A1 A2 ==> C \Cst(C) 1_𝐄
%D ren B0 B1 B2 ==> \CstA \CstB \CstC
%D
%D (( A0 A1 <->
%D    A0 A2 >-> .plabel= l !_C
%D    A1 A2 `-> .plabel= r \can(!_C)
%D    
%D    B0 B1 <->
%D    B0 B2 >->
%D    B1 B2 `->
%D
%D    newnode: C1 at: @B1+v(25,0) .TeX= =12 place
%D    newnode: C2 at: @B2+v(25,0) .TeX= =22 place
%D ))
%D enddiagram
%D
$$\pu
  \diag{CST-def-and-example}
$$

See the property Inc2 of toposes with inclusions in
sec.\ref{inclusions-precisely}.

\msk

\newpage

Let's use an example to understand how the definition of $\chifub$
above works. The diagram below defines objects $A,B∈\SetD$ and an
inclusion $f: A \ito B$, and the map $⊤: 1 \ito Ω$ at its right wall
is the one that we saw in sec.\ref{SetD-true}:

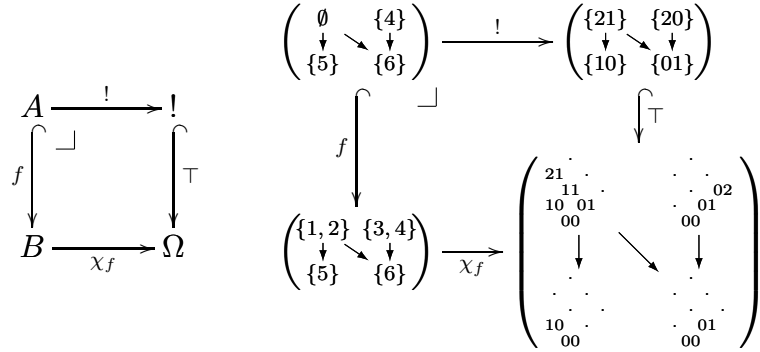
\def\Aa{\littlena{0011}}
\def\Ab{\littlena{1011}}
\def\Ac{\littlena{1111}}
\def\Ba{\littlen{    {∅}       {∅}    {\{5\}} {\{6\}} }}
\def\Bb{\littlen{   {\{2\}}    {∅}    {\{5\}} {\{6\}} }}
\def\Ca{\littlen{    {∅}      {\{4\}} {\{5\}} {\{6\}} }}
\def\Cb{\littlen{ {\{1,2\}} {\{3,4\}} {\{5\}} {\{6\}} }}
\def\Cc{\littlen{    {∅}      {\{4\}}   {∅}   {\{6\}} }}
\def\Cd{\littlen{    {∅}       {∅}      {∅}   {\{6\}} }}
\def\Da{\littlen{ {\{21\}} {\{20\}} {\{10\}} {\{01\}} }}
\def\Db{\littlen{ {↓21} {↓20} {↓10} {↓01} }}
\def\Db{\littlenbig{
  {\olittlen{{21}{·} {10}{11}{·} {00}{01}{·}}}
  {\olittlen{ {·}{·}  {·}{·}{·} {00}{01}{02}}}
  {\olittlen{ {·}{·} {10}{·}{·}   {00}{·}{·}}}
  {\olittlen{ {·}{·}  {·}{·}{·}  {00}{01}{·}}}
  }}

%D diagram PB-LittleN-generic
%D 2Dx     100  +30
%D 2D  100 A1 - A2
%D 2D      |     | 
%D 2D  +30 A4 - A5
%D 2D
%D ren A1 A2 ==> A !
%D ren A4 A5 ==> B Ω
%D
%D (( A1 A2  -> .plabel= a !
%D    A1 A4 `-> .plabel= l f A2 A5 `-> .plabel= r ⊤
%D    A4 A5  -> .plabel= b χ_f 
%D    A1 relplace 7 7 \pbsymbol{7}
%D ))
%D enddiagram
%D
%D diagram PB-LittleN
%D 2Dx     100  +60
%D 2D  100 A1 - A2
%D 2D      |     | 
%D 2D  +45 A4 - A5
%D 2D
%D ren A1 A2 ==> \Ca \Da
%D ren A4 A5 ==> \Cb \Db
%D
%D (( A1 A2  -> .plabel= a !
%D    A1 A4 `-> .plabel= l f A2 A5 `-> .plabel= r ⊤
%D    A4 A5  -> .plabel= b χ_f 
%D    A1 relplace 15 12 \pbsymbol{7}
%D ))
%D enddiagram
%D
\pu
$$\diag{PB-LittleN-generic}
  \qquad
  \diag{PB-LittleN}
$$

\def\Bub{B_{ub}}

It {\sl does not} define the map $χ_f$, but we can use the formula
above to calculate $χ_f(b∈B(u))$ for each pair $(u,b)$ with $u∈𝐃$ and
$b∈B(u)$. Let's do the case $u=2▁$, $b=2$:
$$\def\mini#1{\scalebox{0.7}{$#1$}}
  \def\mini#1{#1}
  \begin{array}{rcl}
    χ_f(b∈B(u))   &=& \Cst(A ∩ {↓}(b∈B(u))) \\
    χ_f(2∈B(2▁))  &=& \Cst(A ∩ {↓}(2∈B(2▁))) \\
                  &=& \Cst( \mini{\Ca} ∩ \mini{\Bb} ) \\
                  &=& \Cst( \mini{\Ba} ) \\
                  &=& \littlenap{0011} \\
                  &=& 11.
  \end{array}
$$

We can also visualize the whole construction at once by drawing the
diagram in the next page.

\newpage

Note that we are abbreviating $\dnbu$ as $\Bub$, and that the left
half of diagram will be different for each choice of $u∈𝐃$ and
$b∈B(u)$.
%
%D diagram classifying-map-0
%D 2Dx     100  +60  +30  +30
%D 2D  100 B0 - A0 - A1 - A2
%D 2D      |    |    |     | 
%D 2D  +25 B1 - A3 - A4 - A5
%D 2D      |
%D 2D  +25 B2   D0 - D1 - D2
%D 2D
%D ren B0 B1 B2 ==>         \Cst(A{∩}\Bub) {↓}u=\Cst(\Bub) 1
%D ren B0 B1 B2 ==> \chifub=\Cst(A{∩}\Bub) {↓}u=\Cst(\Bub) 1
%D ren A0 A1 A2 ==> A{∩}\Bub A 1
%D ren A3 A4 A5 ==> \Bub B Ω
%D
%D (( A0 A1 `-> A1 A2  -> .plabel= a !
%D    A0 A3 `-> A1 A4 `-> .plabel= l f A2 A5 `-> .plabel= r ⊤
%D    A3 A4 `-> A4 A5  -> .plabel= b χ_f 
%D    A0 A3 `->
%D    A0 relplace 9 8 \pbsymbol{7}
%D    A1 relplace 7 7 \pbsymbol{7}
%D    B2 xy+= 0 0
%D    B0 A0 <-> # B0 A3 >->
%D    B0 B1 `->
%D    B1 A3 <-> B1 B2 `-> A3 B2 >->
%D ))
%D enddiagram
%D
$$\pu
  \diag{classifying-map-0}
$$

When $u=2▁$ and $b=2$ the diagram above becomes:
%
%D diagram classifier-LittleN
%D 2Dx     100  +50  +50  +60
%D 2D  100 B0 - A0 - A1 - A2
%D 2D      |    |    |     | 
%D 2D  +45 B1 - A3 - A4 - A5
%D 2D      |
%D 2D  +35 B2
%D 2D
%D ren B0 A0 A1 A2 ==> \Aa \Ba \Ca \Da
%D ren B1 A3 A4 A5 ==> \Ab \Bb \Cb \Db
%D ren B2          ==> \Ac
%D
%D (( A0 A1 `-> A1 A2  -> .plabel= a !
%D    A0 A3 `-> A1 A4 `-> .plabel= l f A2 A5 `-> .plabel= r ⊤
%D    A3 A4 `-> A4 A5  -> .plabel= b χ_f 
%D    A0 A3 `->
%D    A0 relplace 15 12 \pbsymbol{7}
%D    A1 relplace 15 12 \pbsymbol{7}
%D    B2 xy+= 0 0
%D    B0 A0 <-> # B0 A3 >->
%D    B0 B1 `->
%D    B1 A3 <-> B1 B2 `-> A3 B2 >->
%D ))
%D enddiagram
%D
\def\Aa{11=\littlenap{0011}}
\def\Ab{21=\littlenap{1011}}
\def\Ac{22=\littlenap{1111}}
$$\pu
  %\scalebox{0.9}{$
    \diag{classifier-LittleN}
  %$}
$$

\newpage

%  ____       _   ____          _                 
% / ___|  ___| |_|  _ \ _      | | ___ __ _ _ __  
% \___ \ / _ \ __| | | (_)  _  | |/ __/ _` | '_ \ 
%  ___) |  __/ |_| |_| |_  | |_| | (_| (_| | | | |
% |____/ \___|\__|____/(_)  \___/ \___\__,_|_| |_|
%                                                 
% «SetD-Jcan»  (to ".SetD-Jcan")
% (cltp 29 "SetD-Jcan")
% (clta    "SetD-Jcan")

\subsection{The canonical Grothendieck topology $\Jcan$}
\label{SetD-Jcan}

Here's how to define the canonical Grothendieck topology on a
topological space $(X,\Opens(X))$. We will denote open sets of
$\Opens(X)$ by letters like $U,V,W$, and sets of open sets by letters
like $𝓐, 𝓢, 𝓤$ --- so $U,V,W∈\Opens(X)$ and $𝓐, 𝓢, 𝓤⊂\Opens(X)$. We
will regard $\Opens(X)$ as a downward-directed poset, so these are
equivalent:
$$\mat{U \\ \Above \\ V} \;\;,
  \qquad
  \mat{U \\ ↓ \\ V} \;\;,
  \quad
  \text{and}
  \quad
  \mat{U \\ \rotl{⊆} \\ V} \;\;.
$$

A subset $𝓢⊂\Opens(X)$ will be called a {\sl sieve} if it is
downward-closed. A sieve $𝓢⊂\Opens(X)$ will be called a {\sl sieve on}
$U$ if $𝓢⊂\Opens(U)$. We say that a sieve $𝓢⊂\Opens(U)$ {\sl covers}
$U$, or is a {\sl covering sieve on} $U$, if $\bigcup𝓢 = U$, and we
say that a sieve $𝓢$ {\sl r-covers} $U$ if its restriction to
$\Opens(U)$, $𝓢∩\Opens(U)$, covers $U$.

The set of all sieves on $U$ will be denoted by $Ω(U)$, and the set of
all covering sieves on $U$ will be denoted by $\Jcan(U)$. Some extra
notational conventions: The letters $𝓢$ and $𝓡$ will always denote
sieves, and $𝓤$ will always denote a covering sieve (on $U$).

This diagram shows all these notional conventions, plus a few more:

\def\defJanU{ \setofst{𝓤∈Ω(U)}{\bigcup𝓤=U} \hspace*{-2.5cm} }
\def\defJanU{ \setofst{𝓤∈\Downs({↓}U)}{\bigcup𝓤=U} \hspace*{-3.2cm} }

\def\hboxl#1{\hbox to 35pt{$#1$\hss}}
\def\defJanU{\hboxl{\setofst{𝓢∈\Downs({↓}U)}{\bigcup𝓢=U}}}
\def\defTU  {\hboxl{\setofst{𝓢∈\Downs({↓}U)}{𝓢={↓}U}}}

$$\begin{matrix}
  U       &∈& \Opens(X) \\
  \rotl{⊂} && \rotl{⊂} \\
  V       &∈& 𝓢       &∈& Ω(U)      &=& \Downs({↓}U) &⊂& \Pts({↓}U) \\
           &&           && \rotl{⊂}   && \rotl{⊂} \\
  W       &∈& 𝓤       &∈& \Jcan(U)  &=& \defJanU  \\
           &&           && \rotl{⊂}   && \rotl{⊂} \\
           &&           && 1_⊤(U)    &=& \defTU    \\
  \end{matrix}
$$

% (lindp 8 "2")
% (linda   "2")

This $\Jcan$ has this three properties:
$$\scalebox{0.9}{$
  \begin{array}{rccc}
  \hasmax_\Jcan: & ∀U∈\Opens(X).
                 &
                 & {↓}U \Covers U
                   \;,
                 \\[7.5pt]
  \stab_\Jcan:   & ∀\pmat{U \\
                          \Above \\
                          V} ∈ \Opens(X).
                 & ∀𝓢∈Ω(U).
                 & \pmat{𝓢 \Covers U \\
                          ↓ \\
                         𝓢 \RCovers V \\
                        }
                   \;,
                 \\[22pt]
  \trans_\Jcan:  & ∀U∈\Opens(X).∀𝓤∈J(U).
                 & ∀𝓢∈Ω(U).
                 & \pmat{𝓢 \Covers U \\
                          ↑ \\
                         ∀V∈𝓤. \; 𝓢 \RCovers V \\
                        }
                   \;.
                 \\
  \end{array}
  $}
$$

\newpage

%  ____       _   ____          _ 
% / ___|  ___| |_|  _ \ _      | |
% \___ \ / _ \ __| | | (_)  _  | |
%  ___) |  __/ |_| |_| |_  | |_| |
% |____/ \___|\__|____/(_)  \___/ 
%                                 
% «SetD-J»  (to ".SetD-J")
% (cltp 30 "SetD-J")
% (clta    "SetD-J")
\subsection{Other Grothendieck topologies}
\label{SetD-J}

Here's how to define what is a Grothendieck topology on an arbitrary
downward-directed posed $𝐃$. The topology $\Opens(X)$ of the previous
case will become the poset $𝐃$; we will refer to this $𝐃$ as our {\sl
  ex-topology} and to the points $u,v,w∈𝐃$ as {\sl ex-open sets}.
These are equivalent:
$$\mat{u \\ \Above \\ v} %\;\;,
  %\qquad
  %\mat{U \\ ↓ \\ V} \;\;,
  \qquad
  \text{and}
  \qquad
  \mat{u \\ ↓ \\ v} \;\;.
$$

We will denote ex-open sets by letters like $u$, $v$, $w$, and sets of
ex-open sets by letters like $𝓐$, $𝓢$, and $𝓤$. A subset $𝓢⊂𝐃$ will be
called a {\sl sieve} if it is downward-closed. A sieve $𝓢⊂𝐃$ will be
called a {\sl sieve on} $u$ if $𝓢⊂{↓}u$. The down-sets $\Opens(U)$ of
the previous section will be replaced by ${↓}u$ here.

In this context the set of all sieves is exactly $H$ (see
sec.\ref{SetDs}) and the set of all sieves on $u$ is exactly $Ω(u)$
(see sec.\ref{SetD-classifier}).

% --- modulo a small abuse of language that we will explain and
% --- correct soon.

A {\sl Grothendieck topology} is an object $J∈\SetD$ obeying $1_⊤⊂J⊂Ω$
and having the properties $\hasmax_J$, $\stab_J$, and $\trans_J$, that
are defined as:
$$\scalebox{0.9}{$
  \begin{array}{rccc}
  \hasmax_J: & ∀u∈𝐃.
             &
             & {↓}u \JCovers u
               \;,
             \\[7.5pt]
  \stab_J:   & ∀\pmat{u \\
                      \Above \\
                      v} ∈ 𝐃.
             & ∀𝓢∈Ω(u).
             & \pmat{𝓢 \JCovers u \\
                      ↓ \\
                     𝓢 \RJCovers v \\
                    }
               \;,
             \\[22pt]
  \trans_J:  & ∀u∈𝐃.∀𝓤∈J(u).
             & ∀𝓢∈Ω(u).
             & \pmat{𝓢 \JCovers u \\
                      ↑ \\
                     ∀v∈𝓤. \; 𝓢 \RJCovers v \\
                    }
               \;.
             \\
  \end{array}
  $}
$$

Here we say that a sieve $𝓢$ on $u$ {\sl $J$-covers} $u$ if $𝓢∈J(u)$,
and that a sieve $𝓢$ {\sl r-$J$-covers} $v$ if its restriction to
${↓}v$, $𝓢∩{↓}v$, $J$-covers $v$ --- i.e., if $𝓢∩{↓}v∈J(v)$.

Here our notational conventions are that $𝓡$, $𝓢$, and $𝓤$ are sieves,
and that $𝓤$ is a $J$-covering sieve on $u$. Other calligraphic
capitals, like $𝓐$, may denote subsets of $𝐃_0$ that don't need to be
downward-closed.

This diagram shows all these notional conventions, plus a few more:

\def\hboxl#1{\hbox to 35pt{$#1$\hss}}
\def\defJU  {\hboxl{\setofst{𝓢∈\Downs({↓}U)}{\bigcup𝓢=U}}}
\def\defTU  {\hboxl{\setofst{𝓢∈Ω(u)}{𝓢={↓}u}}}
$$\begin{matrix}
  u       &∈& 𝐃 \\
  ↓       && \rotl{⊂} \\
  v       &∈& 𝓢       &∈& Ω(u)      &=& \Downs({↓}u) &⊂& \Pts({↓}u) \\
           &&           && \rotl{⊂}  \\
  w       &∈& 𝓤       &∈& J(u)      \\
           &&           && \rotl{⊂}   \\
           &&           && 1_⊤(u)    &=& \defTU    \\
  \end{matrix}
$$

% [TODO: explain the abuse of notation]

% (lindp 8 "2")
% (linda   "2")

%  ____       _   ____          _  __   __     __ _ _ _            
% / ___|  ___| |_|  _ \ _      | |/ /   \ \   / _(_) | |_ ___ _ __ 
% \___ \ / _ \ __| | | (_)  _  | | | | | | | | |_| | | __/ _ \ '__|
%  ___) |  __/ |_| |_| |_  | |_| | | |_| | | |  _| | | ||  __/ |   
% |____/ \___|\__|____/(_)  \___/| |\__,_| | |_| |_|_|\__\___|_|   
%                                 \_\   /_/                        
%                                 
% «SetD-Ju-filter»  (to ".SetD-Ju-filter")
% (cltp 32 "SetD-Ju-filter")
% (clta    "SetD-Ju-filter")
\subsection{Every $J(u)$ is a filter}
\label{SetD-Ju-filter}

Remember that if $𝐏$ is an upward-directed poset with terminal object
$⊤$ and binary meet $∧$ then a subset $𝐅⊂𝐏$ is a filter if it contains
$⊤$ and is closed upwards and by binary meets. More formally, $𝐅$ is a
filter in $𝐏$ if:
$$\begin{array}{ccc}
              & \ph{mmmmmm} ⊤∈𝐅, \\
              [5pt]
  ∀𝓡,𝓢∈𝐏.   & \pmat{𝓢 \\ \Above \\ 𝓡}
                →
                \pmat{𝓢∈𝐅 \\ ↑ \\ 𝓡∈𝐅}, \\
              [20pt]
  ∀𝓡,𝓢∈𝐏.   & \ph{mmmmm} 
                \pmat{𝓡,𝓢∈𝐅 \\ ↓ \\ 𝓡∧𝓢∈𝐅} . \\
  \end{array}
$$

Also, a filter $𝐅⊂𝐏$ is {\sl principal} when it contains the meet of
all its elements; when this happens we have $𝐅 = {↑}_𝐏 \{\bigwedge𝐅\}
= {↑} \bigwedge𝐅$, and we say that $𝐅$ is generated by its bottom
element $\bigwedge𝐅$. When $𝐏$ is a finite poset all the filters on
$𝐏$ are principal.

\msk

\Theoremsubsection
\label{thm:every-Ju-is-a-filter}
If $J$ is a Grothendieck topology on a $\SetD$ then every $J(u)$ is a
filter on $Ω(u)$.

{\bf Proof.} We can re-state this as 1) $J(u)$ contains the top
element of $Ω(u)$, 2) $J(u)$ is closed upwards, 3) $J(u)$ is closed by
binary meets, and we can re-state that again in a more visual way as:
$$\begin{array}{rccc}
  1) &
               & \ph{mmmmm} {↓}u∈J(u), \\
               [5pt]
  2) &
  ∀𝓡,𝓢∈Ω(u). & \pmat{𝓢 \\ \rotl{⊂} \\ 𝓡}
                →
                \pmat{𝓢∈J(u) \\ ↑ \\ 𝓡∈J(u)}, \\
              [20pt]
  3) &
  ∀𝓡,𝓢∈Ω(u).   & \ph{mmmmm} 
                \pmat{𝓡,𝓢∈J(u) \\ ↓ \\ 𝓡∩𝓢∈J(u)} . \\
  \end{array}
$$

\def\R#1{\ColorRed{#1}}
\def\H{\hspace}

Part (1) is an obvious consequence of $\{{↓}u\} = 1_⊤(u) ⊂ J(u) ⊂
Ω(u)$. The proofs of (2) and (3) are quite technical and difficult to
understand intuitively, so we will present them in Natural Deduction
form and let the reader check that every step is correct.

This is the proof that $J(u)$ is upwards-closed:
%:
%:    [v∈𝓡]^1  𝓡⊂𝓢                        𝓡∈Ω(u)
%:    ---------------                        -------
%:           v∈𝓢         𝓢∈Ω(u)   [v∈𝓡]^1   𝓡⊂𝐃         
%:           --------------------   ----------------
%:                 {↓}v⊂𝓢            v∈𝐃             \hasmax_J
%:             ---------------        --------------------------
%:             𝓢∩{↓}v={↓}v           {↓}v∈J(v)
%:             ------------------------------
%:                                𝓢∩{↓}v∈J(v)
%:                                -------------
%:                                    \SrJcoversv
%:                                    ----------------1
%:  [𝓡∈J(u)]^2  𝓢∈Ω(u)  \H{-4.5cm}  ∀v∈𝓡.\SrJcoversv   \H{-3.5cm}   \trans_J
%:  ---------------------------------------------------------------------------
%:           𝓢∈J(u)
%:           ---------------2
%:           𝓡∈J(u)→𝓢∈J(u)
%:
%:           ^J(u)-is-a-filter-2
%:
\pu
$$\def\SaboveR{𝓢 \Above 𝓡}
  \def\SrJcoversv{𝓢 \RJCovers v}
  \scalebox{0.95}{$
  \ded{J(u)-is-a-filter-2}
  $}
$$

% (grsp 15 "lindenhovius-filter")
% (grsa    "lindenhovius-filter")

And this is the proof that $J(u)$ is closed by binary meets:
%:
%:
%:                                      [v∈𝓡]^1     𝓡∈J(u)          
%:                                     ---------    -------          
%:                                     {↓}v⊂{↓}𝓡   {↓}𝓡=𝓡        
%:                                     ---------------------
%:             𝓢∈J(u)                     {↓}v⊂𝓡                 
%:             -------                  -----------               
%:             𝓢\jcu    \stab_J        𝓡∩{↓}v={↓}v               
%:             ------------------        ----------------------
%:   𝓡∈J(u)   ∀v∈{↓}u.𝓢\rjcv             𝓡∩𝓢∩{↓}v=𝓢∩{↓}v
%:   -------   -------------------       ------------------------
%:   𝓡⊂{↓}u   ∀v∈{↓}u.𝓢∩{↓}v∈J(v)          𝓢∩{↓}v=(𝓡∩𝓢)∩{↓}v
%:   -----------------------------       ----------------------1
%:      ∀v∈𝓡.𝓢∩{↓}v∈J(v)                 ∀v∈𝓡.𝓢∩{↓}v=(𝓡∩𝓢)∩{↓}v
%:      -------------------------------------------------------   
%:                                   ∀v∈𝓡.(𝓡∩𝓢)∩{↓}v∈J(v)
%:             ========                --------------
%:   𝓡∈J(u)   𝓡∩𝓢∈Ω(u)  \H{-4cm}    ∀v∈𝓡.(𝓡∩𝓢)\rjcv      \H{-3.5cm}  \trans_J
%:   ------------------------------------------------------------------------------
%:                                  𝓡∩𝓢∈J(u)
%:
%:                                  ^2.4b
%:
\pu
$$\def\jcu {\JCovers u}
  \def\rjcv{\RJCovers v}
  \scalebox{0.9}{$
  \ded{2.4b}
  $}
$$

% \standout{Here}

Remember that we established (in sec.\ref{SetDs}) that $𝐃$ always
stands for a {\sl finite} downward-directed poset. So in a $\SetD$ all
`$Ω(u)$'s are finite upward-directed posets, and all filters on each
$Ω(u)$ are principal and can recovered from their bottom elements.
This gives us a very compact way to represent Grothendieck topologies.
For example, if $𝐃=\littlenbulletsp$ and
$$\def\foo#1#2{\bigwedge J(#1) = #2}
  \begin{array}{lll}
    \foo{2▁}{10}, && \foo{▁2}{00}, \\
    \foo{1▁}{10}, && \foo{▁1}{00}  \\
  \end{array}
$$
then $1_⊤ ⊂ J ⊂ Ω(u)$ is:
$$\def\DA{\littlenbig{
           {\olittlen{{21}{·}  {·}{·}{·}    {·}{·}{·}}}
           {\olittlen{ {·}{·}  {·}{·}{·}   {·}{·}{02}}}
           {\olittlen{ {·}{·} {10}{·}{·}    {·}{·}{·}}}
           {\olittlen{ {·}{·}  {·}{·}{·}   {·}{01}{·}}}
         }}
  \def\DB{\littlenbig{
           {\olittlen{{21}{·} {10}{11}{·}   {·}{·}{·}}}
           {\olittlen{ {·}{·}  {·}{·}{·} {00}{01}{02}}}
           {\olittlen{ {·}{·} {10}{·}{·}    {·}{·}{·}}}
           {\olittlen{ {·}{·}  {·}{·}{·}  {00}{01}{·}}}
         }}
  \def\DC{\littlenbig{
           {\olittlen{{21}{·} {10}{11}{·} {00}{01}{·}}}
           {\olittlen{ {·}{·}  {·}{·}{·} {00}{01}{02}}}
           {\olittlen{ {·}{·} {10}{·}{·}   {00}{·}{·}}}
           {\olittlen{ {·}{·}  {·}{·}{·}  {00}{01}{·}}}
         }}
  \DA ⊂ \DB ⊂ \DC
$$

\newpage

%  ____       _   ____                       _                
% / ___|  ___| |_|  _ \ _   _ __  _   _  ___| | ___ _   _ ___ 
% \___ \ / _ \ __| | | (_) | '_ \| | | |/ __| |/ _ \ | | / __|
%  ___) |  __/ |_| |_| |_  | | | | |_| | (__| |  __/ |_| \__ \
% |____/ \___|\__|____/(_) |_| |_|\__,_|\___|_|\___|\__,_|___/
%                                                             
% «SetD-nucleus»  (to ".SetD-nucleus")
% (cltp 34 "SetD-nucleus")
% (clta    "SetD-nucleus")
\subsection{Every closure operator $\clop$ induces a nucleus $\nuc$}
\label{SetD-nucleus}

In this section we will suppose that our topos $\SetD=𝐄$ has a closure
operator $\clop$.

If $𝓢$ is a truth-value in $𝐄$ then we will denote the closure of $𝓢
\ito 1$ by $𝓢^* \ito 1$, and if $𝓡$ is a truth-value ``smaller than
$𝓢$'', in the sense that we have an inclusion $𝓡 \ito 𝓢$, then we will
denote the closure of the inclusion $𝓡 \ito 𝓢$ by $𝓡^𝓢 \ito 𝓢$.
Formally,
%
%D diagram S^*-and-R^S-defs
%D 2Dx     100 +20 +40 +20
%D 2D  100 A0      B0
%D 2D
%D 2D  +20     A1      B1
%D 2D
%D 2D  +20 A2      B2
%D 2D
%D ren A0 A1 A2 ==> 𝓢 𝓢^* 1
%D ren B0 B1 B2 ==> 𝓡 𝓡^𝓢 𝓢
%D
%D (( A0 A2 `-> .plabel= l \inc(𝓢,1)
%D    A1 A2 `-> .plabel= r \ovl{\inc(𝓢,1)}
%D    A0 A1 `->
%D
%D    B0 B2 `-> .plabel= l \inc(𝓡,𝓢)
%D    B1 B2 `-> .plabel= r \ovl{\inc(𝓡,𝓢)}
%D    B0 B1 `->
%D ))
%D enddiagram
%D
$$\pu
  \diag{S^*-and-R^S-defs}
  \qquad
  \begin{array}{rccl}
    \nuc: & H & → & H \\ 
          & 𝓢 & → & \dom(\ovl{\inc(𝓢,1)}) \\ 
    \\
    \nus: & {↓}𝓢 & → & {↓}𝓢 \\ 
          & 𝓡 & → & \dom(\ovl{\inc(𝓡,𝓢)}) \\ 
  \end{array}
$$

% A {\sl modality} on a Heyting Algebra $H$ --- we will use the
% terminology from \cite[p.163]{BellLST} instead of the current term
% ``nucleus'' --- is an operation $\nuc: H→H$ that obeys M1, M2, and
% M3 below, and if $𝓢∈H$ then an operation $\nus: {↓}𝓢 → {↓}𝓢$ is a
% {\sl modality on the down-set ${↓}𝓢$} if it obeys MD1, MD2, and MD4
% below:

A {\sl nucleus} on a Heyting Algebra $H$ is an operation $\nuc: H→H$
that obeys M1, M2, and M3 below, and if $𝓢∈H$ then an operation $\nus:
{↓}𝓢 → {↓}𝓢$ is a {\sl modality on the down-set ${↓}𝓢$} if it obeys
MD1, MD2, and MD4 below:
$$\begin{array}{l}
  \text{M1) } 𝓢 ⊂ 𝓢^*,                 \\
  \text{M2) } 𝓢^* = 𝓢^{**},            \\
  \text{M3) } 𝓡^* ∩ 𝓢^* = (𝓡∩𝓢)^*,   \\
  \end{array}
  \qquad
  \begin{array}{l}
  \text{MD1) } 𝓡 ⊂ 𝓡^𝓢,                 \\
  \text{MD2) } 𝓡^𝓢 = 𝓡^{𝓢𝓢},          \\
  \text{MD3) } 𝓠^𝓢 ∩ 𝓡^𝓢 = (𝓠∩𝓡)^𝓢. \\
  \end{array}
$$

\Theoremsubsection
\label{thm:modality}
The operation $\nuc$ induced by the closure operator $\clop$ is a
modality on $H$, and the operation $\nus$ is a modality on the
down-set ${↓}𝓢$.

\def\oc{\ovl{c}}
\def\od{\ovl{d}}
\def\ocd{\ovl{c∩d}}

{\bf Proof.} We will only prove MD1, MD2, and MD3. The proofs for MD1
and MD2 are essentially the same as the proofs for C1 and C2 in
sec.\ref{clops}; the diagram is the one at the left below. The proof
of MD3 is based on the proof of C4 in sec.\ref{clops}. Let $c:𝓠 \ito
𝓢$ and $d:𝓡 \ito 𝓢$. We know from C4 that $χ_{(\oc∩\od)}$ =
$χ_{\ocd}$, so $\dom(\oc∩\od)=𝓠^𝓢∩𝓡^𝓢$ and $\dom(\ocd)=(𝓠∩𝓡)^𝓢$ are
the same object.
%
%D diagram MD1-MD2
%D 2Dx     100  +25  +25  +25  +25  +20  +25  +25
%D 2D  100 A0        B0 - B1 - B2   C0 - C1
%D 2D      |  \      |              |  \    \
%D 2D  +25 |    A1   |              |    C2 - C3
%D 2D      |  /      |              |  /
%D 2D  +25 A2        B3             C4
%D 2D
%D ren A0 A1 A2       ==> 𝓡 \ovl{D} 𝓢
%D ren B0 B1 B2 B3    ==> 𝓡 𝓡^𝓢 𝓡^{𝓢𝓢} 𝓢
%D ren C0 C1 C2 C3 C4 ==> C D \ovl{C} \ovl{D} E
%D
%D (( B1 xy+= 0 5
%D    B2 xy+= -5 20
%D
%D    B0 B1 `->
%D    B1 B2 =
%D    B0 B3 `-> # .plabel= r d
%D    B1 B3 `-> # .plabel= r \ovl{d}
%D    B2 B3 `-> # .plabel= r \ovl{\ovl{d}}
%D ))
%D enddiagram
%D
$$\pu
  \diag{MD1-MD2}
  \qquad
  \begin{array}{rcl}
  \oc&: & 𝓠^𝓢 \ito 𝓢 \\
  \od&: & 𝓡^𝓢 \ito 𝓢 \\
  \oc∩\od&: & (𝓠^𝓢 ∩ 𝓡^𝓢) \ito 𝓢 \\
  c∩d&:     & (𝓠 ∩ 𝓡) \ito 𝓢 \\
  \ocd&:    & (𝓠 ∩ 𝓡)^𝓢 \ito 𝓢 \\
  \dom(\oc∩\od) &=& 𝓠^𝓢∩𝓡^𝓢  \\
  \dom(\ocd)    &=& (𝓠∩𝓡)^𝓢   \\ 
  \end{array}
$$

\newpage

%  ____       _   ____                             __      _ 
% / ___|  ___| |_|  _ \ _   _ __  _   _  ___       \ \    (_)
% \___ \ / _ \ __| | | (_) | '_ \| | | |/ __|  _____\ \   | |
%  ___) |  __/ |_| |_| |_  | | | | |_| | (__  |_____/ /   | |
% |____/ \___|\__|____/(_) |_| |_|\__,_|\___|      /_/   _/ |
%                                                       |__/ 
%
% «SetD-nuc-to-j»  (to ".SetD-nuc-to-j")
% (cltp 35 "SetD-nuc-to-j")
% (clta    "SetD-nuc-to-j")
\subsection{Every nucleus $\nuc$ induces a topology $j$}
\label{SetD-nuc-to-j}

Let's start with a lemma.

\Lemmasubsection
\label{thm:downmodality}
The operations $\nuc$ and $\nus$ induced by the closure operator
$\clop$ are related in this way: $𝓡^𝓢 = 𝓡^*∩𝓢$.

{\bf Proof.} This is a corollary of Theorem \ref{thm:restr-clop-1}.
Let's call our inclusions $m: 𝓡 \ito 𝓢$, $d: 𝓢 \ito 1$, and $c: 𝓡 \ito
1$. Our diagrams are:
%
% (cltp 13 "restricting-a-clop")
% (clta    "restricting-a-clop")
%
%D diagram downmodality
%D 2Dx     100  +20 +20  +20 +20  +20 +22  +20 
%D 2D  100 A0 ----- B0       C0 ----- D0       
%D 2D      |  \     |  \     |  \     |  \     
%D 2D  +20 |    A1 -|--- B1  |    C1 -|--- D1  
%D 2D      |  /     |  /     |  /     |  /     
%D 2D  +20 A2 ----- B2       C2 ----- D2       
%D 2D
%D ren A0 A1 A2 ==> d^{-1}(𝓡) \ovl{d^{-1}(𝓡)} 𝓢
%D ren B0 B1 B2 ==> 𝓡 𝓡^* 1
%D ren C0 C1 C2 ==> 𝓡 𝓡^𝓢{=}𝓡^*{∩}𝓢 𝓢
%D ren D0 D1 D2 ==> 𝓡 𝓡^* 1
%D
%D (( A0 B0 `->
%D    A1 B1 `->
%D    A2 B2 `-> .plabel= b d
%D    
%D    A0 A2 `-> .PLABEL= _(0.25) d^{-1}(c)
%D    A0 A1 `->
%D    A1 A2 `-> .PLABEL= ^(0.30) \ovl{d^{-1}(c)}
%D    
%D    B0 B2 `-> .PLABEL= _(0.25) c
%D    B0 B1 `->
%D    B1 B2 `-> .PLABEL= ^(0.30) \ovl{c}
%D ))
%D (( C0 D0 `->
%D    C1 D1 `->
%D    C2 D2 `-> .plabel= b d
%D    
%D    C0 C2 `-> .PLABEL= _(0.25) m
%D    C0 C1 `->
%D    C1 C2 `-> .PLABEL= ^(0.30) \!\!\ovl{m}=d^{-1}(\ovl{c})
%D    
%D    D0 D2 `-> .PLABEL= _(0.25) c
%D    D0 D1 `->
%D    D1 D2 `-> .PLABEL= ^(0.30) \ovl{c}
%D ))
%D enddiagram
%D
$$\pu
  \diag{downmodality}
$$

\Theoremsubsection
\label{thm:nuc-to-j}

a) If $f:𝓡 \ito 𝓢$ and $𝓢={↓}u$ then $χ_f(*∈𝓢(u)) = (𝓡∈Ω(u))$.

b) If $g:𝓡^𝓢 \ito 𝓢$ and $𝓢={↓}u$ then $χ_g(*∈𝓢(u)) = (𝓡^𝓢∈Ω(u))$.

c) If $h:𝓡 \ito 𝓢$ and $𝓢={↓}u$ then $χ_{\ovl{h}}(*∈𝓢(u)) = (𝓡^𝓢∈Ω(u))$.

d) If $h:𝓡 \ito 𝓢$ and $𝓢={↓}u$ then $(j∘χ_h)(*∈𝓢(u)) = (𝓡^𝓢∈Ω(u))$.

e) If $h:𝓡 \ito 𝓢$ and $𝓢={↓}u$ then $(j∘χ_h)(*∈𝓢(u)) = (𝓡^*∩𝓢∈Ω(u))$.

f) If $h:𝓡 \ito 𝓢$ and $𝓢={↓}u$ then $j(𝓡∈Ω(u)) = (𝓡^*∩𝓢∈Ω(u))$.

g) If $𝓢={↓}u$ then $j(𝓡∈Ω(u)) = (𝓡^*∩𝓢∈Ω(u))$.

h) $j(𝓡∈Ω(u)) = (𝓡^*∩{↓}u∈Ω(u))$.

i) $j(u)(𝓡) = 𝓡^*∩{↓}u$.

j) $j = λu∈𝐃.λ𝓡∈Ω(u).(𝓡^*∧{↓}u)$.

\msk

{\bf Proof.} Item (a) is a consequence of the formula in section
\ref{SetD-chi}:
$$\begin{array}{rcl}
   χ_f(b∈B(u)) &=& (\Cst(A  ∩ ({↓}(b∈B(u)))  ∈Ω(u)) \\[5pt]
  χ_f(*∈𝓢(u)) &=& (\Cst(𝓡 ∩ ({↓}(*∈𝓢(u))) ∈Ω(u)) \\
               &=& (\Cst(𝓡 ∩  {↓}u)         ∈Ω(u)) \\
               &=& (\Cst(𝓡)                 ∈Ω(u)) \\
               &=& (     𝓡                  ∈Ω(u)) \\
  \end{array}
$$

We can visualize it as:
%
%D diagram modality-to-j-item
%D 2Dx     100   +45
%D 2D  100 A0 -- A1
%D 2D
%D 2D  +30 A2 -- A3
%D 2D
%D 2D  +10 B0 -- B1
%D 2D
%D ren A0 A1 A2 A3 ==> 𝓡 1 𝓢 Ω
%D ren A0 A1 A2 A3 ==> \ga{R} 1 𝓢 Ω
%D ren B0 B1 ==> (*∈𝓢(u)) (\ga{R}∈Ω(u))
%D
%D (( A0 A1  -> .plabel= a !
%D    A0 A2 `-> .plabel= l \ga{f}
%D    A1 A3 `-> .plabel= r ⊤
%D    A2 A3  -> .plabel= a \ga{chif}
%D    newnode: A2' at: @A2+v(-12,0) .TeX= {↓}u= place
%D    A0 relplace 7 7 \pbsymbol{7}
%D    B0 B1 |->
%D ))
%D enddiagram
%D
\pu
\def\modalitytojitem#1#2#3{{
    \sa{R}{#1}
    \sa{f}{#2}
    \sa{chif}{#3}
    \diag{modality-to-j-item}
  }}
$$\modalitytojitem{𝓡}{f}{χ_f}
$$

\newpage

The items (b) and (c) are substitution instances of (a) and (b). These
substituions become easy to understand if we draw their diagrams:
$$%\text{b)}
  \modalitytojitem{𝓡^𝓢}{g}{χ_g}
  \quad
  %\text{c)}
  \modalitytojitem{𝓡^𝓢}{\ovl{h}}{χ_{\ovl{h}}}
$$

To draw a diagram for (d), (e), (f), (g), (h), (i), and (j) we
transform the arrow $\ovl{h}$ in the diagram for (c) into a arrow
going southwest in a diagram like the ones in sec.\ref{top-to-clop}:
%
%D diagram nuc-to-j-big-diagram
%D 2Dx     100  +20  +25   +60  +20
%D 2D  100 A0 ------ A1
%D 2D      |          |
%D 2D  +20 |    A2 ---|-------> A3
%D 2D      v  /       v       /
%D 2D  +20 A4 -----> A5 -> A6
%D 2D
%D 2D  +15 B0 -----------> B2
%D 2D  +16 C0 -----> C1
%D 2D  +8            D1 -> D2
%D 2D
%D ren A0 A1 ==> 𝓡 1
%D ren A2 A3 ==> 𝓡^𝓢 1
%D ren A4 A5 A6 ==> 𝓢 Ω Ω
%D ren B0    B2 ==> (*∈𝓢(u))           (𝓡^𝓢∈Ω(u))
%D ren C0 C1    ==> (*∈𝓢(u)) (𝓡∈Ω(u))
%D ren    D1 D2 ==>           (𝓡∈Ω(u)) (𝓡^*∩𝓢∈Ω(u))\ph{mm}
%D
%D (( A0 A1 ->
%D    A0 A4 >-> .PLABEL= _(0.30) h
%D    A1 A5 >->
%D    A2 A3 ->
%D    A2 A4 >-> .PLABEL= _(0.40) \ovl{h}
%D    A3 A6 >->
%D    A4 A5 -> .plabel= a χ_h
%D    A5 A6 -> .plabel= a j
%D    A4 A6 -> sl__ .plabel= b χ_{\ovl{h}}=j∘χ_h
%D    newnode: A4' at: @A4+v(-12,0) .TeX= {↓}u= place
%D
%D    B0 B2 |->
%D    newnode: B2' at: @B2+v(-9,8) .TeX= =(𝓡^*∩𝓢∈Ω(u)) place
%D    C0 C1 |->
%D    D1 D2 |->
%D ))
%D enddiagram
%D
$$\pu
  \diag{nuc-to-j-big-diagram}
$$

TODO: Write the proof of (e)$→$(f).

\newpage

%  ____       _   ____      _     _  _     
% / ___|  ___| |_|  _ \ _  | |__ (_)(_)___ 
% \___ \ / _ \ __| | | (_) | '_ \| || / __|
%  ___) |  __/ |_| |_| |_  | |_) | || \__ \
% |____/ \___|\__|____/(_) |_.__/|_|/ |___/
%                                 |__/     
%
% «SetD-bijections»  (to ".SetD-bijections")
% (cltp 34 "SetD-bijections")
% (clta    "SetD-bijections")

\subsection{Some bijections}
\label{SetD-bijections}

Let $𝐃$ be a finite 2-column graph. As usual, let $𝐄$ be the topos
$\SetD$, and let $H=\Subsets(1_𝐄)$ be its Heyting Algebra of
truth-values. Let's denote the set of subsets of $𝐃_0$ by $\Pts(𝐃)$,
the set of nuclei on $\SetD$ by $\Nucs(𝐃)$, the set of Grotendieck
topologies on $\SetD$ by $\GrTops(𝐃)$, the set of Lawvere-Tierney
topologies on $\SetD$ by $\LTTops(𝐃)$, and the set of closure
operators on $\SetD$ by $\Clops(𝐃)$.

% (linfp 1 "title")
% (linfa   "title")

We have bijections between the sets $\Pts(𝐃)$, $\Nucs(𝐃)$, and
$\GrTops(𝐃)$. They are proved, together with the commutatity of the
triangle, in \cite{Lindenhovius} and \cite{LindenhoviusOchs}. The
bijection between $\GrTops(𝐃)$ and $\LTTops(𝐃)$ is proved in
\cite{MacLaneMoerdijk}, theorem V.1.2 and section V.4, and the
bijection between $\LTTops(𝐃)$ and $\Clops(𝐃)$ has brief proofs in
several standard books and a detailed proof in sections
\ref{clops}--\ref{without-inclusions} here.

The diagram below shows all these bijections and where more
information about them can be found. The right half shows how we will
refer to the here; for example, $(J↔j)$ is the bijection between
$\GrTops(𝐃)$ and $\LTTops(𝐃)$, and we call its components $(J↦j)$ and
$(j↦J)$.
%
%D diagram bijections-references
%D 2Dx     100  +50 +25 +40  +50 +25
%D 2D  100 A0 - A1      B0 - B1
%D 2D  +25 |  /     A4  |  /     B4
%D 2D  +25 A2 - A3      B2 - B3
%D 2D
%D ren A0 A1 A2 A3 A4 ==> \Pts(𝐃) \Nucs(𝐃) \GrTops(𝐃) \LTTops(𝐃) \Clops(𝐃)
%D ren B0 B1 B2 B3 B4 ==> 𝓨 \nuc J j \clop
%D
%D (( A0 A1  -> sl^ .plabel= a \smt{C.4.2,\;C.2,}
%D    A0 A1 <-  sl_ .plabel= b \smt{\cite{LindenhoviusOchs}}
%D    A0 A2  -> sl_ .plabel= l \smtt{2.8,}{C.4.1}
%D    A0 A2 <-  sl^ .plabel= r \smt{2.9}
%D    A2 A1 <->     .plabel= m \smtt{B.8,}{B.25}
%D    A2 A3 <-> sl_ .plabel= b \smt{\cite{MacLaneMoerdijk}}
%D    A3 A4 <-> sl_ .plabel= r \smt{[here]}
%D ))
%D (( # B0 B1  |->   sl^ .plabel= a \sm{(𝓨↦\nuc)}
%D    # B0 B1 <-|    sl_ .plabel= b \sm{(\nuc↦𝓨)}
%D    B0 B1  <->       .plabel= m \sm{(𝓨↦\nuc),\\(\nuc↦𝓨)}
%D    B0 B2  <->       .plabel= m \sm{(𝓨↦J),\\(J↦𝓨)}
%D    B2 B1  <->       .plabel= m \sm{(\nuc↦J),\\(J↦\nuc)}
%D    B2 B3  <->       .plabel= m \sm{(J↦j),\\(j↦J)}
%D    B3 B4  <->       .plabel= r \sm{(j↦\clop),\\(\clop↦j)}
%D ))
%D enddiagram
%D
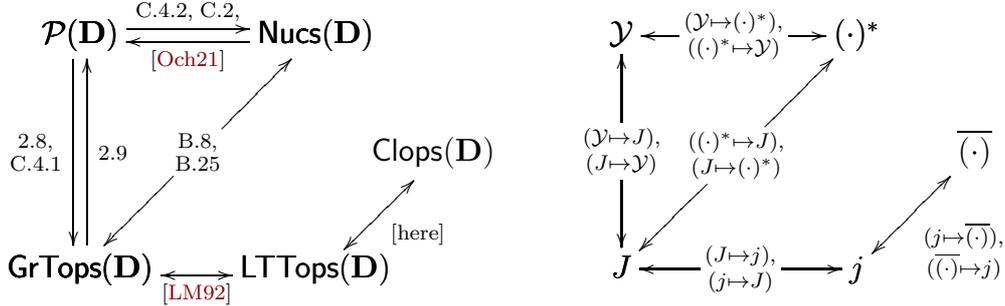
$$\pu
  \diag{bijections-references}
$$

These are the components of the bijections in the triangle. The
references like ``C.4.2'' point to sections, theorems, and definitions
in \cite{Lindenhovius}.
$$\begin{array}{rrcll}
  (𝓨↦\nuc): & 𝓢^* &=&  𝓨→𝓢                         & \text{C.4.2, C.2} \\
 %(\nuc↦𝓨): & 𝓨  &=& \setofst {u∈𝐃} {u\not∈\dnous}  & \text{\cite{LindenhoviusOchs}} \\
  (\nuc↦𝓨): & 𝓨  &=& \setofst {u∈𝐃} {\dnus\not=\dnous}  & \text{\cite{LindenhoviusOchs}} \\
  [5pt]
  (𝓨↦J):  & J(u)   &=& \setofst{𝓢∈Ω(u)}{𝓨∩\dnu⊆𝓢}  & \text{2.8, C.4.1} \\
  (J↦𝓨):  & X_J    &=& \setofst{u∈𝐃}{J(u)=\{{↓}u\}}  & \text{2.9} \\
  [5pt]
  (\nuc↦J): & J(u) &=& \setofst {𝓢∈Ω(u)} {u∈𝓢^*}     & \text{B.8, B.25} \\
  (J↦\nuc): & 𝓢^* &=& \setofst {u∈𝐃} {𝓢∩{↓}u∈J(u)}  & \text{B.8, B.25}  \\
  \end{array}
$$

These are the components of the bijection $(J↔j)$:
$$\begin{array}{rrcll}
  (J↦j): & j        &=& χ_{(J \ito Ω)}                  \\
  (j↦J): & J        &=& \dom(σ(j))                      \\
  \end{array}
$$

These are the components of the bijection $(j↔\clop)$:
$$\begin{array}{rrcll}
  (j↦\clop): & \ovl{(f:A \ito B)} &=& σ(j∘χ_f)  \\
  (\clop↦j): & j &=& χ_{\ovl{(1_⊤ \ito Ω)}}       \\
  \end{array}
$$

We have also defined these two operations in sections
\ref{SetD-nucleus} and \ref{SetD-nuc-to-j}:
$$\begin{array}{rrcll}
  (\clop↦\nuc): & 𝓢^* &=& \dom(\ovl{(𝓢 \ito 1)}) \\
  (\nuc↦j):     & j(u)(𝓡) &=& 𝓡^*∩{↓}u \\
  \end{array}
$$

We did not prove that our formulas for them yield exactly the
corresponding components of the bijections obtained by composing the
bijections that we know. Let state that as two conjectures:

\msk

\Conjecturesubsection
$(\nuc↦j) = (\nuc↦J);(J↦j)$.

\ssk

\Conjecturesubsection
$(\clop↦\nuc) = (\clop↦j);(j↦J);(J↦\nuc)$.

\newpage

%  ____       _   ____      _  _                     _           
% / ___|  ___| |_|  _ \ _  | || |        _   _ _ __ | | ___  ___ 
% \___ \ / _ \ __| | | (_) | || |_ _____| | | | '_ \| |/ _ \/ __|
%  ___) |  __/ |_| |_| |_  |__   _|_____| |_| | |_) | |  __/\__ \
% |____/ \___|\__|____/(_)    |_|        \__,_| .__/|_|\___||___/
%                                             |_|                
%
% «SetD-4-uples»  (to ".SetD-4-uples")
% (cltp 38 "SetD-4-uples")
% (clta    "SetD-4-uples")
\subsection{Valid 4-uples}
\label{SetD-4-uples}

\def\eqQ{∼_𝓠}
\def\eqnuc{∼_{\nuc}}
\def\PAQ{((𝐏_0,A),𝓠)}
\def\senwarrow{\rotatebox{45}{$↔$}}
\def\FourUple {\FourUpleMeta{𝓨}{\nuc}{J}{j}}
\def\FourUpleMeta#1#2#3#4{
  \mat{ #1 &      ↔     & #2 \\
        ↕  & \senwarrow &    \\
        #3 &      ↔     & #4 \\
      }
  }

One way to start to get some visual intuition on what nuclei and
topologies ``mean'' is to choose a 2-column graph $𝐃$ and then choose
either a $\nuc$, a $𝓨$, a $J$, or a $j$ on it, and then calculate the
other three elements of the 4-uple $(\nuc, 𝓨, J, j)$ from it using the
bijections from the previous section. We will always draw these
4-uples in this shape:
$$\FourUple$$

The main result of \cite{PH2} is a way to visualize the bijection
$(𝓨↔\nuc)$ when $𝐃$ is a 2-column graph. We will need to change its
notation a bit. In \cite{PH2} the notation for a 2CG with question
marks is $((P,A),Q)$ and the notation for a ZHA with a nucleus is
$(H,J)$; here we will use $\PAQ$ and $(H,\nuc)$. The set of arrows $A$
of the 2CG is a subset of $𝐏_0×𝐏_0$, and $𝓠$ is the set of points of
$𝐏_0$ ``with question marks''. Here is an example of a $\PAQ$ with its
associated $(H,\nuc)$:
%
%L tdims = TCGDims {qrh=5, q=15, crh=12, h=60, v=25, crv=7}   -- with v arrows
%L tspec_PAQ = TCGSpec.new("46; 11 22 34 45, 25", ".???", "???.?.")
%L tspec_PAQ:mp  ({zdef="O_A(P),J"}):addlrs():print()            :output()
%L tspec_PAQ:tcgq({tdef="(P,A),Q", meta="1pt p"}, "lr q h v ap") :output()
\pu

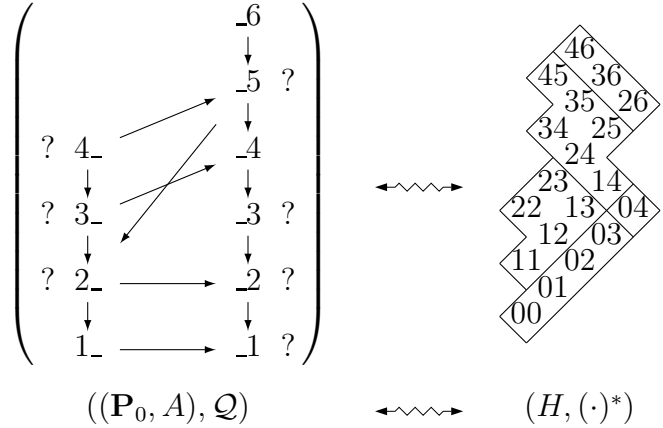
$$\begin{array}{ccc}
  \tcg{(P,A),Q} & \squigbij & \zha{O_A(P),J} \\
  [70pt]
  \PAQ & \squigbij & (H,\nuc) \\
  \end{array}
$$

\msk

We will regard two-digit numbers as ``piles'' of elements of $𝐏_0$:
$$ab \;\;=\;\; \pile(ab) \;\;=\;\; \pileelements{a}{b}$$

A set of question marks $𝓠$ induces an equivalence relation $\eqQ$ on
$H$: we have $ab \eqQ cd$ when $ab$ and $cd$ become indistinguishable
when we erase the information on the question marks; formally, $ab
\eqQ cd$ means $\pile(ab)∖𝓠 = \pile(cd)∖𝓠$. The operation $\nuc$ takes
each $ab∈H$ to the topmost element in its region --- in the example
above we have $12^*=23$ --- and it also induces an equivalence
relation: $ab \eqnuc cd$ when $ab^*=cd^*$. We say that a set of
question marks $𝓠$ ``is associated to'' a nucleus $\nuc$ when $(\eqQ)
= (\eqnuc)$.

\newpage

% (ph9p 31)
% (lindp 73 "C.2")
% (linda    "C.2")

Let's now expand the definitions of the two components of the
bijection $(𝓨↔\nuc)$. The ``Heyting Implication'' `$→$' defined in
\cite[def.C.2]{Lindenhovius} is exactly the intuitionistic implication
from \cite[sec.16]{PH1}, that is calculated using the topological
interior, and $\dnou$ is a shorthand for ${↓}u∖\{u\}$. It turns out
that $𝓨$ is exactly the set of points of $𝐏_0$ without questions
marks, i.e., $𝓠=𝐏_0∖𝓨$ and $𝓨=𝐏_0∖𝓠$:
$$\begin{array}{rrcll}
  (𝓨↦\nuc): & 𝓢^* &=&  𝓨→𝓢                           \\
             &      &=& \Int(\setofst{u∈𝐏}{u∈𝓨 → u∈𝓢}) \\
             &      &=& \Int(\setofst{u∈𝐏}{u\not∈𝓨 ∨ u∈𝓢}) \\
             &      &=& \Int(\setofst{u∈𝐏}{u∈(𝐏_0∖𝓨) ∨ u∈𝓢}) \\
             &      &=& \Int(\setofst{u∈𝐏}{u∈(𝐏_0∖𝓨)∪𝓢}) \\
             &      &=& \Int((𝐏_0∖𝓨)∪𝓢) \\
             &      &=& \Int(𝓠∪𝓢) \\
  (\nuc↦𝓨): & 𝓨  &=& \setofst {u∈𝐃} {\dnus\not=\dnous} \\
%            &     &=& \setofst {u∈𝐃} {\dnus\not∈\dnous} \\
  \end{array}
$$

We have:
$$\begin{array}{rrcll}
  (𝓠↦\nuc): & 𝓢^* &=& \Int(𝓠∪𝓢) \\
  (\nuc↦𝓠): & 𝓠  &=& \setofst {u∈𝐃} {\dnus=\dnous} \\
  \end{array}
$$

Let's see some examples to make this more concrete. Suppose that $𝓠$
is the set of question marks of the example above. Then we can
calculate $12^*$ by doing $\Int(𝓠∪\pile(12))$; the `$\Int$' discards
the points $4▁$, $3▁$, and $▁5$ from $𝓠∪\pile(12)$, and $12^* =
\Int(𝓠∪\pile(12)) = \pile(23) = 23$. Now let's start by the $\nuc$ of
the example and try to obtain $𝓠$. Let $u=3▁$. Then $\dnu=34$,
$\dnou=24$, and 34 and 24 are in the same region of the ZHA, so
$3▁∈𝓠$. Let's try $u=▁4$. Then $\dnu=04$, $\dnou=03$, and 04 and 03
are in different regions the ZHA, so $▁4\not∈𝓠$.

Using these methods we can easily convert a $𝓠$ to a $\nuc$ and
vice-versa; but we want to use $𝓨$ instead of $𝓠$ in our 4-uples, and
we will draw our `$𝓨$'s as subsets of $𝐏_0$, like this:
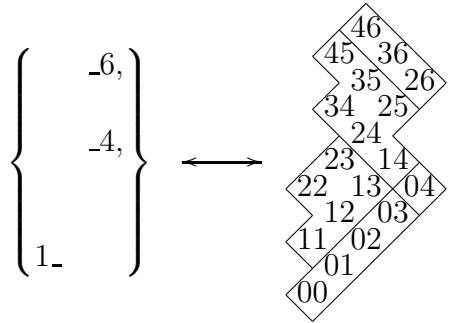
\def\RelevantPoints{
  \cmat{     & ▁6, \\
             &     \\
             & ▁4, \\
             &     \\
             &     \\
         1▁  &     \\
       }}

$$\RelevantPoints \;\; \diagxyto/<->/<250> \;\; \zha{O_A(P),J}$$

\newpage

In our section \ref{SetD-nuc-to-j} we saw that each component
$j(u):Ω(u)→Ω(u)$ of a Lawvere-Tierney topology $j$ is a modality on
the down-set $Ω(u)=\Downs({↓}u)$, and we have $j(u)(𝓡) = {↓}u∧𝓡^* =
𝓡^{{↓}u}$; so each $j(u)$ is ``$\nuc$ truncated to $Ω(u)$''. We will
use this idea to draw our `$j$'s in a nice way --- each component
$j(u)$ will be drawn as a slashing on the corresponding $Ω(u)$, and we
will drawing `$·$'s on all points of $H$ that are ``out of the
domain'', i.e., outside that $Ω(u)$. Remember that at the end of our
section \ref{SetD-Ju-filter} we drew $1_⊤ ⊂ J ⊂ Ω$ as:
$$\def\DA{\littlenbig{
           {\olittlen{{21}{·}  {·}{·}{·}    {·}{·}{·}}}
           {\olittlen{ {·}{·}  {·}{·}{·}   {·}{·}{02}}}
           {\olittlen{ {·}{·} {10}{·}{·}    {·}{·}{·}}}
           {\olittlen{ {·}{·}  {·}{·}{·}   {·}{01}{·}}}
         }}
  \def\DB{\littlenbig{
           {\olittlen{{21}{·} {10}{11}{·}   {·}{·}{·}}}
           {\olittlen{ {·}{·}  {·}{·}{·} {00}{01}{02}}}
           {\olittlen{ {·}{·} {10}{·}{·}    {·}{·}{·}}}
           {\olittlen{ {·}{·}  {·}{·}{·}  {00}{01}{·}}}
         }}
  \def\DC{\littlenbig{
           {\olittlen{{21}{·} {10}{11}{·} {00}{01}{·}}}
           {\olittlen{ {·}{·}  {·}{·}{·} {00}{01}{02}}}
           {\olittlen{ {·}{·} {10}{·}{·}   {00}{·}{·}}}
           {\olittlen{ {·}{·}  {·}{·}{·}  {00}{01}{·}}}
         }}
  \DA ⊂ \DB ⊂ \DC
$$

We will draw our `$j$'s as slashings on the components of our `$Ω$',
and of `$J$'s exactly as we drew above. Here is an example:
%
% From:
% (p2ap 27 "a-particular-case")
% (p2aa    "a-particular-case")
%
% «ArtDecoN»  (to ".ArtDecoN")
%
%L ArtDecoN_ts   = TCGSpec.new("33; 32,"):LRcolstrs("!ga{L1} !ga{L2} !ga{L3}",
%L                                                  "!ga{R1} !ga{R2} !ga{R3}")
%L ArtDecoN_td_0 = TCGDims {h=15,  v=8,  q=15, crh=3.5,  crv=7, qrh=5}
%L ArtDecoN_td_1 = TCGDims {h=25, v=22,  q=15, crh=7.5,  crv=7, qrh=5}
%L ArtDecoN_td_2 = TCGDims {h=65, v=50,  q=15, crh=20,  crv=15, qrh=5}
%L ArtDecoN_td_3 = TCGDims {h=85, v=70,  q=15, crh=30,  crv=30, qrh=5}
%L ArtDecoN_td_4 = TCGDims {h=85, v=80,  q=15, crh=35,  crv=35, qrh=5}
%L ArtDecoN_tq   = TCGQ.newdsoa(ArtDecoN_td_0, ArtDecoN_ts,
%L                              {tdef="ArtDecoNSmall", meta="1pt s"},
%L                              "h ap LR o")
%L ArtDecoN_tq   = TCGQ.newdsoa(ArtDecoN_td_1, ArtDecoN_ts,
%L                              {tdef="ArtDecoNMed", meta="1pt s"},
%L                              "h v ap LR o")
%L ArtDecoN_tq   = TCGQ.newdsoa(ArtDecoN_td_2, ArtDecoN_ts,
%L                              {tdef="ArtDecoNBig", meta="1pt"},
%L                              "h v ap LR o")
%L ArtDecoN_tq   = TCGQ.newdsoa(ArtDecoN_td_3, ArtDecoN_ts,
%L                              {tdef="ArtDecoNBigg", meta="1pt"},
%L                              "h v ap LR o")
%L ArtDecoN_tq   = TCGQ.newdsoa(ArtDecoN_td_4, ArtDecoN_ts,
%L                              {tdef="ArtDecoNBigg", meta="1pt"},
%L                              "h v ap LR o")
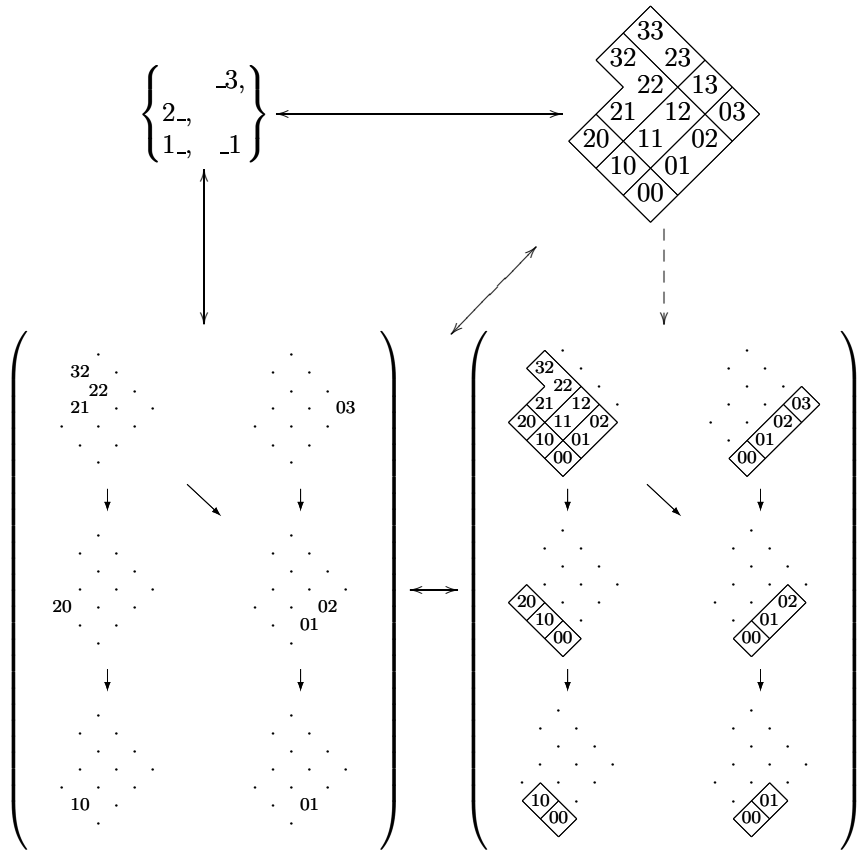
\pu
\def\ArtDecoNSetargs#1#2#3#4#5#6{
  \sa{L3}{#1}\sa{R3}{#2}%
  \sa{L2}{#3}\sa{R2}{#4}%
  \sa{L1}{#5}\sa{R1}{#6}%
  }
%
% «ArtDecoN-shortdefs»  (to ".ArtDecoN-shortdefs")
%
\def\adnsetargs#1{\ArtDecoNSetargs#1}
\def\adn       #1{{\adnsetargs#1        \tcg{ArtDecoNSmall}         }}
\def\padn      #1{{\adnsetargs#1 \left( \tcg{ArtDecoNSmall} \right) }}
\def\badn      #1{{\adnsetargs#1 \left[ \tcg{ArtDecoNSmall} \right] }}
\def\padnmed   #1{{\adnsetargs#1 \left( \tcg{ArtDecoNMed}   \right) }}
\def\padnbig   #1{{\adnsetargs#1 \left( \tcg{ArtDecoNBig}   \right) }}
\def\padnbigg  #1{{\adnsetargs#1 \left( \tcg{ArtDecoNBigg}  \right) }}
\def\padnbiggg #1{{\adnsetargs#1 \left( \tcg{ArtDecoNBigg}  \right) }}
\pu
\def\QMarks {\tcg{ArtDecoN-qmarks}}
\def\RelevantPoints{
  \cmat{     & ▁3, \\
         2▁, &     \\
         1▁, & ▁1  \\
       }}
\def\Nucleus{\zha{ArtDecoN-nucleus}}
\def\GrTopology{
  \padnbig{
    {{\badn{?·1?11}}}  {\badn{·1·?·1}}
     {\badn{··1·1·}}   {\badn{···?·1}}
     {\badn{····1·}}   {\badn{·····1}}
    }}
\def\GrTopology{
  \padnbiggg{
    {{\zha{OADN:J:3_}}} {\zha{OADN:J:_3}}
     {\zha{OADN:J:2_}}  {\zha{OADN:J:_2}}
     {\zha{OADN:J:1_}}  {\zha{OADN:J:_1}}
    }}
\def\LTTopology{
  \padnbiggg{
    {{\zha{OADN:j:3_}}} {\zha{OADN:j:_3}}
     {\zha{OADN:j:2_}}  {\zha{OADN:j:_2}}
     {\zha{OADN:j:1_}}  {\zha{OADN:j:_1}}
    }}
%
%D diagram bijections-particular-case
%D 2Dx     100  +75 +115
%D 2D  100 A00  A0  A1
%D 2D
%D 2D  +120     A2  A3
%D 2D
%D ren A00 A0 A1 ==> \QMarks  \RelevantPoints \Nucleus
%D ren     A2 A3 ==>          \GrTopology     \LTTopology
%D
%D (( A0 A1 <->
%D    A0 A2 <->
%D  # A1 A2 <->
%D    newnode: A1' at: tow(@A1,@A2,0.25)
%D    newnode: A2' at: tow(@A1,@A2,0.5)
%D    A1' A2' <->
%D    A1 A3 -->
%D    A2 A3 <->
%D ))
%D enddiagram
%D
$$\pu
  \scalebox{0.85}{$
  \diag{bijections-particular-case}
  $}
$$

\newpage

The diagram above is a particular case this one,
$$\FourUple$$
in the sense that we replaced each expression in its corners by its
value on a certain particular case. If we draw many such particular
cases --- by choosing a 2-column graph $𝐃$, then a nucleus $\nuc$ as a
slashing, at then calculating the corresponding `$𝓨$'s, `$j$', and
`$J$'s by hand or by computer, we will see that those same patterns
always occur: each element of $𝓨$ corresponds to one of the diagonal
cuts in the slashing of $\nuc$, $j$ is always obtained by doing
truncations of $\nuc$, and each $J(u)$ is made of the element in the
topmost equivalence class of the corresponding $j(u)$.

\bsk

(TODO: prove this last affirmation about `$J$'s and `$j$'s!)

%  ____       _   ____            _                 _ _     _             
% / ___|  ___| |_|  _ \ _  __   _(_)___ _   _  __ _| (_)___(_)_ __   __ _ 
% \___ \ / _ \ __| | | (_) \ \ / / / __| | | |/ _` | | |_  / | '_ \ / _` |
%  ___) |  __/ |_| |_| |_   \ V /| \__ \ |_| | (_| | | |/ /| | | | | (_| |
% |____/ \___|\__|____/(_)   \_/ |_|___/\__,_|\__,_|_|_/___|_|_| |_|\__, |
%                                                                   |___/ 
% «SetD-visualizing»  (to ".SetD-visualizing")
% (cltp 38 "SetD-visualizing")
% (clta    "SetD-visualizing")
\subsection{Visualizing nuclei and topologies}
\label{visualizing-nuclei-and-tops}

TODO: Explain how to use the diagrams for particular cases of the last
section to complement standard texts about Topos Theory; compare with
\cite[Section 5.5]{FavC}.

% (cltp 10 "top-to-clop")
% (clta    "top-to-clop")

\newpage

%  ____       _   ____        _    __   __        _ 
% / ___|  ___| |_|  _ \ _    (_)  / /   \ \      | |
% \___ \ / _ \ __| | | (_)   | | / /_____\ \  _  | |
%  ___) |  __/ |_| |_| |_    | | \ \_____/ / | |_| |
% |____/ \___|\__|____/(_)  _/ |  \_\   /_/   \___/ 
%                          |__/                     
% Deleted: «SetD-j-and-J»  (to ".SetD-j-and-J")
% (cltp 34 "SetD-j-and-J")
% (clta    "SetD-j-and-J")

% \subsection{Lawvere-Tierney subsumes Grothendieck: $j↔J$}
% \label{SetD-j-and-J}
%
% (find-books "__cats/__cats.el" "maclane-moerdijk")
% (find-maclanemoerdijkpage (+   11 219) "V.1 Lawvere-Tierney Topologies")
% (find-maclanemoerdijkpage (+   11 222)   "Theorem 2")
% (find-maclanemoerdijkpage (+   11 233) "V.4 Lawvere-Tierney Subsumes Grothendieck")

%L write_dnt_file()
\pu

\printbibliography

\GenericWarning{Success:}{Success!!!}  % Used by `M-x cv'

\end{document}